\documentclass[12pt,centertags]{amsart}
\usepackage{amsmath,amstext,amsthm,a4,amssymb,amscd}
\usepackage{newcent}
\usepackage[mathscr]{eucal}
\usepackage{mathrsfs}
\numberwithin{equation}{section}

\newcommand{\field}[1]{\mathbb{#1}}
\newcommand{\bZ}{\field{Z}}
\newcommand{\bR}{\field{R}}
\newcommand{\bC}{\field{C}}
\newcommand{\bN}{\field{N}}
\newcommand{\bD}{\field{D}}
\newcommand{\proj}{\field{P}}
\def\bk{{\mathbf k}}

\def\br{{\mathbf r}}

\newcommand{\cali}[1]{\mathscr{#1}}
\newcommand{\cC}{\cali{C}} 
\newcommand{\cO}{\cali{O}} 
\newcommand{\cH}{\cali{H}} 
 \newcommand{\cL}{\cali{L}}
\newcommand{\cR}{\cali{R}} 

\newcommand{\calig}[1]{\mathcal{#1}}
\newcommand\mA{\calig{A}} \newcommand\mB{\calig{B}}
\newcommand\mC{\calig{C}} \newcommand\mH{\calig{H}}
\newcommand\mL{\calig{L}} \newcommand\mO{\calig{O}}
\newcommand\mQ{\calig{Q}} \newcommand\mR{\calig{R}}
\newcommand\mS{\calig{S}}
\def\mP{{\mathcal{P}}}\def\mJ{{\mathcal{J}}}

\newcommand{\til}[1]{\widetilde{#1}}
\newcommand{\tx}{\til{X}} \newcommand{\tj}{\til{J}} 
\newcommand{\tbj}{\til{\mathbf J}}
  
 \newcommand{\tom}{\til{\omega}} 
\newcommand{\tphi}{\til{\Phi}}
\newcommand{\tl}{\til{L}}  \newcommand{\te}{\til{E}}
  \newcommand{\tdel}{\til{\Delta}_{p,\tphi}}

\DeclareMathOperator{\End}{End} 
\DeclareMathOperator{\Ker}{Ker} \DeclareMathOperator{\Dom}{Dom}
 
  \DeclareMathOperator{\spec}{Spec}
\DeclareMathOperator{\Id}{Id} \DeclareMathOperator{\supp}{supp}
\DeclareMathOperator{\tr}{Tr} 
\DeclareMathOperator{\td}{Td} \DeclareMathOperator{\vol}{vol}
\DeclareMathOperator{\ch}{ch} 
\DeclareMathOperator{\db}{\overline\partial}

\newcommand{\norm}[1]{\lVert#1\rVert} \newcommand{\abs}[1]{\lvert#1\rvert}
\newcommand{\om}{\omega} \newcommand{\g}{\Gamma}

\newtheorem{thm}{Theorem}[section]
\newtheorem{lemma}[thm]{Lemma}
\newtheorem{prop}[thm]{Proposition}
\newtheorem{cor}[thm]{Corollary}
\theoremstyle{definition}
\newtheorem{defn}[thm]{Definition}
\theoremstyle{remark}
\newtheorem{rem}[thm]{Remark}
\newcommand{\be}{\begin{equation}}
\newcommand{\ee}{\end{equation}}
\newcommand{\wi}{\widetilde}
\newcommand{\var}{\varepsilon}
\newcommand{\ov}{\overline}
\newcommand{\comment}[1]{}

\begin{document}
\title{Generalized Bergman kernels on symplectic manifolds}
\date{\today}
\author{Xiaonan Ma}
\address{ Centre de Math\'ematiques Laurent Schwartz,
UMR 7640 du CNRS, \'Ecole Polytechnique, 91128 Palaiseau Cedex,
France}
\email{ma@math.polytechnique.fr}

\author{George Marinescu}
\address{Humboldt-Universit{\"a}t zu Berlin, Institut f{\"u}r Mathematik, 
Rudower Chaussee 25, 12489 Berlin, Germany}
\email{george@mathematik.hu-berlin.de}


\begin{abstract}
We study the near diagonal asymptotic expansion of the generalized Bergman kernel
of the renormalized Bochner-Laplacian
on high tensor powers of a positive line bundle 
over a compact symplectic manifold. We show how to compute the coefficients
of the expansion by recurrence and give a closed formula for the
 first two of them.
As consequence, we calculate the density of states function of the 
Bochner-Laplacian and establish a symplectic version 
of the convergence of the induced Fubini-Study metric. 
We also discuss generalizations of the asymptotic expansion for
 non-compact or singular manifolds as well as their applications.
Our approach is inspired by the analytic localization 
techniques of Bismut-Lebeau.
\end{abstract}

\maketitle
\setcounter{section}{-1}
\section{Introduction} \label{s1}

The Bergman kernel for complex projective
manifolds is the smooth kernel of the orthogonal projection 
from the space of smooth sections of 
a positive line bundle $L$ on the space of holomorphic sections of $L$,
or, equivalently, on the kernel of the Kodaira-Laplacian 
$\Box^L= \overline{\partial} ^L\overline{\partial} ^{L*}
+ \overline{\partial} ^{L*}\overline{\partial} ^L$ on $L$.
It is studied in \cite{Tian,Ru,Zelditch,Catlin,BSZ,Lu,Wang1, KS01}, 
 in various generalities,
 establishing the  diagonal asymptotic expansion 
for high powers of $L$. Moreover, the
coefficients in the diagonal asymptotic expansion encode geometric
information
 about the underlying complex projective manifolds. 
This  diagonal asymptotic expansion
plays a crucial role in the recent work of Donaldson \cite{D} where the
existence of K\"ahler metrics with constant scalar curvature is
shown to be closely related to Chow-Mumford stability.

In \cite{DLM}, Dai, Liu and Ma
 studied the asymptotic expansion of the Bergman kernel
of the spin$^c$ Dirac operator associated to a positive line bundle
on a compact symplectic manifold, 
and related it to that of the corresponding heat kernel. As a by
product, they gave a new proof of the above results.
This approach is inspired by Local Index Theory,
especially by the analytic localization 
techniques of Bismut-Lebeau \cite[\S 11]{BL}.

Another natural generalization of the operator $\Box^L$ 
in symplectic geometry was initiated by Guillemin and Uribe 
\cite{GU}. In this very interesting short paper,
they introduce a renormalized Bochner--Laplacian (cf.\,\eqref{laplace}) 
which is exactly $2\Box ^L$ in the K\"ahler case.
 The asymptotic of the 
spectrum of the renormalized Bochner--Laplacian on $L^p$ when $p\to \infty$
is studied in various generalities in 
\cite{BU,Bra1,GU} by applying
 the analysis of Toeplitz structures
of Boutet de Monvel--Guillemin \cite{BoG}, 
and in \cite{MM} as a direct application of Lichnerowicz formula.

Of course, there exists also a replacement of the $\db$--operator and of the  
notion of holomorphic section based on a construction of Boutet de  
Monvel--Guillemin \cite{BoG} of a first order pseudodifferential operator  
$D_b$ which mimic the $\db_b$ operator on the circle bundle associated  
to $L$. However, $D_b$ is neither canonically defined nor unique. This  
point of view was adopted in a series of papers 
\cite{BU1,SZ02,BSZ1}.

In this paper, we will study the asymptotic expansion of the 
generalized Bergman kernel of the renormalized Bochner-Laplacian, 
namely the smooth kernel of  the projection on its bound states 
as $p\to \infty$. The advantage of this approach is that the renormalized  
Bochner-Laplacian has geometric meaning and is canonically defined.
Moreover, it does not require the passage to the associated circle bundle 
as we can work directly on the base manifold.   
Let's explain our results in detail.

Let $(X,\om)$ be a compact symplectic manifold of real dimension $2n$.
Assume that there exists a Hermitian line bundle $L$ over $X$ endowed with
a Hermitian connection $\nabla^L$ with the property that
$\frac{\sqrt{-1}}{2\pi}R^L=\omega$, where $R^L=(\nabla^L)^2$ 
is the curvature of $(L,\nabla^L)$. Let $(E,h^E)$ be a Hermitian vector
bundle  on $X$ with Hermitian connection $\nabla^E$ and curvature $R^E$.
Let $g^{TX}$ be a Riemannian metric on $X$ 
and ${\bf J}:TX\longrightarrow TX$ be the skew--adjoint 
linear map which satisfies the relation
\be \label{0.1}
\om(u,v)=g^{TX}({\bf J}u,v)\quad\text{for}\quad u,v \in TX.
\ee
Let $J$ be an almost complex structure
such that $g^{TX}(J\cdot, J\cdot)=g^{TX}(\cdot, \cdot)$, 
$\om(J\cdot, J\cdot)=\om(\cdot, \cdot)$ 
and that $\om(\cdot, J\cdot)$ defines a metric on $TX$.
Then $J$ commutes with ${\bf J}$, and $-J{\bf J}\in \End(TX)$ is positive,
thus $J= {\bf J} (-{\bf J}^2)^{-1/2}$.

We introduce the Levi-Civita connection $\nabla^{TX}$ on $(TX, g^{TX})$  
with its curvature $R^{TX}$ and scalar curvature $r^X$. 
Let $\nabla ^XJ \in T^*X \otimes \End(TX)$ be the covariant 
derivative of $J$ induced by $\nabla ^{TX}$.  Let $\Delta ^{L^p\otimes E}$ 
be the induced Bochner-Laplacian  acting on $\cC^\infty (X,L^p\otimes E)$.
We fix a smooth Hermitian section $\Phi$ of $\End (E)$ on $X$.
Let $\{e_i\}_i$ be an orthonormal frame of $(TX, g^{TX})$. Set
\begin{align}\label{0.2}
&\tau(x)=-\pi \tr_{|TX} [J{\bf J}]= \frac{\sqrt{-1}}{2} R^L (e_j, J e_j) >0,\\
&\mu_0=\displaystyle\inf_{{u\in T_x X,\,x\in X}}
\sqrt{-1} R^L_x(u, Ju)/|u|^2_{g^{TX}} >0 , \label{0.21}\\
& \Delta_{p, \Phi} = \Delta ^{L^p\otimes E} -p\tau  + \Phi.\label{laplace}
\end{align}
By \cite[Cor. 1.2]{MM} (also cf. \cite{BU,GU,Bra1,BiV}),
there exists $C_L>0$ (which can be estimated precisely by using the $\cC^0$-norms of 
$R^{TX}$, $R^E$, $R^L$, $\nabla ^XJ$ and $\Phi$ cf. \cite[p.656\,--\,658]{MM}) 
independent of $p$ such that
\begin{equation}\label{0.3}
\spec\,\Delta_{p, \Phi} \subset [-C_L, C_L] \cup [2p\mu_0 -C_L, + \infty[\,.
\end{equation}
where we denote by $\spec(A)$ the spectrum of any operator $A$.

Let $\mH_p$ be the eigenspace of $\Delta_{p,\Phi}$ with the eigenvalues in
$[-C_L, C_L]$.
 Then for $p$ large enough,
again by \cite[Cor. 1.2]{MM} (also cf. \cite{BU,GU}  when $E$ is
trivial and ${\bf J}=J$)
\begin{align}\label{0.0}
&\dim \mH_p =d_p = \int_X \ch(L^p\otimes E)\td(TX)\\
&\hspace*{3mm} 
={\rm rk} (E) \int_X \frac{c_1(L)^n}{n!} p^n 
+ \int_X \Big(c_1(E) + \frac{{\rm rk} (E)}{2} c_1(TX)\Big)
\frac{c_1(L)^{n-1}}{(n-1)!} p^{n-1} + \cO(p^{n-2}), \nonumber
\end{align}
where $\ch(\cdot), c_1(\cdot), \td(\cdot)$ are the Chern character,
the first Chern class and the Todd class of the corresponding complex 
vector bundles ($TX$ is a complex vector bundle with complex structure $J$). 

Let $\{S^p_i\}_{i=1}^{d_p}$  be any orthonormal basis of
$\mH_p$ with respect to the inner
product \eqref{0c2} such that $\Delta_{p, \Phi} S^p_i = \lambda_{i,p} S^p_i$.
For $q\in \bN$, we define $B_{q,p}\in \cC ^\infty (X, \End (E))$ as follows,
\begin{align} \label{0.4}
B_{q,p}(x) &= \sum_{i=1}^{d_p} \lambda_{i,p}^q S^p_i (x) \otimes (S^p_i(x))^*,
\end{align}
here we denote by $\lambda_{i,p}^0=1$. 
Clearly, $B_{q,p}(x)$ does not depend on the choice of $\{S^p_i\}$.
Let $\det {\bf J}$ be the determinant function of ${\bf J}_x\in \End(T_xX)$.
A corollary of Theorem \ref{t3.8} is one of our main results:
\begin{thm}\label{t0.1} There exist smooth coefficients
$b_{q,r}(x)\in \End (E)_x$ which  are polynomials in $R^{TX}$,
$R^E$ {\rm(}and $R^L$, $\Phi${\rm)}, their derivatives of order
$\leqslant 2(r+q)-1$ {\rm(}resp. $2(r+q)${\rm)}, and reciprocals
of linear combinations of eigenvalues of ${\bf J}$ at $x$, with 
\begin{equation}\label{0.5}
b_{0,0}=(\det {\bf J})^{1/2}  \Id_E,
\end{equation}
such that for any $k,l\in \bN$, there exists
$C_{k,\,l}>0$ such that for any $x\in X$, $p\in \bN$,
\begin{align}\label{0.6}
&\Big |\frac{1}{p^n}B_{q,p}(x)
- \sum_{r=0}^{k} b_{q,r}(x) p^{-r} \Big |_{\cC ^l} \leqslant C_{k,\,l}\: p^{-k-1}.
\end{align}
Moreover, the expansion is uniform in the following sense{\rm:} for any fixed $k,l\in \bN$, 
assume that the derivatives of $g^{TX}$, $h^L$, $\nabla ^L$, $h^E$, $\nabla ^E$, $J$ and $\Phi$
with order $\leqslant 2n+2k+2q+l+4$ run over a set bounded in the $\cC^l$--\,norm
taken with respect to the parameter $x\in X$ and, moreover, $g^{TX}$ runs over a set
bounded below. Then the constant $C_{k,\,l}$ is independent of $g^{TX}${\rm;}
and the $\cC ^l$-norm in \eqref{0.6} includes also the  
derivatives on the parameters. 
\end{thm}

\noindent
By derivatives with respect to the parameters we mean directional derivatives in the spaces
of all appropriate $g^{TX}$, $h^L$, $\nabla ^L$, $h^E$, $\nabla ^E$, $J$ and $\Phi$ (on which
$B_{q,p}$ and $b_{q,r}$ implicitly depend).

We calculate further
the coefficients $b_{0,1}$ and $b_{q,0}$\,, $q\geqslant1$ as follows\footnote{Here 
$|\nabla ^X J |^2= \sum_{i j}|(\nabla_{e_i} ^X J)e_j |^2$
which is two times the corresponding $|\nabla ^X J |^2$ from \cite{MM04a}.}.  
\begin{thm}\label{t0.2} If $J={\bf J}$, then for $q\geqslant1$,
\begin{align}\label{0.8}
&b_{0,1}= \frac{1}{8\pi}\Big[r^X + \frac{1}{4} |\nabla ^X J |^2
+ 2\sqrt{-1} R^E (e_j,Je_j)\Big],\\
&b_{q,0}=\Big (\frac{1}{24}   |\nabla ^{X}J| ^2 
+\frac{\sqrt{-1}}{2} R^E (e_j,Je_j) + \Phi\Big )^q.\label{0.7}
\end{align}
\end{thm}
Let us check our formulas with the help of 
the Atiyah-Singer formula \eqref{0.0}.
Let $T^{(1,0)}X = \{ v\in TX\otimes_\bR \bC; Jv=\sqrt{-1}v\}$ be the almost 
complex tangent bundle on $X$ and let $P^{1,0}= \frac{1}{2} (1-\sqrt{-1}J)$ 
be the natural projection from $TX\otimes_\bR \bC$ onto $T^{(1,0)}X$. 
Then $\nabla ^{1,0}= P^{1,0}\nabla ^{TX} P^{1,0}$ is a
Hermitian connection on $T^{(1,0)}X$,  and 
the Chern-Weil representative of  $c_1(TX)$ is $c_1(T^{(1,0)}X, \nabla ^{1,0})
= \frac{\sqrt{-1}}{2 \pi} \tr_{|T^{(1,0)}X} (\nabla^{1,0})^2$.
By (\ref{g11}), 
\begin{align}\label{0.10}
(\nabla^{1,0})^2 = P^{1,0} \Big[R^{TX}
-\frac{1}{4} (\nabla ^XJ)\wedge (\nabla ^XJ)\Big]P^{1,0}.
\end{align}
Thus if $J={\bf J}$,
then by (\ref{0.10}), (\ref{g29}), (\ref{g30}), (\ref{g36}) and (\ref{g38}), 
\begin{align}\label{0.11}
\left\langle c_1(T^{(1,0)}X, \nabla ^{1,0}), \omega\right\rangle
&= \frac{1}{4\pi} \Big(r^X + \frac{1}{4} |\nabla ^X J |^2\Big).
\end{align}
Therefore, by integrating over $X$ the expansion \eqref{0.6} for $k=1$ we obtain \eqref{0.0}, so
\eqref{0.8} is compatible with \eqref{0.0}.

Theorem \ref{t0.1} for $q=0$ and \eqref{0.8} generalize the results of 
\cite{Catlin,Zelditch,Lu} and \cite{Wang1}
to the symplectic case. 
The term $r^X + \frac{1}{4} |\nabla ^X J |^2$ 
in \eqref{0.8} is called the Hermitian scalar curvature in the literature
\cite[Chap.\,10]{Ga04} 
and is a natural substitute for the 
Riemannian scalar curvature in the almost-K\"ahler case. 
It was used by Donaldson \cite{D99} to define the moment map on the space of 
compatible almost-complex structures.
We can view  \eqref{0.7}  as an extension and refinement of 
the results of \cite{BU2}, \cite[\S 5]{GU} about the density of 
states function of $\Delta_{p, \Phi}$ (cf. Remark \ref{t5.2} for
the details).

In \cite{DLM}, Dai, Liu and Ma also focused on the {\em full off-diagonal}  
asymptotic expansion 
(cf. \cite[Theorem 4.18]{DLM}) which is needed to study the Bergman kernel
on orbifolds, and the only small eigenvalue of the operator is $0$
when $p\to\infty$, thus they had the key equation \cite[(4.89)]{DLM}. 
In the current situation, we have small eigenvalues (cf. \eqref{0.3}) 
and we are interested to prove Theorem \ref{t3.8}, that is, the 
{\em near diagonal} expansion of the generalized Bergman kernels. 
This result is enough for most of applications.  
At first, the spectral gap \eqref{0.3} and the finite propagation speed of 
solutions of hyperbolic equations allow to localize the problem. 
Then we will combine the Sobolev norm estimates as in \cite{DLM}
and a formal power series trick to obtain Theorem \ref{t3.8},
and in this way, we get a method to compute the coefficients 
(cf. \eqref{c90}, \eqref{1c54}) which is new also in the case of \cite{DLM}.

In a forthcoming paper \cite{MM04c}, we will find the {\em full off-diagonal}
asymptotic expansion of the  generalized Bergman kernels by combining 
the results here and in \cite{DLM}, and as a direct  application, 
we will study the Toeplitz operators on symplectic manifolds and 
Donaldson Theorem \cite{D96} for the Kodaira map $\Phi_{p}$  \eqref{sz0}.

Let us provide a short road-map of the paper.
  In  Section \ref{s3}, we  prove  Theorem \ref{t0.1}.
  In  Section \ref{s4}, we compute the coefficients $b_{q,r}$,
and thus establish Theorem \ref{t0.2}.
  In Section \ref{s5}, we explain some applications of our results.
Among others, we give a symplectic version 
of the convergence of the induced Fubini-Study metric \cite{Tian}, and
we show how to handle the first-order pseudo-differential operator $D_b$
of Boutet de Monvel and Guillemin \cite{BoG}, which was studied 
extensively by Shiffman and  Zelditch \cite{SZ02}, and the operator 
$\ov{\partial}+\ov{\partial}^*$ when $X$ is K\"ahler but ${\bf J}\neq J$.
We include also generalizations for non-compact or singular manifolds
and as a consequence we obtain an unified treatment 
of the convergence of the induced Fubini--Study metric, 
the holomorphic Morse inequalities
and the characterization of Moishezon spaces.
Some results of this paper have been announced in \cite{MM04a}.
We refer also the readers our recent book \cite{MM05b} for our approach.


\section{Generalized Bergman kernels }\label{s3}

As pointed out in Introduction, we will apply the strategy of the proof
in \cite{DLM}. However, we have small eigenvalues when $p\to \infty$ 
(cf. (\ref{0.3})), thus we cannot use directly the key equation 
\cite[(4.89)]{DLM} to get a {\em full off-diagonal} asymptotic expansion
of the generalized Bergman kernels. 
After localizing the problem, we will adapt the Sobolev norm estimates 
developed in \cite{DLM} to our problem in Section \ref{s3.3}. 
To complete the proof of Theorem \ref{t0.1},  
we need to prove the vanishing of the coefficients $F_{q,r}$ ($r<2q$) in the expansion
\eqref{0ue45}. 
We will introduce a formal power series trick to overcome this difficulty 
and give a method to compute the coefficients in (\ref{0.6}). 
The ideas used here are inspired by the technique of Local Index Theory,
especially by \cite[\S 10, 11]{BL}.

This Section is organized as follows. In Section \ref{s3.1}, we
explain that the asymptotic expansion of the generalized Bergman kernel 
$P_{q,p}(x,x')$ is local on $X$ by using the spectral gap (\ref{0.3}) 
and the finite propagation speed of solutions of hyperbolic equations. 
In Section \ref{s3.2}, we obtain an asymptotic expansion of
$\Delta_{p,\Phi}$ in normal coordinates. In Section \ref{s3.3}, we
study the uniform estimate of the generalized Bergman kernels  of the 
renormalized Bochner-Laplacian $\cL_t$. In Section \ref{s3.4}, we
study the Bergman kernel of the limit operator $\cL_0$. In Section
\ref{s3.5}, we compute some coefficients $F_{q,r}(r\leqslant 2q)$ in
the asymptotic expansion in Theorem \ref{tue14}. Finally, in
Section \ref{s3.6}, we prove Theorem \ref{t0.1}.


\subsection{Localization of the problem}\label{s3.1}
Let $a^X$ be the injectivity radius of $(X, g^{TX})$.
We fix  $\var\in (0,a^X/4)$.
We denote by $B^{X}(x,\varepsilon)$ and  $B^{T_xX}(0,\varepsilon)$ the open balls
in $X$ and $T_x X$ with center $x$ and radius $\varepsilon$, respectively.
Then the map $ T_x X\ni Z \to \exp^X_x(Z)\in X$ is a
diffeomorphism from $B^{T_xX} (0,\varepsilon)$  on $B^{X} (x,\varepsilon)$ for
$\varepsilon\leqslant a^X$.  From now on, we identify $B^{T_xX}(0,\varepsilon)$
with $B^{X}(x,\varepsilon)$ for $\varepsilon \leqslant a^X$.

Let $\langle\,,\,\rangle_{L^p\otimes E}$ be the metric on 
$L^p\otimes E$ induced by  $h^L$ and $h^E$ 
and $dv_X$ be the Riemannian volume form of $(TX, g^{TX})$.
The $L^2$--scalar product on $\cC^\infty (X,L^p\otimes E)$, 
the space of smooth sections of $L^p\otimes E$, is given by
\begin{equation}\label{0c2}
\langle s_1,s_2 \rangle =\int_X\langle s_1(x),
s_2(x)\rangle_{L^p\otimes E}\,dv_X(x)\,.
\end{equation}
We denote the corresponding norm with $\norm{\,\cdot\,}_{L^2}$.

Let $\nabla^{TX}$ be the Levi--Civita connection
of the metric $g^{TX}$ and $\nabla^{L^p\otimes E}$ be the  
connection on $L^p\otimes E$ induced by $\nabla^L$ and  $\nabla^E$.
Let $\{e_i\}_i$ be an orthonormal frame  of $TX$.
Then the Bochner-Laplacian on $L^p\otimes E$ is given by
\be\label{0c1}
\Delta ^{L^p\otimes E} =- \sum_i \Big [ (\nabla ^{L^p\otimes E}_{e_i} )^2 -
\nabla ^{L^p\otimes E}_{\nabla ^{TX}_{e_i}{e_i}}\Big ].
\ee
Let $P_{\mH_p}$ be the orthonormal projection from
$\cC ^\infty (X,L^p\otimes E)$ onto $\mH_p$, the span of eigensections of
$\Delta_{p,\Phi}= \Delta ^{L^p\otimes E} -p \tau +\Phi$ corresponding to  
eigenvalues in $[-C_L, C_L]$.  
 
\begin{defn} \label{d3.0}
The smooth kernel of $(\Delta_{p,\Phi})^qP_{\mH_p}\,,\,q\geqslant 0$ (where
$(\Delta_{p,\Phi})^0=1$), with respect to $dv_X(x')$ is  
denoted $P_{q,p}(x,x')$ and is called a {\em generalized Bergman kernel\/}  
of $\Delta_{p,\Phi}$\,.
\end{defn} 
The kernel $P_{q,p}(x,x')$ is a section of 
$\pi^*_1(L^p\otimes E)\otimes\pi^*_2(L^p\otimes E)^*$ over $X\times X$,  
where $\pi_1$ and $\pi_2$ are the projections of $X\times X$ on the first  
and second factor. Using the notations of \eqref{0.4} we can write  
\begin{equation}\label{0.42} 
P_{q,p}(x,x')=\sum^{d_p}_{i=1}\lambda^q_{i,p}S^p_i(x)\otimes(S^p_i(x'))^* 
\in(L^p\otimes E)_x\otimes(L^p\otimes E)^*_{x'} .
\end{equation} 
 Since $L^p_x\otimes(L^p_x)^*$ is canonically isomorphic to $\bC$, the  
restriction of $P_{q,p}$ to the diagonal $\{(x,x):x\in X\}$ can be  
identified to $B_{q,p}\in\cC^\infty(X,E\otimes E^*)=\cC^\infty(X,\End(E))$.  

Let $f : \bR \to [0,1]$ be a smooth even function such that
$f(v)=1$ for $|v| \leqslant  \var/2$, and $f(v) = 0$ for $|v| \geqslant \var$.
Set
 \begin{align} \label{0c3}
F(a)= \Big(\int_{-\infty}^{+\infty}f(v) dv\Big)^{-1} 
\int_{-\infty}^{+\infty} e ^{i v a} f(v) dv.
\end{align}
Then $F(a)$ is an even function and
 lies in the Schwartz space $\mathcal{S} (\bR)$ and $F(0)=1$.
Let $\wi{F}$ be the holomorphic function on
 $\bC$ such that $\wi{F}(a ^2) =F(a)$.
The restriction of $\wi{F}$
 to $\bR$ lies in the Schwartz space $\mS (\bR)$. Then there exists
 $\{c_j\}_{j=1}^{\infty}$ such that for any $k\in \bN$, the function
\begin{align} \label{0c5}
F_k(a)= \wi{F}(a) - \sum_{j=1}^k c_j a ^j \wi{F}(a),
 \end{align}
verifies
\begin{align} \label{0c6}
F_k^{(i)}(0)= 0 \quad \mbox{\rm for any} \  0< i\leqslant k.
 \end{align}

\begin{prop}\label{0t3.0}
 For any $k,m\in \bN$, there exists $C_{k,m}>0$ such that for $p\geqslant1$
\begin{align}\label{0c7}
\left|F_k\big(\tfrac{1}{\sqrt{p}}\Delta_{p,\Phi}\big)(x,x') 
- P_{0,p}(x,x')\right|_{\cC ^m(X\times X)}
\leqslant C_{k,m} p^{-\frac{k}{2} +2(2m+2n+1)}.
\end{align}
Here the $\cC ^m$ norm is induced by $\nabla^L$ and $\nabla^E$.
\end{prop}
\begin{proof}
By \eqref{0c3}, for any $m\in \bN$, there exists $C'_{k,m}>0$ such that
\be\label{c9}
\sup_{a\in \bR} |a|^m |F_k(a) | \leqslant C'_{k,m}.
\ee
Set
\begin{align} \label{0c8}
G_{k,p}(a)= 1_{[\sqrt{p}\mu_0, +\infty[} (a) F_k(a),
\quad H_{k,p}(a)= 1_{[0, \frac{C_L}{\sqrt{p}}]} (|a|) F_k(a),
\end{align}
By (\ref{0.3}), for $p$ big enough,
\begin{equation} \label{0c9}
F_k\big(\tfrac{1}{\sqrt{p}}\Delta_{p,\Phi}\big)
= G_{k,p}\big(\tfrac{1}{\sqrt{p}}\Delta_{p,\Phi}\big)
+  H_{k,p}\big(\tfrac{1}{\sqrt{p}}\Delta_{p,\Phi}\big).
\end{equation}

As $X$ is compact, there exist $\{x_i\}_{i=1}^r$ such that
 $\{U_i = B^X(x_i,\var)\}_{i=1}^r$ is a covering of $X$.
We identify $B^{T_{x_i}X}(0,\var)$ with $B^{X}(x_i,\var)$
by the exponential map as above. We identify
$ (L^p\otimes E)_Z$ for $Z\in B^{T_{x_i}X}(0,\var)$
to $ (L^p\otimes E)_{x_i}$  by parallel transport with respect
to the connection $\nabla^{L^p\otimes E}$ along the curve
$\gamma_Z: [0,1]\ni u \to \exp^X_{x_i} (uZ)$. Let $\{e_i\}_i$
be an orthonormal basis of $T_{x_i}X$. Let $\wi{e}_i (Z)$ be the parallel
transport of ${e}_i$ with respect to $\nabla^{TX}$ along the above curve.
Let $\Gamma ^E, \Gamma ^L$
be the corresponding connection forms of $\nabla^E$, $\nabla^L $
 with respect to any fixed frame for
$E,L$ which is parallel along the curve $\gamma_Z$
under the trivialization on $U_i$.
Denote by $\nabla_U$ is the ordinary differentiation
 operator on $T_{x_i}X$ in the direction $U$. Then
\be\label{c10}
\nabla ^{L^p\otimes E}_{e_j}
= \nabla _{e_j}+ p \Gamma ^L(e_j) + \Gamma ^E(e_j).
\ee
Let $\varphi_i$ be a partition function associated to $\{U_i\}$. We define
 an Sobolev norm on the $l$-th Sobolev space $H^l(X,L^p\otimes E)$ by
\be\label{c11}
\| s\| _{H^l_p}^2 = \sum_i \sum_{k=0}^l \sum_{i_1 \cdots i_k=1} ^{2n}
\|\nabla_{e_{i_1}}\cdots  \nabla_{e_{i_k}}(\varphi _i s)\|_{L^2}^2
\ee
Then by \eqref{0.2}, \eqref{0c1}, \eqref{c10},
there exists $C>0$ such that for $p\geqslant1$, $s\in H^2(X, L^p\otimes E)$,
\be\label{c12}
\|s\|_{H^2_p}   \leqslant C(\|\Delta_{p,\Phi} s\|_{L^2} + p^2\|s\|_{L^2}).
\ee
Let $Q$ be a differential operator of order $m\in \bN$ with
scalar principal symbol and with  compact support
in $U_i$, then by (\ref{c12}) and $[\Delta_{p,\Phi},Q]$
is a   differential operator of order $m+1$, we get 
\begin{equation}\label{c13}
\begin{split}
\|Qs\|_{H^2_p}   
&\leqslant C(\|\Delta_{p,\Phi}Q s\|_{L^2} + p^2 \|Qs\|_{L^2})\\
 &\leqslant C(\|Q \Delta_{p,\Phi} s\|_{L^2} + p^2 \|s\|_{H^{m+1}_p}
+ p^2\|Qs\|_{L^2}).
\end{split}
\end{equation}
This means
\be\label{c17}
\|s\|_{H^{2m+2}_p} 
\leqslant C_m p^{4m} \sum_{j=0}^{m} \|\Delta_{p,\Phi}^j s\|_{L^2}.
\ee
Moreover for ${\bf G}_{k,p}=  G_{k,p}$ or $H_{k,p}$,
$\langle \Delta_{p,\Phi}^{m'}
{\bf G}_{k,p}(\frac{1}{\sqrt{p}}\Delta_{p,\Phi})Q s,s'\rangle
=\langle s,Q^ * {\bf G}_{k,p}(\frac{1}{\sqrt{p}}\Delta_{p,\Phi})
\Delta_{p,\Phi}^{m'} s'\rangle$, so
from \eqref{0c6}, \eqref{c9}, we know that for $l,m'\in \bN$,
there exist $C, C'>0$ such that for $p>1$,
\begin{equation}\label{c18}
\begin{split}
&\Big\|\Delta_{p,\Phi}^{m'} G_{k,p}\big(\tfrac{1}{\sqrt{p}}\Delta_{p,\Phi}\big)Qs\Big\|_{L^2}
  \leqslant C p^{-l}\|s\|_{L^2},\\
& \Big \|\Delta_{p,\Phi}^{m'}
\Big (H_{k,p}\big(\tfrac{1}{\sqrt{p}}\Delta_{p,\Phi}\big)-P_{\mH_p}\Big )Qs
\Big \|_{L^2}  \leqslant C' p^{2m-\frac{k}{2}}\|s\|_{L^2}. 
\end{split}
\end{equation}
We deduce from \eqref{c17} and \eqref{c18} that if $P,Q$
are differential operators with compact support in
$U_i$, $U_j$ respectively, then for any $l\in \bN$, there exists
$ C>0$ such that for $p>1$,
\begin{equation}\label{c19}
\begin{split}
&\Big\|P G_{k,p}\big(\tfrac{1}{\sqrt{p}}\Delta_{p,\Phi}\big) Q s\Big\|_{L^2}
\leqslant C p^{-l} \|s\|_{L^2},\\
&\Big \|P \Big (H_{k,p}\big(\tfrac{1}{\sqrt{p}}\Delta_{p,\Phi}\big)
-P_{\mH_p}\Big ) Q s\Big \|_{L^2}
\leqslant C    p^{2(m+m')-\frac{k}{2}} \|s\|_{L^2}.
\end{split}
\end{equation}

On $U_i\times U_j$, we use Sobolev inequality, we know for any $l\in \bN$, 
there exists $ C>0$ such that for $p>1$,
\begin{equation}\label{c20}
\begin{split}
&\Big |G_{k,p}\big(\tfrac{1}{\sqrt{p}}\Delta_{p,\Phi}\big)(x,x')\Big|_{\cC ^m}
 \leqslant C_{l,m} p^{-l},\\
&\Big | \Big (H_{k,p}\big(\tfrac{1}{\sqrt{p}}\Delta_{p,\Phi}\big)
-P_{0,p}\Big )(x,x')\Big|_{\cC ^m}
\leqslant C p^{2(2m+2n+1)-\frac{k}{2}}. 
\end{split}
\end{equation}
By \eqref{0c9} and \eqref{c20}, we get our Proposition \ref{0t3.0}.
\end{proof}

Using \eqref{0c3}, \eqref{0c5} and the finite propagation speed 
\cite[\S 7.8]{CP}, \cite[\S 4.4]{T1}, it is clear that for $x,x'\in X$, 
$F_k\big(\frac{1}{\sqrt{p}}\Delta_{p,\Phi}\big)(x,\cdot)$ only depends on 
the restriction of $\Delta_{p,\Phi}$ to $B^X(x,\var p^{-\frac{1}{4}})$, 
and $F_k\big(\frac{1}{\sqrt{p}}\Delta_{p,\Phi}\big)(x,x')= 0$,
if $d(x, x') \geqslant\var p^{-\frac{1}{4}}$.
 This means that the asymptotic of $\Delta_{p,\Phi}^q P_{\mH_p}(x,\cdot)$
when $p\to +\infty$,
modulo $\cO(p^{-\infty})$ (i.e. terms whose $\cC^m$ norm is 
$\cO(p^{-l})$ for any $l,m\in\bN$),
only depends on the restriction of $\Delta_{p,\Phi}$ to
$B^X(x,\var p^{-\frac{1}{4}})$.

\subsection{Rescaling and a Taylor expansion of the
operator $\Delta_{p,\Phi}$}\label{s3.2}

We fix $x_0\in X$. From now on, we identify $B^{T_{x_0}X}(0,\var)$
with $B^{X} (x_0,\var)$.
For $Z\in B^{T_{x_0}X}(0,\var)$ we identify $L_Z, E_Z$ and $(L^p\otimes E)_Z$
to $L_{x_0}, E_{x_0}$ and $(L^p\otimes E)_{x_0}$
by parallel transport with respect to the connections
$\nabla ^L$, $\nabla ^E$ and $\nabla^{L^p\otimes E}$ along the curve 
$\gamma_Z :[0,1]\ni u \to \exp^X_{x_0} (uZ)$.
 Let $\{e_i\}_i$ be an oriented orthonormal
basis of $T_{x_0}X$, and let $\{e ^i\}_i$ be its dual basis.

For $\var >0$ small enough, we will extend the geometric
objects from $B^{T_{x_0}X}(0,\var)$ to $\bR^{2n} \simeq T_{x_0}X$
where the identification is given by
\begin{equation}\label{0c11}
 (Z_1,\cdots, Z_{2n}) \in \bR^{2n} \longrightarrow \sum_i
Z_i e_i\in T_{x_0}X
\end{equation}
 such that $\Delta_{p,\Phi}$ is the
restriction of a renormalized Bochner-Laplacian on $\bR^{2n}$ associated
to a Hermitian line bundle with positive curvature. In this way,
we  replace $X$ by  $\bR^{2n}$.

At first, we denote by $L_0$, $E_0$ the trivial bundles with fiber $L_{x_0},
E_{x_0}$ on $X_0= \bR^{2n}$. We still denote by $\nabla ^L,
\nabla ^E$, $h^L$ etc. the connections and metrics on   $L_0$,
$E_0$ on $B^{T_{x_i}X}(0,4\var)$ induced by the above
identification. Then $h^L$, $h^E$ is  identified to the constant
metrics $h^{L_0}=h^{L_{x_0}}$,  $h^{E_0}=h^{E_{x_0}}$. 

Let $\rho: \bR\to [0,1]$ be a smooth even function such that
\begin{align}\label{1c14}
\rho (v)=1  \  \  {\rm if} \  \  |v|<2;
\quad \rho (v)=0 \   \   {\rm if} \  |v|>4.
\end{align}
 Let $\varphi_\var : \bR^{2n} \to \bR^{2n}$ is the map defined by
$\varphi_\var(Z)= \rho(|Z|/\var) Z$.
Then $\Phi_0=\Phi\circ \varphi_\var$ is a smooth self-adjoint
 section of $\End(E_0)$ on $X_0$.
Let $g^{TX_0}(Z)= g^{TX}(\varphi_\var(Z))$, $J_0(Z)= J(\varphi_\var(Z))$
be the metric and complex structure on $X_0$.
 Set $\nabla ^{E_0}= \varphi_\var ^* \nabla ^{E}$. Then  $\nabla ^{E_0}$
is the extension of $ \nabla ^{E}$ on $B^{T_{x_0}X}(0,\var)$.
If $\mR =
\sum_i Z_i e_i=Z$ denotes radial vector field on $\bR^{2n}$,
we define the Hermitian connection $\nabla ^{L_0}$ on $(L_0, h^{L_0})$ by
\begin{align}\label{1c15}
&\nabla ^{L_0}|_Z = \varphi_\var ^* \nabla ^{L} +\frac{1}{2}
(1-\rho ^2(|Z|/\var) )  R^{L}_{x_0} (\mR,\cdot).
\end{align}
Then we calculate easily that its curvature $R^{L_0}= (\nabla ^{L_0})^2$ is
\begin{equation}\label{1c16}
\begin{split}
R^{L_0}(Z) &= \varphi_\var ^* R^L + \frac{1}{2}d \Big((1-\rho ^2(|Z|/\var))
  R^{L}_{x_0} (\mR,\cdot)\Big)\\
&= \Big(1-\rho ^2(|Z|/\var)\Big)  R^{L}_{x_0}
+ \rho ^2(|Z|/\var) R^L_{\varphi_\var(Z)}\\
&- (\rho \rho')(|Z|/\var) \frac{Z_i e ^i}{\var |Z|}\wedge \big[
R^{L}_{x_0} (\mR,\cdot)- R^{L}_{\varphi_\var(Z)} (\mR,\cdot)\big].
\end{split}
\end{equation}
Thus $R^{L_0}$ is positive in the sense of (\ref{0.2}) for $\var$ small enough,
 and  the corresponding constant $\mu_0$ for  $R^{L_0}$ is bigger
than $\frac{4}{5}\mu_0$. From now on, we fix  $\var$ as above.

Let $\Delta_{p,\Phi_0}^{X_0}=\Delta^{L^p_0\otimes E_0}-p\tau_0-\Phi_0$ 
be the renormalized Bochner-Laplacian on $X_0$
associated to the above data, as in \eqref{laplace}. Observe that
$R^{L_0}$ is uniformly positive on $\bR^{2n}$,
so by the relations (3.2), (3.11) and (3.12) in \cite[p. 656\,--\,658]{MM}, we know that \eqref{0.3}
still holds for $\Delta_{p,\Phi_0}^{X_0}$. 
Especially, there exists $C_{L_0}>0$ such that
\begin{align}\label{1c17}
&\spec \Delta_{p,\Phi_0}^{X_0} 
\subset [-C_{L_0},C_{L_0}]\cup [\frac{8}{5} p\mu_0 -C_{L_0},+\infty[\,.
\end{align}
We note that $\Delta_{p,\Phi_0}^{X_0}$ has not necessarily discrete spectrum.

Let $S_L$ be an unit vector of $L_{x_0}$. Using $S_L$ and the above 
discussion, we get an isometry $E_0\otimes L_0^p \simeq  E_{x_0}$.
Let $P_{0,\mH_p}$ be the spectral projection of $\Delta_{p,\Phi_0}^{X_0}$ 
from $\cC^\infty (X_0,L^p_0\otimes E_0)\simeq \cC^\infty (X_0,E_{x_0})$
corresponding to the interval
$[-C_{L_0},C_{L_0}]$, and let $P_{0,q,p}(x,x')$
$(q\geqslant0)$ be the smooth kernels of 
$P_{0,q,p}=(\Delta_{p,\Phi_0}^{X_0})^q P_{0,\mH_p}$ 
(we set $(\Delta_{p,\Phi_0}^{X_0})^0=1$)
with respect to the volume form $dv_{X_0}(x')$.
The following Proposition shows that $P_{q,p}$ and 
$P_{0,q,p}$ are asymptotically close on $B^{T_{x_0}X}(0,\var)$ in the $\cC^\infty$--\,topology, as $p\to\infty$.
\begin{prop} \label{p3.2} For any $l,m\in \bN$, there exists $C_{l,m}>0$
 such that for $x,x' \in B^{T_{x_0}X}(0,\var)$,
\begin{equation}\label{1c19}
\left |(P_{0,q,p} -P_{q,p})(x,x')\right |_{\cC ^m}
\leqslant C_{l,m} p^{-l}.
\end{equation}
\end{prop}
\begin{proof} Using \eqref{0c3} and \eqref{1c17}, we know that
for $x,x' \in B^{T_{x_0}X}(0,\var)$,
\begin{align}\label{1c20}
\left|F_k\big(\tfrac{1}{\sqrt{p}}\Delta_{p,\Phi}\big)(x,x') 
- P_{0,0,p}(x,x')\right|_{\cC ^m} \leqslant C_{k,m} p^{-\frac{k}{2} +2(m+n+1)}.
\end{align}
 Thus from \eqref{0c7} and \eqref{1c20} for $k$ big enough,
we infer \eqref{1c19} for $q=0$; Now from the definition of
$P_{0,q,p}$ and $P_{q,p}$, we get \eqref{1c19} from \eqref{c10}
and \eqref{1c19} for $q=0$.
\end{proof}
It suffices therefore to study the kernel $P_{0,q,p}$ and for this purpose we
rescale the operator $\Delta_{p,\Phi_0}^{X_0}$. 
Let $dv_{TX}$ be the Riemannian volume form of
$(T_{x_0}X, g^{T_{x_0}X})$.
Let $\kappa (Z)$ be the smooth positive function defined by the equation
\be\label{c22}
dv_{X_0}(Z) = \kappa (Z) dv_{TX}(Z),
\ee
with $\kappa(0)=1$.
 Denote by  $\nabla_U$ the ordinary differentiation
 operator on $T_{x_0}X$ in the direction $U$, 
and set $\partial_i=\nabla _{e_i}$.
If $\alpha = (\alpha_1,\cdots, \alpha_{2n})$ is a multi-index,
set $Z^\alpha = Z_1^{\alpha_1}\cdots Z_{2n}^{\alpha_{2n}}$. 
We also denote by
$(\partial ^\alpha R^L)_{x_0}$ the tensor $(\partial^\alpha R^L)_{x_0}(e_i,e_j)
=\partial ^\alpha( R^L(e_i,e_j))_{x_0}$. 
Denote by $t=\frac{1}{\sqrt{p}}$.
For $s \in \cC ^{\infty}(\bR^{2n}, E_{x_0})$ and $Z\in \bR^{2n}$, set
\begin{equation}\label{c27}
\begin{split}
(S_{t} s ) (Z) = & s (Z/t),
\quad  \nabla_{t} =  tS_t^{-1}\kappa ^{\frac{1}{2}}
\nabla ^{L^p_0\otimes E_0}\kappa ^{-\frac{1}{2}}S_t,\\
& \cL_t= S_t^{-1} \,\tfrac{1}{p}\,\kappa ^{\frac{1}{2}}\,
\Delta_{p,\Phi_0}^{X_0} \,\kappa ^{-\frac{1}{2}}\,S_t. 
\end{split}
\end{equation}
The operator $\cL_t$ is the rescaled operator, which we now develop in Taylor series. 
\begin{thm}\label{t3.3} There exist polynomials $\mA_{i,j,r}$ 
{\rm(} resp. $\mB_{i,r}$, $\mC_{r}${\rm)}
{\rm(}$r\in \bN, i,j\in \{1,\cdots, 2n\}${\rm)}
in $Z$ with the following properties:

-- their coefficients are polynomials in
$R^{TX}$  {\rm(}resp. $R^{TX}$,  $R^L$, $R^{E}$, $\Phi${\rm)}
and their derivatives at $x_0$ up to order $r-2$ {\rm(}resp. $r-1$, $r$,
 $r-1$, $r${\rm)}\,, 

-- $\mA_{i,j,r}$ is a monomial in $Z$ of degree $r$, 
the degree in $Z$ of $\mB_{i,r}$ {\rm(}resp. $\mC_{r}${\rm)} has the same parity 
with $r-1$ {\rm(}resp. $r${\rm)}\,, 

-- if we denote by
\begin{align}\label{0c35}
\mO_{r} =  \mA_{i,j,r}\nabla_{e_i}\nabla_{e_j}
+ \mB_{i,r}\nabla_{e_i}+ \mC_{r},
\end{align}
 then
\begin{align}\label{c30}
\cL_t=  \cL_0+ \sum_{r=1}^m t^r \mO_r + \cO(t^{m+1}).
\end{align}
and there exists $m'\in \bN$ such that for any $k\in \bN$, $t\leqslant 1$
the derivatives of order $\leqslant k$ of the coefficients of the operator
 $\cO(t^{m+1})$ are dominated by $C t^{m+1} (1+|Z|)^{m'}$. Moreover
\begin{align}\label{c31}
&\cL_0 = -\sum_j  \Big(\nabla_{e_j}+\frac{1}{2} R^L_{x_0}(Z, e_j) \Big)^2 
-\tau_{x_0},\\
&\mO_1(Z)= -\frac{2}{3}   ( \partial_j R^L)_{x_0} (\mR,e_i)Z_j
\Big(\nabla _{e_i}+ \frac{1}{2}  R^L_{x_0}(\mR, e_i)\Big)
 -\frac{1}{3} (\partial_i R^L)_{x_0} (\mR,
e_i)
- (\nabla _\mR \tau)_{x_0},\nonumber\\
&\mO_2(Z)=    \frac{1}{3} \left \langle R^{TX}_{x_0} (\mR,e_i)
\mR, e_j\right \rangle_{x_0} \Big( \nabla_{e_i}  +  \frac{1}{2}
R^L_{x_0}(\mR,e_i)\Big)\Big( \nabla_{e_j}
+  \frac{1}{2} R^L_{x_0}(\mR,e_j)\Big)\nonumber\\
 &\hspace*{3mm} +\Big [\frac{2}{3} \left \langle R^{TX}_{x_0}
(\mR, e_j) e_j,e_i\right \rangle_{x_0} 
- \Big(\frac{1}{2}\sum_{|\alpha|=2}(\partial ^{\alpha}R^L)_{x_0} 
\frac{Z^\alpha}{\alpha !} + R^E_{x_0} \Big)(\mR,e_i)\Big ]
\Big(\nabla_{e_i} + \frac{1}{2} R^L_{x_0}(\mR,e_i)\Big)\nonumber\\
&\hspace*{3mm} -\frac{1}{4} \nabla_{e_i}\Big(\sum_{|\alpha|=2}
(\partial ^{\alpha}R^L)_{x_0} \frac{Z^\alpha}{\alpha !}(\mR,e_i)\Big)
 -\frac{1}{9}\sum_i \Big[\sum_j (\partial_j R^L)_{x_0} (\mR,e_i)Z_j\Big]^2
\nonumber\\
&\hspace*{3mm}-\frac{1}{12}\Big[\cL_0, \left \langle R^{TX}_{x_0} (\mR,e_i)
\mR, e_i\right \rangle_{x_0} \Big]
-\sum_{|\alpha|=2}(\partial ^{\alpha}\tau)_{x_0}
\frac{Z^\alpha}{\alpha !}+  \Phi_{x_0} .\nonumber
\end{align}
\end{thm}
\begin{proof}
Set $g_{ij}(Z)= g^{TX}(e_i,e_j)(Z) =  \langle e_i,e_j\rangle_Z$ 
and let $(g^{ij}(Z))$ be the inverse of
the matrix $(g_{ij}(Z))$.
By \cite[Proposition 1.28]{BeGeV}, the Taylor expansion of $g_{ij}(Z)$
with respect to the basis $\{e_i\}$ to order $r$ is a polynomial of
 the Taylor expansion of $R^{TX}$
to order $r-2$, moreover
\begin{equation}\label{0c30}
\begin{split}
&g_{ij}(Z) =
\delta_{ij} +  \frac{1}{3}
\left \langle R^{TX}_{x_0} (\mR,e_i) \mR, e_j\right \rangle_{x_0}
 + \cO (|Z|^3),\\
&\kappa(Z)= |\det (g_{ij}(Z))|^{1/2}  = 1 + 
\frac{1}{6} \left \langle R^{TX}_{x_0} (\mR,e_i) \mR, e_i\right \rangle_{x_0}
 + \cO (|Z|^3).
\end{split}
\end{equation}
If $\Gamma _{ij}^l$ is the connection form of $\nabla ^{TX}$
with respect to the basis $\{e_i\}$, we have $(\nabla
^{TX}_{e_i}e_j)(Z) = \Gamma _{ij}^l (Z) e_l$. Owing to \eqref{0c30},
\begin{equation}\label{0c31}
\begin{split}
\Gamma _{ij}^l (Z)& =  \frac{1}{2} g^{lk} (\partial_i g_{jk}
+ \partial_j g_{ik}-\partial_k g_{ij})(Z)\\
&= \frac{1}{3}\Big  [ 
\left \langle R^{TX}_{x_0} (\mR, e_j) e_i, e_l\right \rangle _{x_0}
+ \left \langle R^{TX}_{x_0} (\mR, e_i) e_j, e_l\right \rangle_{x_0}\Big ]
 + \cO(|Z|^2). 
\end{split}
\end{equation}
Now by \eqref{0c1},
\begin{align}\label{c32}
\Delta_{p,\Phi} = - g^{ij} ( \nabla ^{L^p\otimes E}_{e_i}
\nabla ^{L^p\otimes E}_{e_j}- \nabla ^{L^p\otimes E}_{\nabla ^{TX}_{e_i}e_j} )
-p \tau +\Phi.
\end{align}
so from \eqref{c27} and \eqref{c32} we infer the expression
\begin{align}\label{0c37}
\cL_t=-g^{ij} (tZ) \Big [ \nabla_{t, e_i}\nabla_{t, e_j} - t
\Gamma _{ij}^l \nabla_{t, e_l} \Big ](tZ) - \tau (tZ)
+t^2\Phi(tZ).
\end{align}
Let $\Gamma ^E$, $\Gamma ^L$
be the connection forms of $\nabla^E$ and $\nabla^L$
 with respect to any fixed frames for $E$, $L$
 which are parallel along the curve $\gamma_Z$
under our trivializations on $B^{T_{x_0}X}(0,\var)$.
\eqref{c27} yields on $B^{T_{x_0}X}(0, \var/t)$
\begin{align}\label{0c36}
&\nabla_{t, e_i}|_{Z} = \kappa ^{\frac{1}{2}}(tZ)\Big(\nabla_{e_i}
+ \frac{1}{t} \Gamma ^L (e_i)(tZ) + t \Gamma ^E (e_i) (tZ)\Big)
\kappa ^{-\frac{1}{2}}(tZ).
\end{align}

Let $\Gamma ^\bullet= \Gamma ^E, \Gamma ^L$ and
$R^\bullet= R^E, R^L$, respectively.
By \cite[Proposition 1.18]{BeGeV}
the Taylor coefficients of $\Gamma ^\bullet (e_j) (Z)$ at $x_0$
to order $r$ are only determined by those of $R^\bullet$ to order $r-1$, and
\begin{align}\label{0c39}
\sum_{|\alpha|=r}  (\partial^\alpha
 \Gamma ^\bullet ) _{x_0} (e_j) \frac{Z^\alpha}{\alpha !}
=\frac{1}{r+1} \sum_{|\alpha|=r-1}
(\partial^\alpha R^\bullet ) _{x_0}(\mR, e_j) 
  \frac{Z^\alpha}{\alpha !}.
\end{align}

Owing to \eqref{0c30}, \eqref{0c39}
\begin{align}\label{0c41}
& \cL_t = - \Big ( \delta_{ij} - \frac{t^2}{3}
\left \langle R^{TX}_{x_0} (\mR,e_i) \mR, e_j\right \rangle
 + \cO (t^3)\Big ) \kappa ^{\frac{1}{2}}(tZ) \Big \{\\
&\Big [\nabla_{e_i}+ \Big(\frac{1}{2} R^L_{x_0} 
+ \frac{t}{3} (\partial_k R^L)_{x_0} Z_k +\frac{t^2}{4}\sum_{|\alpha|=2}
(\partial ^{\alpha}R^L)_{x_0} \frac{Z^\alpha}{\alpha !}
+ \frac{t^2}{2} R^E_{x_0}\Big ) (\mR,e_i)+ \cO (t^3)\Big]  \nonumber \\
&\Big [\nabla_{e_j}+\Big (\frac{1}{2} R^L_{x_0} + \frac{t}{3}
(\partial_k R^L)_{x_0} Z_k +\frac{t^2}{4}\sum_{|\alpha|=2}
(\partial ^{\alpha}R^L)_{x_0}  \frac{Z^\alpha}{\alpha !}
+ \frac{t^2}{2} R^E _{x_0} \Big )(\mR,e_j) + \cO (t^3)\Big]   \nonumber \\
&\left.- t \Gamma ^l_{ij}(tZ) \Big (\nabla_{e_l}+
\frac{1}{2} R^L_{x_0}(\mR,e_l)
+ \cO (t)\Big ) \right \} \kappa ^{-\frac{1}{2}}(tZ) \nonumber \\
&-\tau_{x_0} -t (\nabla_{\mR} \tau)_{x_0}
-t^2\sum_{|\alpha|=2}(\partial ^{\alpha}\tau)_{x_0}
\frac{Z^\alpha}{\alpha !}+ t^2 \Phi_{x_0} + \cO (t^3). \nonumber
\end{align}
Relations \eqref{0c30} and \eqref{0c37}\,--\,\eqref{0c41} settle 
our Theorem.
\end{proof}

\subsection{Uniform estimate of the generalized Bergman kernels}
\label{s3.3}

We shall estimate the Sobolev norm of the resolvent of $\cL_t$ so we introduce 
the following norms.
We denote by $\left \langle\,\cdot\,,\cdot\,\right\rangle_{0,L^2}$ and $\norm{\,\cdot\,}_{0,L^2}$
the scalar product and the $L^2$ norm on $\cC ^\infty (X_0, E_{x_0})$
induced by $g^{TX_0}, h^{E_0}$ as in \eqref{0c2}.
For $s\in \cC ^{\infty}(X_0, E_{x_0}) $, set
\begin{align}\label{u0}
&\|s\|_{t,0}^2= \|s\|_{0}^2 
= \int_{\bR^{2n}} |s(Z)|^2_{h^{E_{x_0}}}dv_{TX}(Z),\\
&\| s \|_{t,m}^2 = \sum_{l=0}^m \sum_{i_1,\cdots, i_l=1}^{2n}
\| \nabla_{t,e_{i_1}} \cdots \nabla_{t,e_{i_l}} s\|_{t,0}^2. \nonumber
\end{align}
We denote by $\left \langle s ', s \right\rangle_{t,0}$
 the inner product on $\cC^\infty (X_0, E_{x_0})$
 corresponding to $\norm{\,\cdot\,}^2_{t,0}$\,.
Let $H^m_t$ be the Sobolev space of order $m$ with norm $\norm{\,\cdot\,}_{t,m}$.
Let $H^{-1}_t$ be the Sobolev space of order $-1$ and 
let $\norm{\,\cdot\,}_{t,-1}$ be the norm on  $H^{-1}_t$ defined by
$\|s\|_{t,-1} = \sup_{0\neq s'\in  H^1_t }$
$|\left \langle s,s'\right \rangle_{t,0}|/\|s'\|_{t,1}$.
If $A\in \cL (H^{m}, H^{m'})$ $(m,m' \in \bZ)$, we denote by 
$\norm{A}^{m,m'}_t$ the norm of $A$ with respect to the norms
$\norm{\,\cdot\,}_{t,m}$ and $\norm{\,\cdot\,}_{t,m'}$.

\begin{rem}\label{0t3.3} Note that $\Delta_{p,\Phi_0}^{X_0} $ is 
self-adjoint with respect to $\norm{\,\cdot\,}_{0}$, thus
by \eqref{c22}, \eqref{c27}, \eqref{u0}, 
 $\cL_t$ is a formal self adjoint elliptic operator with respect to  
$\|\quad\|_{0}$,
and is a smooth family of operators with the parameter $x_0\in X$.
Thus $\cL_0$ and $\mO_r$ are also formal self-adjoint with respect to  
$\norm{\,\cdot\,}_{0}$. This will simplify the computation of the coefficients
$b_{0,1}$ in (\ref{0.7}) (cf.  \S \ref{s4.3})
and explains why we prefer to conjugate with $\kappa ^{1/2}$ 
comparing to \cite[(3.38)]{DLM}.
\end{rem}

\begin{thm}  \label{tu1} There exist constants $ C_1, C_2, C_3>0$
such that for $t\in ]0,1]$ and any
$s,s'\in \cC ^{\infty}_0(\bR^{2n}, E_{x_0})$,
\begin{equation}\label{u1}
\begin{split}
& \left \langle \cL_t s,s\right \rangle_{t,0}
\geqslant C_1\|s\|_{t,1}^2 -C_2 \|s\|_{t,0}^2  , \\
& |\left \langle \cL_t s, s'\right \rangle_{t,0}|
 \leqslant C_3 \|s\|_{t,1}\|s'\|_{t,1}.
\end{split}
\end{equation}
\end{thm}

\begin{proof}
Relations \eqref{laplace} and \eqref{0c1} yield
\begin{align}\label{1u1}
 \left \langle \Delta_{p,\Phi}s,s\right \rangle_{0,L^2} =
\|\nabla ^{ L^p_0\otimes  E_0}s\|_{0,L^2}^2 - \left \langle \left (
p\tau-\Phi\right )s,s\right \rangle_{0,L^2}.
\end{align}
Thus from (\ref{c27}), (\ref{u0}) and (\ref{1u1}) we get
\begin{align}\label{1u2}
\left \langle  \cL_t s,s\right \rangle_{t,0} = \|\nabla_ts\|_{t,0}^2
- \left \langle\left  (
S_t^{-1}\tau  -t^2\Phi\right )s,s\right \rangle_{t,0}.
\end{align}
which implies \eqref{u1}.
\end{proof}

Let $\delta$ be the counterclockwise oriented circle in $\bC$
of center $0$ and radius $\mu_0/4$.

\begin{thm}\label{tu4} The resolvent $(\lambda- \cL_t)^{-1}$ exists for 
$\lambda \in \delta$, and there exists $C>0$ such that for $t\in ]0,1]$,
$\lambda \in  \delta$, and $x_0\in X$,
\begin{align}\label{ue2}
& \| (\lambda- \cL_t)^{-1}\|^{0,0}_t \leqslant C,
& \| (\lambda- \cL_t)^{-1}\|^{-1,1}_t
\leqslant C .
\end{align}
\end{thm}
\begin{proof} By \eqref{1c17}, \eqref{c27}, for $t$ small enough,
\begin{align}\label{1u3}
\spec \, \cL_t\subset 
\big[-C_{L_0} t^2, C_{L_0} t^2\big] \cup \Big[\, \,\mu_0,+\infty\Big[.
\end{align}
Thus the resolvent $(\lambda- \cL_t)^{-1}$ exists for $\lambda \in \delta$,
and we get the  first inequality of (\ref{ue2}).
By \eqref{u1}
$(\lambda_0- \cL_t)^{-1}$ exists for $\lambda_0\in\bR$, 
$\lambda_0\leqslant -2C_2$, and
$\|(\lambda_0- \cL_t)^{-1}\|^{-1,1}_t \leqslant \frac{1}{C_1}$. Now,
\be\label{ue7}
(\lambda- \cL_t)^{-1}= (\lambda_0- \cL_t)^{-1}
- (\lambda-\lambda_0) (\lambda- \cL_t)^{-1}(\lambda_0- \cL_t)^{-1}.
\ee
Thus for  $\lambda\in \delta$, from (\ref{ue7}), we get
\be\label{ue8}
\|(\lambda-\cL_t)^{-1}\|^{-1,0}_t  \leq
\frac{1}{C_1} \left(1+\frac{4}{\mu_0}|\lambda-\lambda_0|\right).
\ee
Changing the last two factors in \eqref{ue7} and applying \eqref{ue8}
we get
\begin{align}\label{ue9}
\|(\lambda-\cL_t)^{-1}\|^{-1,1}_t  \leqslant  \frac{1}{C_1}+
 \frac{|\lambda-\lambda_0|}{{C_1}^2} 
 \left(1+\frac{4}{\mu_0}|\lambda-\lambda_0|\right)\leqslant C.
\end{align}
The proof of our Theorem is complete.
\end{proof}

\begin{prop} \label{tu5} Take $m \in \bN^*$. There exists $C_m>0$ such that
 for   $t\in ]0,1]$,  $Q_1, \cdots, Q_m$
$\in \{ \nabla_{t,e_i}, Z_i\}_{i=1}^{2n}$ and  $s,s'\in C
^{\infty}_{0}(X_0, E_{x_0})$, \be\label{ue11} \left |\left \langle
[Q_1, [Q_2,\ldots, [Q_m,  \cL_t]] \ldots ]s, s'\right
\rangle_{t,0} \right | \leqslant C_m  \|s\|_{t,1} \|s'\|_{t,1}. \ee
\end{prop}
\begin{proof}
Note that $[\nabla_{t,e_i}, Z_j]=\delta_{ij}$, hence
\eqref{0c37} implies that $[Z_j, \cL_t]$ verifies (\ref{ue11}).
On the other hand, we obtain from \eqref{c27} 
\begin{align}\label{1ue2}
[\nabla_{t,e_i},\nabla_{t,e_j}]= \left(R^{L_0}(tZ)
+ t^2  R^{E_0}(tZ)\right)(e_i,e_j).
\end{align}
Thus from (\ref{0c37}) and (\ref{1ue2}), 
we know that $[\nabla_{t,e_k},  \cL_t]$
has the same structure as  $\cL_t$ for $t \in ]0,1]$,
i.e. $[\nabla_{t,e_k},  \cL_t]$ has the same type as
\begin{align}\label{1ue3}
\sum_{ij} a_{ij}(t,tZ) \nabla_{t,e_i}\nabla_{t,e_j}
+\sum_{i}b_{i}(t,tZ) \nabla_{t,e_i} + c(t,tZ),
\end{align}
and $a_{ij}(t,Z),b_{i}(t,Z),  c(t,Z)$ and their derivatives in $Z$
are uniformly bounded for $Z\in \bR ^{2n}, t\in [0,1]$.
Moreover they are polynomials in $t$.

If $(\nabla_{t,e_i})^*$ is the adjoint of $\nabla_{t,e_i}$ with respect to
$\left \langle \,\cdot\, ,\,\cdot\,\right \rangle_{t,0}$\,, \eqref{u0} yields
\begin{align}\label{1ue4}
(\nabla_{t,e_i})^* =- \nabla_{t,e_i}- t  (\kappa ^{-1}(e_i \kappa))(tZ).
\end{align}
Thus by (\ref{1ue3}) and (\ref{1ue4}), (\ref{ue11}) is verified for $m=1$.

By recurrence, it transpires that
$[Q_1,[Q_2,\ldots, [Q_m,  \cL_t]] \ldots ]$
has the same structure \eqref{1ue3} as $\cL_t$\,, so from \eqref{1ue4}
we get the required assertion.
\end{proof}

\begin{thm}\label{tu6} For any $t\in ]0,1]$, $\lambda \in \delta$,
 $m \in \bN$, the resolvent $(\lambda-\cL_t)^{-1}$ maps $H^m_t$
into $H^{m+1}_t$. Moreover for any $\alpha\in \bZ^{2n}$,
 there exists $C_{\alpha, m}>0$
such that for  $t\in ]0,1]$, $\lambda \in \delta$,
 $s\in \cC ^\infty (X_0,E_{x_0})$,
\be\label{ue12}
\|Z^\alpha (\lambda-\cL_t)^{-1} s\|_{t,m+1}
\leqslant C_{\alpha, m} \sum_{\alpha' \leqslant \alpha}
\|Z^{\alpha'} s\|_{t,m} .
\ee
\end{thm}
\begin{proof} For $Q_1, \cdots, Q_m\in \{\nabla_{t,e_i}\}_{i=1}^{2n}$,
$Q_{m+1},\cdots, Q_{m+|\alpha|}\in\{Z_i\}_{i=1}^{2n}$,  we can
express $Q_1\cdots$ $Q_{m+|\alpha|}(\lambda-\cL_t)^{-1}$ as a
linear combination of operators of the type \be\label{ue13} [Q_1 [
Q_2,\ldots [Q_{m'},(\lambda-\cL_t)^{-1}]]\ldots ]Q_{m'+1} \cdots
Q_{m+|\alpha|}, \quad m'\leqslant m+|\alpha|. \ee Let $\cR_{t}$ be the
family of operators $\cR_{t} = \{ [ Q_{j_1}[Q_{j_2},\ldots
[Q_{j_l},\cL_t]]\ldots ] \}$.  Clearly, any commutator $[Q_1 [
Q_2,\ldots [Q_{m'},(\lambda-\cL_t)^{-1}]]\ldots ]$ is a linear
combination of operators of the form \be\label{ue14}
(\lambda-\cL_t)^{-1}R_1(\lambda-\cL_t)^{-1}R_2 \cdots
R_{m'}(\lambda-\cL_t)^{-1} \ee with  $R_1, \cdots, R_{m'} \in
\cR_{t}$.

From Proposition \ref{tu5} we deduce that the norm $\norm{\,\cdot\,}_t^{1,-1}$ of the
 operators $R_j\in \cR_{t}$ is
uniformly bounded by $C$.
By Theorem \ref{tu4} there exists $C>0$, 
such that the norm $\norm{\,\cdot\,}_t^{0,1}$
of operators \eqref{ue14} is dominated by
$C$.
\end{proof}

The next step is to convert the estimates for the resolvent into estimates for
the spectral projection
$\mP_{0,t}:(\cC^{\infty}(X_0,E_{x_0}), \norm{\,\cdot\,}_0)\to
(\cC^{\infty}(X_0,E_{x_0}), \norm{\,\cdot\,}_0)$  
of $\cL_t$ corresponding to the interval $[-C_{L_0} t^2, C_{L_0}t^2]$. 
Let $\mP_{q,t}(Z,Z')=\mP_{q,t,x_0}(Z,Z')$, (with $Z,Z'\in X_0$, $q\geqslant0$) be the
smooth kernel of $\mP_{q,t}=(\cL_t)^q \mP_{0,t}$ 
(we set $(\cL_t)^0=1$) with respect to $dv_{TX}(Z')$. Note that
$\cL_t$ is a family of differential operators on $T_{x_0}X$ with coefficients
in $\End (E)_{x_0}$. Let $\pi : TX\times_{X} TX \to X$ be the
natural projection from the fiberwise product of $TX$ on $X$. Then
we can view  $\mP_{q,t}(Z,Z')$ as a smooth section of $\pi ^* (\End
(E))$ over  $TX\times_{X} TX$
by identifying a section $S\in \cC^\infty (TX\times_{X}TX,\pi ^* \End (E))$
with the family $(S_x)_{x\in X}$, where
$S_x=S|_{\pi^{-1}(x)}$.
 Let $\nabla ^{\End (E)}$ be the connection on $\End (E)$ induced by
$\nabla ^E$. Then $\nabla^{\pi^*\End (E)}$
induces naturally a $\cC ^m$--\,norm of $S$ for the parameter $x_0\in X$.
\begin{thm}\label{tue8}
 For any $m,m'\in\bN$, $\sigma>0$, there exists $C>0$,
 such that for $t\in ]0,1]$, $Z,Z'\in T_{x_0}X$, $|Z|,|Z'|\leqslant \sigma$,
 \begin{align}\label{ue15}
&\sup_{|\alpha|,|\alpha'|, r\leqslant m} \Big
|\frac{\partial^{|\alpha|+|\alpha'|}} {\partial Z^{\alpha}
{\partial Z'}^{\alpha'}} \frac{\partial^{r}}{\partial t^{r}}
\mP_{q,t}\left (Z, Z'\right )\Big |_{\cC ^{m'}(X)}   \leqslant C.
\end{align}
Here $\cC ^{m'}(X)$ is the $\cC ^{m'}$ norm for the parameter $x_0\in X$.
\end{thm}
\begin{proof} By (\ref{1u3}), for any $k\in \bN^*$, $q\geqslant0$,
\begin{align}\label{1ue15}
&\mP_{q,t}=(\cL_t)^q \mP_{0,t}= \frac{1}{2\pi i} \begin{pmatrix}
q+k-1 \\ k-1 \end{pmatrix} ^{-1} \int_{\delta} \lambda  ^{q+k-1}
(\lambda - \cL_t)^{-k} d \lambda.
\end{align}
For $m\in \bN$, let $\mQ ^m$ be the set of operators
 $\{\nabla_{t,e_{i_1}}\cdots  \nabla_{t,e_{i_j}}\}_{j\leqslant m}$.
From Theorem \ref{tu6}, we deduce that
if $Q\in \mQ ^{m}$, there is $C_m>0$ such that
\be\label{ue18}
\| Q (\lambda-\cL_t)^{-m}\|_t^{0,0}
 \leqslant C_m\,,\quad\text{for all $\lambda\in\delta$ }\, .
\ee
Observe that $\cL_t$ is self--adjoint with respect to $\norm{\,\cdot\,}_{t,0}$,
so after taking the adjoint of (\ref{ue18}), we have
\be\label{ue20}
\| (\lambda-\cL_t)^{-m}Q\|_t^{0,0}
\leqslant C_m\,.
\ee
From  (\ref{1ue15}), (\ref{ue18}) and (\ref{ue20}), we obtain
\begin{align}\label{ue21}
&\|Q \mP_{q,t}Q' \|^{0,0}_t \leqslant C_{m}\,,\quad\text{for $Q,Q' \in \mQ ^{m}$}\,.
\end{align}
Let $|\,\cdot\,|_{(\sigma),m}$ be the usual Sobolev norm on
$\cC^\infty(B^{T_{x_0}X}(0,\sigma+1),
E_{x_0})$ induced by $h^{E_{x_0}}$ and
the volume form $dv_{TX}(Z)$ as in (\ref{u0}).
Let $\|A\|_{(\sigma),m}$ be the  operator norm of $A$
with respect to $|\,\cdot\,|_{(\sigma),m}$.  Observe that by
(\ref{0c36}), (\ref{u0}), for $m>0$, 
there exists $C_\sigma>0$ such that for $s\in
\cC^\infty (X_0, E_{x_0})$, $\supp s \subset B^{T_{x_0}X}(0,\sigma +1)$, 
\begin{align}\label{1ue21}
&\frac{1}{C_\sigma}  \|s\|_{t,m}\leqslant |s|_{(\sigma),m}
\leqslant C_\sigma \|s\|_{t,m}.
\end{align}
Now \eqref{ue21} and \eqref{1ue21} together with Sobolev's inequalities
imply
\begin{equation}\label{ue23}
\sup_{|Z|,|Z'|\leqslant \sigma}| Q_Z Q'_{Z'} \mP_{q,t}(Z,Z')  | \leqslant C\,,\quad\text{for $Q,Q' \in \mQ ^{m}$}\,.
\end{equation}
Thanks to \eqref{0c36} and \eqref{ue23} estimate \eqref{ue15} holds for $r=m'=0$.

\noindent
To obtain (\ref{ue15}) for $r \geqslant1$ and $m'=0$, note that from (\ref{1ue15}), 
 \begin{align}\label{ue28}
\frac{\partial^{r}}{\partial t^{r}} \mP_{q,t}=& \frac{1}{2\pi i}
\begin{pmatrix} q+k-1 \\ k-1 \end{pmatrix} ^{-1} \int_{\delta}
\lambda ^{q+k-1} \frac{\partial^{r}}{\partial
t^{r}}(\lambda-\cL_t)^{-k} d \lambda\,,\quad\text{for $k\geqslant1$}\,.
\end{align}
 Set
\be\label{ue29}
I_{k,r} = \Big \{ (\bk,\br)=(k_i,r_i)| 
\sum_{i=0}^j k_i =k+j, \sum_{i=1}^j r_i =r,\, \,  k_i, r_i \in \bN^*\Big \}.
\ee
Then  there exist $a ^{\bk}_{\br} \in \bR$ such that
\begin{equation}\label{ue30}
\begin{split}
& A^{\bk}_{\br}  (\lambda,t) = (\lambda-\cL_t)^{-k_0}
\frac{\partial^{r_1}\cL_t}{\partial t^{r_1}}  (\lambda-\cL_t)^{-k_1}
\cdots\frac{\partial^{r_j}\cL_t}{\partial t^{r_j}}  (\lambda-\cL_t)^{-k_j},\\
& \frac{\partial^{r}}{\partial t^{r}}
(\lambda-\cL_t)^{-k}=
\sum_{(\bk,\br)\in I_{k,r} }
 a ^{\bk}_{\br}  A ^{\bk}_{\br}  (\lambda,t).
\end{split}
\end{equation}

We claim that $A ^{\bk}_{\br}(\lambda,t)$ is well defined and
for any $m\in \bN$, $k>2(m+r+1)$, $Q,Q'\in \mQ^m$,
there exist $C>0$, $N\in \bN$ such that for $\lambda\in \delta$,
\begin{align}\label{1ue30}
\|Q  A ^{\bk}_{\br}(\lambda,t)Q' s\|_{t,0}
\leq C  \sum_{|\beta|\leq 2r} \|Z^\beta s\|_{t,0}.
\end{align}

In fact, by (\ref{0c37}), $\frac{\partial^{r}}{\partial t^{r}}\cL_t$ is
combination of $\frac{\partial^{r_1}}{\partial t^{r_1}}(g^{ij}(tZ))$
$(\frac{\partial^{r_2}}{\partial t^{r_2}}\nabla_{t,e_i})$
$(\frac{\partial^{r_3}}{\partial t^{r_3}}\nabla_{t,e_j})$,
 $\frac{\partial^{r_1}}{\partial t^{r_1}}(d(tZ))$,
$\frac{\partial^{r_1}}{\partial t^{r_1}}(d_{i}(tZ))$
$(\frac{\partial^{r_2}}{\partial t^{r_2}}\nabla_{t,e_i})$.
Now $\frac{\partial^{r_1}}{\partial t^{r_1}}(d(tZ))$
(resp. $\frac{\partial^{r_1}}{\partial t^{r_1}}\nabla_{t,e_i}$)
 ($r_1\geq 1$),
 are functions of  the type as
 $d'(tZ)Z^\beta$, $|\beta|\leq r_1$ (resp. $r_1+1$)
and $d'(Z)$ and its derivatives on $Z$ are  bounded smooth functions on $Z$.

Let $\cR'_t$ be the family of operators of the type
$$\cR'_{t} = \{ [f_{j_1} Q_{j_1}, 
[f_{j_2} Q_{j_2},\ldots [f_{j_l} Q_{j_l},\cL_t]]\ldots ] \}$$
with $f_{j_i}$ smooth bounded (with its derivatives) functions
and  $Q_{j_i}\in \{\nabla_{t,e_l}\}_{l=1}^{2n}$.

Now for the operator $A ^{\bk}_{\br}(\lambda,t)Q'$, 
we will move first all the term $Z^\beta$ in $d'(tZ)Z^\beta$ 
as above to the right hand side of this operator, to do so, 
we always use the commutator trick, i.e., each time, 
we consider only the commutation for $Z_i$,
 not for $Z^\beta$ with $|\beta|>1$. 
Then $A ^{\bk}_{\br}(\lambda,t)Q'$ is as the form 
$\sum_{|\beta|\leq 2r} L^t_\beta Q''_\beta  Z^\beta$, and $Q''_\beta$ is obtained from 
 $Q'$ and its commutation with $Z^\beta$.
Now we move all 
the terms $\nabla_{t,e_i}$ in 
$\frac{\partial^{r_j} \cL_t}{\partial t^{r_j}}$ 
to the right hand side of the operator $L^t_\beta$.
Then as in the proof of Theorem \ref{tu6}, 
we get finally  that $Q  A ^{\bk}_{\br}(\lambda,t)Q'$
is as the form  $\sum_{|\beta|\leq 2r} \cL^t_\beta Z^\beta$ where $\cL^t_\beta$
is a linear combination of operators of the form 
\begin{align*} 
Q (\lambda-\cL_t)^{-k'_0}R_1(\lambda-\cL_t)^{-k'_1}R_2
\cdots R_{l'}(\lambda-\cL_t)^{-k'_{l'}} Q''' Q'',
\end{align*}
with  $R_1, \cdots, R_{l'} \in \cR'_{t}$,  $Q'''\in \mQ^l$,
 $Q''\in   \mQ^m$, $|\beta|\leq 2 r$,
 and $Q''$ is obtained from $Q'$ and its commutation with $Z^\beta$.
By the argument as in \eqref{ue18} and  \eqref{ue20}, as $k>2(m+r+1)$,
we can split the above operator to two parts 
\begin{align*} 
&Q (\lambda-\cL_t)^{-k'_0}R_1(\lambda-\cL_t)^{-k'_1}R_2
\cdots R_{i}(\lambda-\cL_t)^{-k''_{i}};\\
&(\lambda-\cL_t)^{-(k'_{i}-k''_{i}) }\cdots  
R_{l'}(\lambda-\cL_t)^{-k'_{l'}} Q''' Q'',
\end{align*}
and the $\|\quad \|^{0,0}_t$-norm of each part is bounded
by $C$ for $\lambda\in \delta$. Thus
the proof of  (\ref{1ue30}) is complete.

\comment{
We claim that $A ^{\bk}_{\br}(\lambda,t)$ is well defined and 
for $\lambda \in \delta$, $l\in \bN$,
\begin{equation}\label{1ue30}
\|A ^{\bk}_{\br}(\lambda,t) s\|_{t,l}
\leqslant C \sum_{|\alpha|\leqslant 2r} \|Z^\alpha s\|_{t,l+2r-k}\,.
\end{equation}
In fact, by \eqref{0c37}, $\frac{\partial^{r}}{\partial t^{r}}\cL_t$ is
combination of $\frac{\partial^{r_1}}{\partial t^{r_1}}(g^{TX_0}_{ij}(tZ))
(\frac{\partial^{r_2}}{\partial t^{r_2}}\nabla_{t,e_i})
(\frac{\partial^{r_3}}{\partial t^{r_3}}\nabla_{t,e_j})$,
 $\frac{\partial^{r_1}}{\partial t^{r_1}}(b(tZ))$,
$\frac{\partial^{r_1}}{\partial t^{r_1}}(a_{i}(tZ))
(\frac{\partial^{r_2}}{\partial t^{r_2}}\nabla_{t,e_i})$.
Now $\frac{\partial^{r_1}}{\partial t^{r_1}}(b(tZ))$
(resp. $\frac{\partial^{r_1}}{\partial t^{r_1}}\nabla_{t,e_i}$), 
for $r_1\geqslant1$, are functions of the type $b'(tZ)Z^\beta$, 
$|\beta|\leqslant r_1$ (resp. $r_1+1$) and $b'(Z)$ and its derivatives 
with respect to $Z$ are  bounded smooth functions.
In view \eqref{ue12}, we get \eqref{1ue30}.
}

\noindent
By \eqref{ue28}, \eqref{ue30} and the above argument, we get the
estimate \eqref{ue15} with $m'=0$.

Finally, for any vector $U$ on $X$,
\begin{align}\label{0ue30}
& \nabla ^{\pi ^* \End (E)}_U \mP_{q,t} =\frac{1}{2\pi i}
\begin{pmatrix} q+k-1 \\ k-1 \end{pmatrix} ^{-1} \int_{\delta}
\lambda ^{q+k-1} \nabla ^{\pi ^* \End (E)}_U (\lambda-\cL_t)^{-k}
d \lambda .
\end{align}
Now we use a similar formula as \eqref{ue30} for
$\nabla ^{\pi ^* \End (E)} _U(\lambda-\cL_t)^{-k}$ by replacing
$\frac{\partial^{r_1}\cL_t}{\partial t^{r_1}}$ by
$\nabla ^{\pi ^* \End (E)}_U \cL_t$, and remark that
$\nabla ^{\pi ^* \End (E)}_U \cL_t$ is a differential operator
on $T_{x_0} X$ with the same structure  as $\cL_t$.
Then by the above argument, we conclude that \eqref{ue15} holds for $m'\geqslant1$.
\end{proof}

For $k$ big enough, set
\begin{equation}\label{ue31}
\begin{split}
& F_{q,r}= \frac{1}{2\pi i \, r! } 
\begin{pmatrix} q+k-1 \\ k-1 \end{pmatrix}^{-1}  \int_{\delta}
\lambda ^{q+k-1}   \sum_{(\bk,\br)\in I_{k,r} }
 a ^{\bk}_{\br}  A ^{\bk}_{\br}  (\lambda,0)d \lambda ,\\
&F_{q,r,t} = \frac{1}{r!}\frac{\partial^{r}}{\partial t^{r}}
\mP_{q,t}- F_{q,r}.
\end{split}
\end{equation}
 Let $F_{q,r}(Z,Z')$ $(Z,Z'\in T_{x_0}X)$ be the
smooth kernel of $F_{q,r}$  with respect to $dv_{TX}(Z')$.
Then $F_{q,r}\in \cC^\infty (TX\times_{X}TX,\pi ^* \End (E))$.
Certainly, as $t\to0$, the limit of  $\|\quad\|_{t,m}$ exists,
and we denote it by  $\|\quad\|_{0,m}$.
\begin{thm} \label{tue9} For any $r\geq 0$, $k>0$, 
there exists $C>0$ such that for
$t \in [0,1], \lambda \in \delta$,
\begin{align}\label{ue32}
&  \left \|\Big(\frac{\partial^{r}\cL_t}{\partial t^{r}} -
\frac{\partial^{r}\cL_t}{\partial t^{r}} |_{t=0}\Big )s \right \|_{t,-1}
\leqslant Ct \sum_{|\alpha|\leqslant r+3} \|Z^\alpha s\|_{0,1},\\
& \Big \|\Big (\frac{\partial^{r}}{\partial t^{r}} (\lambda-\cL_t)^{-k}
-\sum_{(\bk,\br)\in I_{k,r} }
 a ^{\bk}_{\br}  A ^{\bk}_{\br} (\lambda,0)\Big )s\Big \|_{0,0}
\leqslant C t \sum_{|\alpha|\leqslant 4r+3} \|Z^\alpha s\|_{0,0}.\nonumber
\end{align}
\end{thm}
\begin{proof}  Note that by (\ref{0c36}), (\ref{u0}), for $t\in [0,1]$, $k\geqslant1$,
\begin{align}\label{1ue35}
\|s\|_{t,0}= \|s \|_{0,0},\quad 
\|s\|_{t,k}\leqslant C \sum_{|\alpha|\leqslant k} \|Z^\alpha s\|_{0,k}.
\end{align}
 An application of Taylor expansion for
(\ref{0c37}) leads to the following estimate for
compactly supported $s,s'$:
\be\label{ue33}
\Big | \left \langle  \Big (\frac{\partial^{r}\cL_t}{\partial t^{r}} -
\frac{\partial^{r}\cL_t}{\partial t^{r}} |_{t=0} \Big )s,
s'\right \rangle_{0,0}\Big |
\leqslant C t \|s'\|_{t,1}\sum_{|\alpha|\leqslant r+3} \|Z^\alpha s\|_{0,1}.
\ee
Thus we get the first inequality of (\ref{ue32}).  Note that
\begin{align}\label{ue34}
& (\lambda-\cL_t)^{-1}- (\lambda-\cL_0)^{-1}
=(\lambda-\cL_t)^{-1}(\cL_t-\cL_0) (\lambda-\cL_0)^{-1}.
\end{align}
Now from (\ref{ue33}) and (\ref{ue34}),
\begin{align}\label{0ue34}
& \left \|\left ((\lambda-\cL_t)^{-1}- (\lambda-\cL_0)^{-1} \right)s
 \right \|_{0,0}
\leqslant Ct   \sum_{|\alpha|\leqslant 3} \|Z^\alpha s\|_{0,1}.
\end{align}
After taking the limit, we know that Theorems \ref{tu4}, \ref{tu5} still hold 
for $t=0$. Note that 
$\nabla_{0,e_j} = \nabla_{e_j}+\frac{1}{2} R^L_{x_0}(\mR, e_j)$ 
by (\ref{0c36}). If we denote by $\cL_{\lambda,t}=\lambda -\cL_{t}$, then
\begin{multline}\label{0ue35}
 A^{\bk}_{\br}  (\lambda,t)-  A^{\bk}_{\br}  (\lambda,0)
 = \sum_{i=1}^j \cL_{\lambda,t}^{-k_0}\cdots 
\left ( \frac{\partial^{r_i}\cL_t}{\partial t^{r_i}} 
- \frac{\partial^{r_i}\cL_t}{\partial t^{r_i}} |_{t=0}\right )
\cL_{\lambda,0}^{-k_{i}}\cdots \cL_{\lambda,0}^{-k_j}\\
+ \sum_{i=0}^j \cL_{\lambda,t}^{-k_0}\cdots 
\left (\cL_{\lambda,t}^{-k_i}- \cL_{\lambda,0}^{-k_i} \right)
\left (\frac{\partial^{r_{i+1}}\cL_t}{\partial t^{r_{i+1}}}|_{t=0}\right)
\cdots \cL_{\lambda,0}^{-k_j}.
\end{multline}
From the discussion after (\ref{1ue30}), 
formulas (\ref{ue2}), (\ref{ue30})
and (\ref{0ue34}),  we get (\ref{ue32}).
\end{proof}

\begin{thm} \label{tue12}
For $\sigma >0$, there  exists $C>0$
such that for $t \in ]0,1]$,
 $Z,Z'\in T_{x_0}X$, $|Z|, |Z'|\leqslant \sigma  $,
\begin{align}\label{ue42}
\Big |F_{q,r,t}(Z,Z')\Big |\leqslant & C t^{1/2(2n+1)}.
\end{align}
\end{thm}
\begin{proof}
  By (\ref{ue28}), (\ref{ue31}) and (\ref{ue32}), there exists  $C>0$
such that for $t \in ]0,1]$,
\begin{align}\label{ue39}
& \|F_{q,r,t}\|_{(\sigma),0 }
\leqslant C t .
\end{align}

 Let $\phi : \bR\to [0,1]$ be a smooth function with compact 
support, equal $1$ near $0$, such that  
$\int_{T_{x_0}X} \phi (Z) dv_{TX}(Z)=1$. Take $\nu \in ]0,1]$.
 By the proof of Theorem \ref{tue8} and (\ref{ue31}),
 $F_{q,r}$ verifies the similar inequality as in (\ref{ue15}) 
with $r=0$. Thus by (\ref{ue15}), there exists $C>0$ such that if
$|Z|,|Z'|\leqslant \sigma$, $U,U'\in E_{x_0}$,
\begin{multline}\label{ue43}
\Big |  \left \langle
F_{q,r,t} (Z,Z') U,U'     \right \rangle
-\int_{T_{x_0}X \times T_{x_0}X}  \left \langle
F_{q,r,t}(Z-W,Z'-W') U,U'     \right \rangle\\
  \frac{1}{\nu ^{4n}} \phi (W/\nu) \phi (W'/\nu) dv_{TX}(W)dv_{TX}(W')\Big |
\leqslant C \nu |U||U'|.
\end{multline}
On the other hand, by (\ref{ue39}),
\begin{multline}\label{ue44}
\Big |\int_{T_{x_0}X \times T_{x_0}X}  \left \langle
F_{q,r,t}  (Z-W,Z'-W') U,U'     \right \rangle\\
  \frac{1}{\nu ^{4n}} \phi (W/\nu) \phi (W'/\nu) dv_{TX}(W)dv_{TX}(W')\Big |
\leqslant C       t\frac{1}{\nu ^{2n}}  |U||U'|.
\end{multline}
By taking $\nu = t^{1/2(2n+1)}$, we obtain (\ref{ue42}).
\end{proof}
Finally, we obtain the following off-diagonal estimate for the kernel of 
$\mP_{q,t}$.
\begin{thm} \label{tue14} For $k,m,m'\in \bN$, $\sigma>0$, there exists
$C>0$  such that if
$t\in ]0,1]$, $Z,Z'\in T_{x_0}X$, $|Z|,|Z'| \leqslant \sigma$,
 \begin{align}\label{0ue45}
&\sup_{|\alpha|,|\alpha'|\leqslant m} \Big
|\frac{\partial^{|\alpha|+|\alpha'|}} {\partial Z^{\alpha}
{\partial Z'}^{\alpha'}}\Big (\mP_{q,t} - \sum_{r=0}^k
F_{q,r}t^r\Big ) (Z,Z') \Big |_{\cC ^{m'}(X)} \leqslant C  t^{k+1}.
\end{align}
\end{thm}
\begin{proof} By (\ref{ue31}), and (\ref{ue42}),
 \begin{align}\label{0ue47}
&\frac{1}{r!}  \frac{\partial^{r}}{\partial t^{r}} \mP_{q,t}
|_{t=0} = F_{q,r}.
\end{align}
Now by Theorem \ref{tue8} and (\ref{ue31}), $F_{q,r}$ has the same
estimate as $\frac{1}{r!}\frac{\partial^{r}}{\partial
t^{r}}\mP_{q,t}$, in (\ref{ue15}). Again from (\ref{ue15}), (\ref{ue31}),
 and the Taylor expansion
 \begin{align}\label{0ue46}
G(t)- \sum_{r=0}^k \frac{1}{r !} \frac{\partial ^r G}{\partial t^r}(0) t^r
= \frac{1}{k!}\int_0^t (t-t_0)^k 
\frac{\partial ^{k+1} G}{\partial t^{k+1} }(t_0) dt_0,
\end{align}
we have (\ref{0ue45}).
\end{proof}

\subsection{Bergman kernel of $\cL_0$}
\label{s3.4}
The almost  complex structure $J$ induces a splitting
$T_\bR X\otimes_\bR \bC=T^{(1,0)}X\oplus T^{(0,1)}X$,
where $T^{(1,0)}X$ and $T^{(0,1)}X$ are the eigenbundles of $J$ 
corresponding to the eigenvalues $\sqrt{-1}$ and $-\sqrt{-1}$ respectively. 
We denote by $\det_{\bC}$ the determinant function on the complex
 bundle $T^{(1,0)}X$. 
Set 
\begin{align}\label{0ue62}
\mathcal{J}= -2\pi \sqrt{-1} {\bf J}. 
\end{align}
By \eqref{0.1}, $\mathcal{J}\in \End (T^{(1,0)}X)$
is positive, and $\mathcal{J}$ acting on $TX$ is  skew-adjoint.
For any tensor $\psi$ on $X$, we denote by $\nabla ^{X}\psi$ 
the covariant derivative of $\psi$ induced by $\nabla ^{TX}$.
 Thus  $\nabla ^{X} \mJ, \nabla ^{X} J \in
T^*X \otimes \End(TX)$, $\nabla ^{X} \nabla ^{X} \mJ\in T^*X
\otimes T^*X \otimes \End(TX)$. We also adopt
the convention that all tensors will be evaluated at the base point
$x_0\in X$, and most of the time, we will omit the subscript
$x_0$. 

Let $P^N$ be the orthogonal projection from  
$(L^2 (\bR^{2n}, E_{x_0}), \|\quad\|_0=\|\quad\|_{t,0})$ onto $N=\Ker \cL_0$,
and let $P^N(Z,Z')$ be the smooth kernel of $P^N$ with respect 
to $dv_{T_{x_0}X}(Z)$. Then  $P^N(Z,Z')$ is the
the Bergman kernel of $\cL_0$. For $Z,Z'\in T_{x_0}X$, we have
\begin{align}\label{ue62}
P^N(Z,Z')  =\frac{\det_{\bC}{\mathcal{J}_{x_0}}}{(2\pi)^n}
\exp\Big (- \frac{1}{4} 
\left \langle (\mathcal{J}^2_{x_0})^{1/2}(Z-Z'),(Z-Z') \right \rangle
+\frac{1}{2}    \left \langle \mathcal{J}_{x_0} Z,Z' \right \rangle\Big ).
\end{align}

Now we discuss the eigenvalues and eigenfunctions of $\cL_0$ 
in more precise way. We choose $\{ w_i\}_{i=1}^n$
an orthonormal basis of $T^{(1,0)}_{x_0} X$, such that
\begin{align}\label{0ue52}
\mathcal{J}_{x_0}= {\rm diag} (a_1 , \cdots, a_n)\in {\rm End}
(T^{(1,0)}_{x_0} X),
\end{align}
with $0< a_1\leqslant a_2\leqslant \cdots \leqslant  a_n$, and let
$\{w^j\}_{j=1}^n$ be its dual basis. Then
$e_{2j-1}=\tfrac{1}{\sqrt{2}}(w_j+\overline{w}_j)$ and
$e_{2j}=\tfrac{\sqrt{-1}}{\sqrt{2}}(w_j-\overline{w}_j)\,,
 j=1,\dotsc,n\, $
forms an orthonormal basis of $T_{x_0}X$.
We use the coordinates on $T_{x_0}X\simeq\bR^{2n}$ induced by $\{ e_i\}$ as in \eqref{0c11}
and in what follows we also introduce the complex coordinates $z=(z_1,\cdots,z_n)$
on $\bC^n\simeq\bR^{2n}$. Thus
$Z=z+\overline{z}$, and $w_i=\sqrt{2}\tfrac{\partial}{\partial z_i}$,
$\overline{w}_i=\sqrt{2}\tfrac{\partial}{\partial\overline{z}_i}$.
We will also identify $z$ to $\sum_i z_i\tfrac{\partial}{\partial z_i}$ 
and $\overline{z}$ to
$\sum_i\overline{z}_i\tfrac{\partial}{\partial\overline{z}_i}$ when
we consider $z$ and $\overline{z}$ as vector fields. Remark that
\begin{equation}\label{g0}
\Big\lvert\tfrac{\partial}{\partial z_i}\Big\rvert^2=
\Big\lvert\tfrac{\partial}{\partial\overline{z}_i}\Big\rvert^2
=\dfrac{1}{2}\,,
\quad\text{so that $|z|^2=|\overline{z}|^2=\dfrac{1}{2} |Z|^2$\,.}
\end{equation}
It is very useful to rewrite $\cL_0$ by using  
the creation and annihilation operators. Set
\begin{equation}\label{0g1}
\nabla_{0,\cdot}= \nabla_{\cdot} + \frac{1}{2} R^L_{x_0}(\mR, \cdot),
\quad b_i=-2\nabla_{0,\tfrac{\partial}{\partial z_i}},\quad
b^{+}_i=2\nabla_{0,\tfrac{\partial}{\partial \overline{z}_i}},
\quad b=(b_1,\cdots,b_n)\,.
\end{equation}
Then by \eqref{0ue62}, and \eqref{0ue52}, we have 
\begin{equation}\label{g1}
b_i=-2{\tfrac{\partial}{\partial z
_i}}+\frac{1}{2}a_i\overline{z}_i\,,\quad
b^{+}_i=2{\tfrac{\partial}{\partial\overline{z}_i}}+\frac{1}{2}a_i
z_i,
\end{equation}
and for any polynomial  $g(z,\ov{z})$ on $z$ and $\ov{z}$, 
\begin{align}\label{g2}
&[b_i,b^{+}_j]=b_i b^{+}_j-b^{+}_j b_i =-2a_i \delta_{i\,j},\\
&[b_i,b_j]=[b^{+}_i,b^{+}_j]=0\, ,\nonumber\\
& [g(z,\ov{z}),b_j]=  2 \tfrac{\partial}{\partial z_j}g(z,\ov{z}), 
\quad  [g(z,\ov{z}),b_j^+]
= - 2\tfrac{\partial}{\partial \overline{z}_j}g(z,\ov{z})\,. \nonumber
\end{align}
By  \eqref{0.2},  \eqref{0ue52}, $\tau_{x_0}= \sum_i a_i$. 
Thus from  \eqref{c31}, \eqref{0ue52}, \eqref{0g1}-\eqref{g2},
\begin{equation}\label{g3}
\cL_0=\sum_i b_i b^{+}_i.
\end{equation}
\begin{rem} \label{r3.4} Let $L=\bC$ be the trivial holomorphic line bundle
 on $\bC^n$ with the canonical section $1$. Let $h^L$ be the metric on $L$
 defined by  $|1|_{h^L}(z) = e ^{-\frac{1}{4}\sum_{j=1}^n a_{j}|z_{j}|^2}$
 for $z\in \bC^n$. Let $g^{T\bC ^n}$ be the canonical metric on $\bC^n$. 
Then $\cL_{0}$ is twice the corresponding Kodaira-Laplacian 
$\overline{\partial}^{L*}\overline{\partial}^{L}$ under the trivialization 
of $L$ by using the unit section
  $e ^{\frac{1}{4}\sum_{j=1}^n a_{j}|z_{j}|^2} 1$.
\end{rem}
\begin{thm}\label{t3.4}
The spectrum of the restriction of $\cL_0$ on $L^2(\bR^{2n})$ is given by
\begin{equation}\label{g4}
\spec{{\cL_0}_{\upharpoonright\,{L^2(\bR^{2n})}}}=
\Big\lbrace2\sum_{i=1}^n\alpha_i a_i\,:\, 
 \alpha =(\alpha_1,\cdots,\alpha_n)\in\bN ^n\Big\rbrace
\end{equation}
and an orthogonal basis of the eigenspace of $2\sum_{i=1}^n\alpha_i a_i$
is given by
\begin{equation}\label{g5}
b^{\alpha}\big(z^{\beta}\exp\big({-    \frac{1}{4}\sum_i
a_i|z_i|^2}\big)\big)\,,\quad\text{with $\beta\in\bN^n$}\,.
\end{equation}
\end{thm}
\begin{proof} At first $z^{\beta}\exp\big({-    \frac{1}{4}\sum_i
a_i|z_i|^2}\big)$, $\beta\in\bN^n$ are annihilated by $b^{+}_i$ 
$(1\leq i\leq n)$, thus they are in the kernel of 
${\cL_0}_{\upharpoonright\,{L^2(\bR^{2n})}}$. 
Now, by  (\ref{g2}), \eqref{g5} are eigenfunctions of 
${\cL_0}_{\upharpoonright\,{L^2(\bR^{2n})}}$ with eigenvalue 
$2\sum_{i=1}^n\alpha_i a_i$. 
But the span of functions \eqref{g5} includes 
all the rescaled Hermite polynomials 
multiplied by $\exp\big({- \frac{1}{4}\sum_i a_i|z_i|^2}\big)$, 
which is an orthogonal basis of $L^2(\bR^{2n})$ by \cite[\S 6]{T2}.  
Thus the  eigenfunctions in \eqref{g5} 
are all the eigenfunctions of ${\cL_0}_{\upharpoonright\,{L^2(\bR^{2n})}}$.
The proof of Theorem \ref{t3.4} is complete.
\end{proof}

Especially an orthonormal basis of 
$\Ker \cL_0\upharpoonright_{L^2(\bR^{2n})}$ is
\begin{align}\label{g6}
\Big(\frac{a ^\beta}{(2\pi)^n 2 ^{|\beta|} \beta!}
\prod_{i=1}^n a_i\Big)^{1/2}z^\beta
\exp\Big (-\frac{1}{4} \sum_{j=1}^n a_j |z_j|^2\Big )\,,\quad \beta\in\bN^n\,.
\end{align}
From \eqref{g6}, we recover \eqref{ue62}:
\begin{equation}\label{g7}
\begin{split}
P^N(Z,Z') &=\frac{1}{(2\pi)^n}\prod_{i=1}^n
a_i\:\:\exp\Big(-\frac{1}{4}\sum_i
a_i\big(|z_i|^2+|z^{\prime}_i|^2 -2z_i\overline{z}_i'\big)\Big).
\end{split}
\end{equation}

Recall that the operators $\mO_1,\mO_2$ were defined in (\ref{c31}). 
 Theorem \ref{t3.5} below is crucial in proving the vanishing result of
 $F_{q,r}$ (cf. Theorem \ref{0t3.6}).
\begin{thm}\label{t3.5}
We have the relation
\be\label{g8} P^N\mO_1 P^N =0.\ee
\end{thm}
\begin{proof} From \eqref{0.1}, for $U,V,W\in TX$, 
$\langle(\nabla ^{X}_U{\bf J})V,W\rangle=(\nabla ^{X}_U\omega)(V,W)$, thus
\begin{equation} \label{g9}
\langle(\nabla ^{X}_U{\bf J})V,W\rangle
+\langle(\nabla ^{X}_V{\bf J})W,U\rangle+\langle
(\nabla ^{X}_W{\bf J})U,V\rangle=d\omega(U,V,W)=0.
\end{equation}
By \eqref{0.1} and \eqref{0.2}, 
\begin{equation}\label{g10}
\begin{split}
&R^L(U,V)=\langle \mJ U,V \rangle,   \\
& (\nabla ^{X}_U R^L) (V,W)
= \langle (\nabla ^{X}_U\mJ) V,W \rangle,\\
&\nabla _U\tau =  
- \frac{\sqrt{-1}}{2}\tr|_{TX}[\nabla ^{X}_U(J\mJ)]. 
\end{split}
\end{equation}

 As $J$, $\mathcal{J}\in \End (TX)$ are
skew-adjoint and commute, $\nabla ^{X}_{U}J$, $\nabla ^{X}_{U}\mathcal{J}$
are skew-adjoint and $\nabla ^{X}_{U}(J\mathcal{J})$ is symmetric.
From $J^2=-\Id$, we know that 
\begin{align}\label{g11}
J(\nabla ^{X} J)+(\nabla ^{X} J)J=0,
\end{align}
thus $\nabla ^{X}_{U}J$ exchanges $T^{(1,0)}X$ and $T^{(0,1)}X$.
 From \eqref{g9} and \eqref{g10}, we have
\begin{equation}\label{g12}
\begin{split}
&(\nabla _\mR\tau)_{x_0}=-2\sqrt{-1} \left\langle
(\nabla ^{X}_{\mR}(J\mJ))\tfrac{\partial}{\partial z_i},
\tfrac{\partial}{\partial\ov{z}_i}\right\rangle
=2\left\langle(\nabla ^{X}_{\mR}\mJ)\tfrac{\partial}{\partial
z_i},
\tfrac{\partial}{\partial\ov{z}_i}\right\rangle,\\
&(\partial_i R^L)_{x_0}(\mR,e_i)
=2\Big\langle(\nabla ^{X}_{\tfrac{\partial}{\partial
z_i}}\mJ) \mR,\tfrac{\partial}{\partial\ov{z}_i}\Big\rangle
+ 2 \Big\langle(\nabla ^{X}_{\tfrac{\partial}{\partial\overline{z}_i}}\mJ)
\mR,\tfrac{\partial}{\partial z_i}\Big\rangle \\
&\hspace*{20mm}= 4\Big\langle(\nabla ^{X}_{\tfrac{\partial}{\partial
z_i}}\mJ) \mR,\tfrac{\partial}{\partial\ov{z}_i}\Big\rangle 
-2 \left\langle(\nabla ^{X}_{\mR}\mJ) \tfrac{\partial}{\partial
z_i},\tfrac{\partial}{\partial\ov{z}_i}\right\rangle . 
\end{split}
\end{equation}
 From \eqref{c31}, \eqref{g2},  \eqref{g10} and \eqref{g12}, we infer
\begin{equation} \label{g13}
\begin{split}
\mO_1&=-\frac{2}{3}\Big[\left\langle(\nabla ^{X}_\mR\mJ)\mR,
\tfrac{\partial} {\partial z_i}\right\rangle b^+_i
- \left\langle(\nabla ^{X}_\mR\mJ)\mR,\tfrac{\partial}
{\partial\ov{z}_i}\right\rangle b_i\\
&\hspace*{10mm}+2\Big\langle(\nabla ^{X}_{\tfrac{\partial}{\partial
z_i}}\mJ)\mR,\tfrac{\partial}{\partial\ov{z}_i}\Big\rangle+
2\left\langle(\nabla ^{X}_\mR\mJ)\tfrac{\partial}{\partial
z_i},\tfrac{\partial} {\partial\ov{z}_i}\right\rangle\Big]\\
&=-\frac{2}{3}\Big[\left\langle(\nabla ^{X}_\mR\mJ)\mR,
\tfrac{\partial} {\partial z_i}\right\rangle b^+_i
-b_i \left\langle(\nabla ^{X}_\mR\mJ)\mR,\tfrac{\partial}
{\partial\ov{z}_i}\right\rangle\Big].
\end{split}
\end{equation}
 Note that by (\ref{g1}),  (\ref{g7}), 
\begin{equation} \label{g14}
(b^+_iP^N)(Z,Z^{\prime}) =0\,,\quad \,
(b_iP^N)(Z,Z^{\prime})=a_i(\ov{z}_i-\ov{z}^{\prime}_i)P^N(Z,Z^{\prime}).
\end{equation}
We learn from \eqref{g14} that for any polynomial $g(z,\ov{z})$ in $z,\ov{z}$ 
we can write $g(z,\ov{z})P^N(Z,Z^{\prime})$ as 
sums of $b^\beta g_{\beta} (z,\ov{z}')P^N(Z,Z^{\prime})$ with 
$g_{\beta} (z,\ov{z}')$ polynomials in $z,\ov{z}'$. 
By Theorem \ref{t3.4},
\begin{equation}\label{g15}
 P^N b^\alpha   g(z,\ov{z})P^N=0\,,\quad\text{for $|\alpha|>0$},
\end{equation}
and relations \eqref{g13}\,--\,\eqref{g15} yield the desired relation \eqref{g8}.
\end{proof}

\subsection{Evaluation of $F_{q,r}$}
\label{s3.5}
For $s\in \bR$, let $[s]$ denote the greatest integer 
which is less than or equal to $s$.
Let $f(\lambda, t)$ be a formal power series with values  in $
\End(L^2(\bR^{2n},E_{x_0}))$ \be\label{c54} f(\lambda, t)=
\sum_{r=0}^\infty t^r f_r(\lambda), \quad f_r(\lambda)\in
\End(L^2(\bR^{2n}, E_{x_0})). \ee
 By (\ref{c30}),  consider
the equation of formal power series for  $\lambda\in \delta$,
\begin{align}\label{c55}
(-\cL_0 +\lambda - \sum_{r=1}^\infty t^r \mO_r) f(\lambda, t) =
\Id_{L^2(\bR^{2n}, E_{x_0})}.
\end{align}
Let $N^\bot$ be the orthogonal space of $N$ in
$L^2(\bR^{2n},E_{x_0})$, and $P^{N^\bot}$ be the orthogonal
projection from $L^2(\bR^{2n},E_{x_0})$ to $N^\bot$. We decompose
$f(\lambda, t)$ according the splitting $L^2(\bR^{2n}, E_{x_0})=
N\oplus N^\bot$,
\begin{align}\label{c56}
 g_r(\lambda) = P^{N}f_r(\lambda),\quad
 f^\bot_r(\lambda)=P^{N^\bot}f_r(\lambda).
\end{align}
Using (\ref{c56}) and identifying the powers of $t$ in
(\ref{c55}), we find that
\begin{equation}\label{c57}
\begin{split}
&g_0(\lambda) = \frac{1}{\lambda}P^{N},
\quad  f^\bot_0(\lambda)= (\lambda -\cL_0)^{-1}P^{N^\bot} ,\\
&f^\bot_r(\lambda)=(\lambda -\cL_0)^{-1}
 \sum_{j=1}^{r} P^{N^\bot}\mO_j   f_{r-j}(\lambda),\\
&  g_r(\lambda)  = \frac{1}{\lambda } \sum_{j=1}^{r} P^{N}\mO_j
f_{r-j}(\lambda). 
\end{split}
\end{equation}
\begin{lemma}\label{t3.6}  For $r\in \bN$,
$\lambda ^{[\frac{r}{2}] +1} g_r(\lambda)$, $\lambda
^{[\frac{r+1}{2}] }f^\bot_r(\lambda)$ are holomorphic functions
for $|\lambda|\leqslant \mu_0/4$ and
\begin{equation}\label{c81}
(\lambda ^{r+1} g_{2r})(0) =( P^N\mO_2 P^N
- P^N \mO_1 \cL_0^{-1} P^{N^\bot} \mO_1 P^N)^r P^N.
\end{equation}
\end{lemma}
\begin{proof}
By \eqref{c57} we know that Lemma \ref{t3.6} is true for $r=0$.
Assume that Lemma \ref{t3.6} is true for $r\leqslant m$. Now, by
Theorem \ref{t3.4}, (\ref{c57}) and the recurrence assumption, it follows that $\lambda
^{[\frac{m}{2}]+1 }f^\bot_{m+1}(\lambda)$ is holomorphic  for
$|\lambda|\leqslant \mu_0/4$, and
\begin{equation}\label{c82}
\lambda ^{[\frac{m+1}{2}] +1} g_{m+1}(\lambda) 
=  \lambda ^{[\frac{m+1}{2}]} \sum_{i=1}^{m+1} P^N \mO_i\Big [  
g_{m+1-i}(\lambda)  + f^\bot_{m+1-i}(\lambda) \Big ] .
\end{equation}
 By our recurrence, $\lambda ^{[\frac{m+1}{2}]-1}
f^\bot_{m+1-j}(\lambda)$, $\lambda ^{[\frac{m+1}{2}]-1}
g_{m-j}(\lambda)$, $\lambda ^{[\frac{m+1}{2}]}
f^\bot_{m}(\lambda)$ are
 holomorphic  for $|\lambda|\leqslant \mu_0/4$, $j\geqslant2$.
Thus by Theorem \ref{t3.5}, and \eqref{c57} and \eqref{c82},
$\lambda ^{[\frac{m+1}{2}] +1} g_{m+1}(\lambda)$
 is also  holomorphic  for $|\lambda|\leqslant \mu_0/4$, and
\begin{multline}\label{c84}
(\lambda ^{[\frac{m+1}{2}] +1} g_{m+1})(0)
=\Big( \lambda ^{[\frac{m+1}{2}]}
 (P^N \mO_1  f^\bot_{m}  +P^N\mO_2 g_{m-1}  ) \Big)(0)\\
= -P^N \mO_1 \cL_0^{-1} P^{N^\bot} \mO_1 
\Big( \lambda ^{[\frac{m+1}{2}]}( g_{m-1}
 + f^\bot_{m-1})\Big)(0)\
+( \lambda ^{[\frac{m+1}{2}]} P^N\mO_2 g_{m-1})(0) . 
\end{multline}
If $m$ is odd, then by (\ref{c84}) and recurrence assumption,
\begin{multline}\label{c85}
(\lambda ^{[\frac{m+1}{2}] +1} g_{m+1})(0)
=P^N (-\mO_1 \cL_0^{-1} P^{N^\bot} \mO_1 + \mO_2
)P^N(\lambda ^{[\frac{m-1}{2}] +1} g_{m-1})(0) \\
= ( P^N\mO_2
P^N- P^N \mO_1 \cL_0^{-1} P^{N^\bot} \mO_1 P^N)^{[\frac{m+1}{2}]}P^N .
\end{multline}
The proof of Lemma \ref{t3.6} is complete.
\end{proof}

\begin{thm} \label{0t3.6} There exist  $J_{q,r}(Z,Z')$  polynomials
in $Z,Z'$ with the same parity as $r$ and $\deg J_{q,r}(Z,Z')\leq 3r$, 
whose coefficients  
are polynomials in $R^{TX}$, $R^E$ {\rm(}\,and $R^L$, $\Phi${\rm)}  
and their derivatives of order $\leqslant r-1$ {\rm(}\,resp. $r${\rm)}, 
and reciprocals of linear combinations of eigenvalues of ${\bf J}$ 
 at $x_0$\,, such that
\begin{align}\label{c86}
&F_{q,r}(Z,Z')= J_{q,r}(Z,Z')P^N(Z,Z').
\end{align}
Moreover,
\begin{equation}\label{1c52}
\begin{split}
&F_{0,0}=P^N,\\ 
&F_{q,r} =0, \quad\text{for $q>0, r< 2q$}\,,   \\
&F_{q,2q} =\big(P^N\mO_2P^N
- P^N \mO_1 \cL_0^{-1} P^{N^\bot} \mO_1 P^N\big)^q\, 
P^N\quad\text{for $q>0$}\,.
\end{split}
\end{equation}
\end{thm}
\begin{proof} Recall that $\mP_{q,t}= (\cL_t)^q \mP_{0,t}$.
By  (\ref{1ue15}),
$\mP_{q,t} = \frac{1}{2 \pi i} \int_{\delta}\lambda ^q (\lambda
-\cL_t)^{-1}d\lambda$ .
Thus by (\ref{ue28}), (\ref{ue31}) and (\ref{c56}),
\begin{equation}\label{c90}
F_{q,r}= \frac{1}{2\pi i} \int_{\delta}\lambda ^q  g_r (\lambda)d
\lambda +  \frac{1}{2\pi i} \int_{\delta}\lambda ^q f_r^\bot
(\lambda)d \lambda.
\end{equation}
From  Lemma \ref{t3.6} and  (\ref{c90}), we get  (\ref{1c52}).
Generally,  from Theorems \ref{t3.3}, \ref{t3.4}, Remark \ref{0t3.3},
\eqref{g7}, \eqref{c57}, \eqref{c90} and the residue formula, we conclude that 
$F_{q,r}$ has the form \eqref{c86}.
\end{proof}

From Theorems \ref{t3.4}, \ref{t3.5}, (\ref{c57}), (\ref{c90}) 
and the residue formula,
we can get $F_{q,r}$ by using the operators $\cL_0^{-1}$, 
$P^N$, $P^{N^\bot}$, $\mO_k (k\leq r)$. 
This give us a direct method to compute  $F_{q,r}$ 
in view of Theorem \ref{t3.4}. In particular, we get
\footnote{The formula $F_{0,2}$ in \cite[(20)]{MM04a} missed the last
two terms here which are zero at $(0,0)$ if ${\bf J}=J$, 
cf. Section \ref{s4.3}.}
\begin{equation}\label{g51}
\begin{split}
 F_{0,1}=& -P^N\mO_1\cL_0^{-1}P^{N^\bot}-P^{N^\bot}\cL_0^{-1}\mO_1P^N,\\
 F_{0,2}=& \frac{1}{2 \pi i} \int_{\delta}
  \Big[(\lambda -\cL_0)^{-1} P^{N^\bot} (\mO_1 f_1 + \mO_2 f_0)(\lambda)
+ \frac{1}{\lambda}  P^N  (\mO_1 f_1+ \mO_2 f_0)(\lambda)\Big]d\lambda \\
=&\cL_0^{-1}P^{N^\bot}\mO_1 \cL_0^{-1}P^{N^\bot}\mO_1 P^N
- \cL_0^{-1}P^{N^\bot}\mO_2 P^N\\
&+ P^N\mO_1\cL_0^{-1}P^{N^\bot}\mO_1\cL_0^{-1}P^{N^\bot}
- P^N \mO_2\cL_0^{-1}P^{N^\bot}\\
&+P^{N^\bot}\cL_0^{-1}\mO_1 P^N\mO_1\cL_0^{-1}P^{N^\bot}
-P^N\mO_1 \cL_0^{-2}  P^{N^\bot}\mO_1  P^N.
\end{split}
\end{equation}

\subsection{Proof of Theorem \ref{t0.1}}
\label{s3.6}
Recall that $P_{0,q,p}= (\Delta_{p,\Phi_0}^{X_0})^q P_{0,\mH_p}$.
 By (\ref{c22}), (\ref{c27}), for $Z,Z'\in
\bR^{2n}$, \be\label{0c53} P_{0,q,p}(Z,Z') =
t^{-2n-2q}\kappa ^{-\frac{1}{2}}(Z)  \mP_{q,t}(Z/t, Z'/t) \kappa ^{-\frac{1}{2}}(Z').
 \ee
By \eqref{1c19}, \eqref{0c53}, Proposition \ref{p3.2},
Theorems \ref{tue14} and \ref{0t3.6},
we get the following main technical result of this paper, the near off-diagonal expansion 
of the generalized Bergman kernels:
\begin{thm} \label{t3.8} For $k,m,m'\in \bN$, $k\geqslant2q$, $\sigma>0$,
 there exists
$C>0$  such that if $p\geqslant1$, $Z,Z'\in T_{x_0}X$, $|Z|,|Z'| \leq
\sigma/\sqrt{p}$,
 \begin{multline}\label{1c53}
\sup_{|\alpha|+|\alpha'|\leqslant m} \Big|
\frac{\partial^{|\alpha|+|\alpha'|}} {\partial Z^{\alpha}
{\partial Z'}^{\alpha'}}\Big ( \frac{1}{p^n}P_{q,p} (Z,Z')\\
 - \sum_{r=2q}^k F_{q,r} (\sqrt{p}Z,\sqrt{p}Z')
\kappa ^{-\frac{1}{2}}(Z)\kappa^{-\frac{1}{2}}(Z')
p^{-\frac{r}{2}+q}\Big ) \Big |_{\cC ^{m'}(X)}
 \leqslant C p^{-\frac{k-m}{2}+q}.
\end{multline}
\end{thm}

\noindent
Set now $Z=Z'=0$ in \eqref{1c53}. By Theorem \ref{0t3.6},
we obtain \eqref{0.6} and
\begin{align}\label{1c54}
b_{q,r}(x_0)= F_{q,2r+2q} (0,0).
\end{align}
Hence \eqref{0.5} follows from \eqref{ue62} and \eqref{1c54}.
The statement about the structure of $b_{q,r}$ 
follows from Theorems \ref{t3.4} and \ref{0t3.6}. 
To prove the uniformity part of Theorem \ref{t0.1}, we notice that 
in the proof of Theorem \ref{tue8}, we only use the derivatives of 
the coefficients of $\cL_t$ with order $\leqslant 2n+2m+m'+2$. 
Thus by \eqref{0ue46},
 the constants in Theorems \ref{tue8}, \ref{tue12} and \ref{tue14} 
are uniformly bounded, if with respect to a fixed metric $g^{TX}_0$, 
the $\cC^{2n+2m+m'+4}$\,--\,norms on $X$ of the data 
{\rm(}$g^{TX}$, $h^L$, $\nabla ^L$, $h^E$, $\nabla ^E$, $J$ and $\Phi${\rm)}
 are bounded, and $g^{TX}$ is bounded below.  
Moreover, taking derivatives with respect to the parameters we obtain 
a similar equation as \eqref{0ue30}, where
$x_0\in X$ plays now a role of a parameter. Thus 
the $\cC^{m'}$--\,norm in \eqref{1c53} can also include the parameters
if the $C^{m'}$--\,norms (with respect to the parameter $x_0\in X$) of the 
derivatives of above data with order $\leqslant 2n+2m+2k+m'+4$ are bounded.
Thus we can take $C_{k,\, l}$ in \eqref{0.6} independent of $g^{TX}$
under our condition.
This achieves the proof of Theorem \ref{t0.1}.

\section{Computing the coefficients $b_{q,r}$}
\label{s4}

In principle, Theorem \ref{t3.4}, the equations \eqref{c57}, \eqref{c90}
and the residue formula give us a direct method to
calculate $b_{q,r}$ by recurrence. Actually, it is computable for
the first few terms $b_{q,r}$ in (\ref{0.6}) in this way.
This Section is organized as follows. 
In Section \ref{s4.1}, we will give a simplified formula
for $\mO_2P^N$ without the assumption ${\bf J}=J$. In Sections  \ref{s4.2}, 
\ref{s4.3},
we will compute $b_{q,0}$ and $b_{0,1}$ under the assumption ${\bf J}=J$,
thus proving Theorem \ref{t0.2}.

In this Section, we use the notation in Section \ref{s3.4}, 
and all tensors will be evaluated at the base point $x_0\in X$.
Recall that the operators $\mO_1,\mO_2$ were defined in (\ref{c31}).

\subsection{A  formula for  $\mO_2P^N$}\label{s4.1}

We will use the following Lemma to evaluate $b_{q,r}$ in (\ref{0.6}).
\begin{lemma}\label{t4.1}
The following relation holds:
\begin{equation}\label{g17}
\begin{split}
&\mO_2P^N = \Big\{\frac{1}{3}b_ib_j 
\left \langle R^{TX} (\mR,\tfrac{\partial}{\partial\ov{z}_i})
\mR, \tfrac{\partial}{\partial\ov{z}_j}\right \rangle
+ \frac{1}{2} b_i \sum_{|\alpha|=2}(\partial ^{\alpha}R^L)_{x_0}
(\mR,\tfrac{\partial}{\partial \ov{z}_i})
\frac{Z^\alpha}{\alpha !} \\
&+  \frac{4}{3} b_j \left[  \left \langle R^{TX} 
(\tfrac{\partial}{\partial z_i},\tfrac{\partial}{\partial\ov{z}_i})
\mR, \tfrac{\partial}{\partial\ov{z}_j}\right \rangle 
- \left \langle R^{TX} (\mR,\tfrac{\partial}{\partial z_i})
\tfrac{\partial}{\partial \ov{z}_i}, 
\tfrac{\partial}{\partial\ov{z}_j}\right \rangle \right]
+R^E (\mR,\tfrac{\partial}{\partial \ov{z}_i})b_i \\
& +
\left \langle (\nabla ^{X} \nabla ^{X}\mJ)_{(\mR,\mR)}
\tfrac{\partial}{\partial z_i},
 \tfrac{\partial}{\partial \ov{z}_i} \right \rangle 
+ 4 \left \langle R^{TX} (\tfrac{\partial}{\partial z_i},
\tfrac{\partial}{\partial z_j}) \tfrac{\partial}{\partial \ov{z}_i}, 
\tfrac{\partial}{\partial \ov{z}_j}\right \rangle \Big\} P^N\\
&+ \Big(-\frac{1}{3}\cL_0 \left \langle R^{TX} 
(\mR,\tfrac{\partial}{\partial z_j})
\mR, \tfrac{\partial}{\partial \ov{z}_j}\right \rangle 
+\frac{1}{9} |(\nabla_\mR^X \mJ) \mR |^2
-\sum_{|\alpha|=2}(\partial ^{\alpha}\tau)_{x_0}
\frac{Z^\alpha}{\alpha !}+  \Phi \Big)  P^N.
\end{split}
\end{equation}
\end{lemma}
\begin{proof}
Set 
\begin{equation}\label{g18}
\begin{split}
I_1& =  \frac{1}{2} \sum_{|\alpha|=2}(\partial ^{\alpha}R^L)_{x_0}
(\mR,\tfrac{\partial}{\partial \ov{z}_i}) 
\frac{Z^\alpha}{\alpha !}b_i\\
&- \frac{1}{2} {\tfrac{\partial}{\partial \ov{z}_i}}
\Big(\sum_{|\alpha|=2}(\partial ^{\alpha}R^L)_{x_0}
(\mR,\tfrac{\partial}{\partial z_i})
\frac{Z^\alpha}{\alpha !}\Big)
-\frac{1}{2}{\tfrac{\partial}{\partial z_i}}
\Big(\sum_{|\alpha|=2}(\partial ^{\alpha}R^L)_{x_0}
(\mR,\tfrac{\partial}{\partial \ov{z}_i})
\frac{Z^\alpha}{\alpha !}\Big), \\
I_2&=\frac{1}{3} \left \langle R^{TX}_{x_0} 
(\mR,\tfrac{\partial}{\partial\ov{z}_i})
\mR, \tfrac{\partial}{\partial\ov{z}_j}\right \rangle b_ib_j\\
&-\frac{4}{3} \left[\left \langle R^{TX}_{x_0} 
(\mR,\tfrac{\partial}{\partial \ov{z}_i} )
\tfrac{\partial}{\partial z_i},
\tfrac{\partial}{\partial \ov{z}_j}\right \rangle
+ \left \langle R^{TX}_{x_0} (\mR,\tfrac{\partial}{\partial z_i} )
\tfrac{\partial}{\partial \ov{z}_i},
\tfrac{\partial}{\partial \ov{z}_j}\right \rangle\right] b_j. 
\end{split}
\end{equation}
From \eqref{c31}, \eqref{0ue62}, \eqref{g1}, \eqref{g2}, \eqref{g10},
\eqref{g18}, and since $\mJ$ is purely imaginary,
\begin{equation}\label{g19}
\begin{split}
\mO_2 &= I_1 +I_2-\frac{1}{3}\Big [\cL_0, \left \langle R^{TX} 
(\mR,\tfrac{\partial}{\partial z_j})
\mR, \tfrac{\partial}{\partial \ov{z}_j}\right \rangle\Big]
+ R^E (\mR,\tfrac{\partial}{\partial \ov{z}_i})b_i\\
&+ \frac{1}{3} \left[
\left \langle R^{TX} (\mR,\tfrac{\partial}{\partial z_i})
\mR, \tfrac{\partial}{\partial  \ov{z}_j}\right \rangle
(-2b_jb_i^+ -2a_i \delta_{ij})
+ \left \langle R^{TX} (\mR,\tfrac{\partial}{\partial z_i})
\mR, \tfrac{\partial}{\partial z_j}\right \rangle b_i^+b_j^+\right]\\
&+\Big( \frac{2}{3} \left \langle R^{TX}
(\mR, e_j) e_j,\tfrac{\partial}{\partial z_i}\right \rangle 
- \Big(\frac{1}{2}\sum_{|\alpha|=2}(\partial ^{\alpha}R^L)_{x_0}
\frac{Z^\alpha}{\alpha !}
+R^E\Big) (\mR,\tfrac{\partial}{\partial z_i})\Big) b_i^+\\
&+\frac{1}{9} |(\nabla_\mR^X \mJ) \mR |^2
-\sum_{|\alpha|=2}(\partial ^{\alpha}\tau)_{x_0}
\frac{Z^\alpha}{\alpha !}+  \Phi.
\end{split}
\end{equation}
In normal coordinates, $(\nabla^{TX}_{e_i}e_j)_{x_0}=0$, 
so from \eqref{0c31}, at $x_{0}$,
\begin{equation}\label{g20}
\begin{split}
\nabla _{e_j} \nabla _{e_i}\left \langle \mJ e_k, e_l \right
\rangle  =& \left \langle (\nabla ^{X}_{e_j} \nabla
^{X}_{e_i} \mJ) e_k
 + \mJ(\nabla ^{TX}_{e_j} \nabla ^{TX}_{e_i}  e_k), e_l \right\rangle 
+  \left \langle  \mJ e_k, 
\nabla ^{TX}_{e_j} \nabla ^{TX}_{e_i}e_l \right \rangle \\
=& \left \langle (\nabla ^{X}_{e_j} \nabla ^{X}_{e_i} \mJ) e_k,
e_l \right \rangle
 -\frac{1}{3}  \left \langle R^{TX}(e_j,e_i)e_k+ R^{TX}(e_j,e_k)e_i, 
\mJ e_l  \right \rangle  \\
&+\frac{1}{3}  \left \langle R^{TX}(e_j,e_i)e_l+
R^{TX}(e_j,e_l)e_i, \mJ e_k \right \rangle .
\end{split}
\end{equation}
From  \eqref{g10}, \eqref{g20},
\begin{equation}\label{g21}
\begin{split}
&\sum_{|\alpha|=2}(\partial ^{\alpha}R^L)_{x_0} (e_k,e_l) 
\frac{Z^\alpha}{\alpha !}
= \frac{1}{2} (\nabla _{e_j} \nabla _{e_i}
\left \langle \mJ e_k, e_l \right \rangle )_{x_0}  Z_iZ_j \\
&=  \frac{1}{2} \left \langle (\nabla ^{X} \nabla ^{X}\mJ)_{(\mR,\mR)} e_k,
e_l \right \rangle  + \frac{1}{6} \left[ \left \langle R^{TX}
(\mR, e_l)\mR, \mJ e_k \right \rangle  - \left \langle R^{TX}
(\mR, e_k)\mR, \mJ e_l \right \rangle \right].
\end{split}
\end{equation}
Thus
\begin{align}\label{g22}
\sum_{|\alpha|=2}(\partial ^{\alpha}R^L)_{x_0}
(\mR, e_l)  \frac{Z^\alpha}{\alpha !}
= \frac{1}{2} \left \langle (\nabla ^{X} \nabla ^{X}\mJ)_{(\mR,\mR)} \mR,
 e_l \right \rangle 
+ \frac{1}{6}\left \langle R^{TX} (\mR, \mJ \mR)\mR, e_l \right
\rangle .
\end{align}
From (\ref{g2}), (\ref{g7}), (\ref{g18}) and (\ref{g22}), we know that
\begin{equation}\label{g23}
\begin{split}
I_1 =\,&\frac{1}{2}\, b_i \sum_{|\alpha|=2}(\partial ^{\alpha}R^L)_{x_0}
(\mR,\tfrac{\partial}{\partial \ov{z}_i})
\frac{Z^\alpha}{\alpha !} \\
&+\frac{1}{12}\left[{\tfrac{\partial}{\partial z_i}}
\left \langle R^{TX}_{x_0} (\mR, \mJ \mR)\mR, 
\tfrac{\partial}{\partial \ov{z}_i} \right\rangle
- {\tfrac{\partial}{\partial \ov{z}_i}}
\left \langle R^{TX}_{x_0}  (\mR, \mJ \mR)\mR, 
\tfrac{\partial}{\partial z_i} \right\rangle \right] \\
&+ \frac{1}{4} \left[ {\tfrac{\partial}{\partial z_i}} 
\left \langle (\nabla ^{X} \nabla ^{X}\mJ)_{(\mR,\mR)} \mR,
 \tfrac{\partial}{\partial \ov{z}_i} \right \rangle
  -  {\tfrac{\partial}{\partial \ov{z}_i}} 
  \left \langle (\nabla ^{X} \nabla ^{X}\mJ)_{(\mR,\mR)} \mR,
 \tfrac{\partial}{\partial z_i} \right \rangle  \right] .
\end{split}
\end{equation}
The definition of $\nabla ^{X}\nabla ^{X}\mJ$, $R^{TX}$ and \eqref{g9} imply,
for $U,V,W,Y\in TX$, 
\begin{equation}\label{g24}
\begin{split}
&\left \langle R^{TX}(U,V)W,Y \right \rangle
=\left \langle R^{TX}(W,Y) U,V\right \rangle,\\
&R^{TX}(U,V)W+ R^{TX}(V,W)U+ R^{TX}(W,U)V=0, \\
& (\nabla ^{X}\nabla ^{X}\mJ)_{(U,V)}- (\nabla ^{X}\nabla ^{X}\mJ)_{(V,U)}
=[R^{TX}(U,V),\mJ],\\
&\left \langle (\nabla ^{X}\nabla ^{X}\mJ)_{(Y,U)}V,W \right \rangle
+ \left \langle (\nabla ^{X}\nabla ^{X}\mJ)_{(Y,V)}W,U \right \rangle
+\left \langle (\nabla ^{X}\nabla ^{X}\mJ)_{(Y,W)}U,V \right \rangle =0. 
\end{split}
\end{equation}
Note that $\mJ,(\nabla ^{X}\nabla ^{X}\mJ)_{(Y,U)}$ are skew-adjoint, 
by \eqref{0ue52} and \eqref{g24}, 
\begin{equation}\label{g25}
\begin{split}
{\tfrac{\partial}{\partial z_i}}& 
\left \langle (\nabla ^{X} \nabla ^{X}\mJ)_{(\mR,\mR)} \mR,
 \tfrac{\partial}{\partial \ov{z}_i} \right \rangle 
 - {\tfrac{\partial}{\partial \ov{z}_i}} 
 \left \langle (\nabla ^{X} \nabla ^{X}\mJ)_{(\mR,\mR)} \mR,
 \tfrac{\partial}{\partial z_i} \right \rangle \\
=&  \left \langle 2 (\nabla ^{X} \nabla ^{X}\mJ)_{(\mR,\mR)}
\tfrac{\partial}{\partial z_i}
+ 2(\nabla ^{X} \nabla ^{X}\mJ)_{(\mR,\tfrac{\partial}{\partial z_i})} \mR
+\Big[R^{TX} (\tfrac{\partial}{\partial z_i},  \mR), \mJ \Big] \mR,
 \tfrac{\partial}{\partial \ov{z}_i} \right \rangle \\
 -& \left \langle 
2(\nabla ^{X} \nabla ^{X}\mJ)_{(\mR,\tfrac{\partial}{\partial \ov{z}_i})} \mR
+ \Big[R^{TX} (\tfrac{\partial}{\partial \ov{z}_i},  \mR), 
\mJ \Big] \mR,
 \tfrac{\partial}{\partial z_i} \right \rangle  \\
=&4 \left \langle (\nabla ^{X} \nabla ^{X}\mJ)_{(\mR,\mR)}
\tfrac{\partial}{\partial z_i},
 \tfrac{\partial}{\partial \ov{z}_i} \right \rangle
+ \left \langle  2 a_i  R^{TX} (\mR, \tfrac{\partial}{\partial z_i}) \mR
- R^{TX} (\mR,\mJ  \mR) \tfrac{\partial}{\partial z_i},
 \tfrac{\partial}{\partial \ov{z}_i} \right \rangle,
\end{split}
\end{equation}
\begin{equation*}
\begin{split}
{\tfrac{\partial}{\partial z_i}}&\left \langle R^{TX} (\mR, \mJ \mR)\mR, 
\tfrac{\partial}{\partial \ov{z}_i} \right\rangle
- {\tfrac{\partial}{\partial \ov{z}_i}}
\left \langle R^{TX} (\mR, \mJ \mR)\mR, 
\tfrac{\partial}{\partial z_i} \right\rangle \\
&=  2\left \langle R^{TX} (\mR, \mJ \mR)\tfrac{\partial}{\partial z_i}, 
\tfrac{\partial}{\partial \ov{z}_i} \right\rangle
+2a_i \left \langle R^{TX} (\mR, \tfrac{\partial}{\partial \ov{z}_i})\mR, 
\tfrac{\partial}{\partial z_i} \right\rangle \\
&+\left \langle R^{TX} (\tfrac{\partial}{\partial z_i}, \mJ \mR)\mR, 
\tfrac{\partial}{\partial \ov{z}_i} \right\rangle
-\left \langle R^{TX} (\tfrac{\partial}{\partial \ov{z}_i}, \mJ \mR)\mR, 
\tfrac{\partial}{\partial z_i} \right\rangle \\
&=3\left \langle R^{TX} (\mR, \mJ \mR)\tfrac{\partial}{\partial z_i}, 
\tfrac{\partial}{\partial \ov{z}_i} \right\rangle
+2a_i \left \langle R^{TX} (\mR, \tfrac{\partial}{\partial z_i})\mR, 
\tfrac{\partial}{\partial \ov{z}_i} \right\rangle .
\end{split}
\end{equation*}
Thus by (\ref{g23})-(\ref{g25}),
\begin{multline}\label{g26}
I_1 = \frac{1}{2} b_i \sum_{|\alpha|=2}(\partial ^{\alpha}R^L)_{x_0}
(\mR,\tfrac{\partial}{\partial \ov{z}_i})\frac{Z^\alpha}{\alpha !} 
+\left \langle (\nabla ^{X} \nabla ^{X}\mJ)_{(\mR,\mR)}
\tfrac{\partial}{\partial z_i},
 \tfrac{\partial}{\partial \ov{z}_i} \right \rangle \\
+ \frac{2}{3}a_i \left \langle R^{TX} 
(\mR, \tfrac{\partial}{\partial z_i})\mR, 
\tfrac{\partial}{\partial \ov{z}_i} \right\rangle.
\end{multline}
Now by (\ref{g2}), (\ref{g18}) and (\ref{g24}),
\begin{multline}\label{g27}
I_2 =\,\frac{4}{3}\, b_j \left[  \left \langle R^{TX} 
(\tfrac{\partial}{\partial z_i},\tfrac{\partial}{\partial\ov{z}_i})
\mR, \tfrac{\partial}{\partial\ov{z}_j}\right \rangle 
- \left \langle R^{TX} (\mR,\tfrac{\partial}{\partial z_i} )
\tfrac{\partial}{\partial \ov{z}_i},
\tfrac{\partial}{\partial \ov{z}_j}\right \rangle\right] \\
+ \frac{1}{3}\,b_ib_j 
\left \langle R^{TX} (\mR,\tfrac{\partial}{\partial\ov{z}_i})
\mR, \tfrac{\partial}{\partial\ov{z}_j}\right \rangle
- \frac{8}{3}   \left \langle R^{TX} 
(\tfrac{\partial}{\partial z_j},\tfrac{\partial}{\partial z_i})
\tfrac{\partial}{\partial \ov{z}_i}, 
\tfrac{\partial}{\partial\ov{z}_j}\right \rangle \\
+ \frac{4}{3}\,\left[  \left \langle R^{TX} 
(\tfrac{\partial}{\partial z_i},\tfrac{\partial}{\partial\ov{z}_i})
\tfrac{\partial}{\partial z_j}, 
\tfrac{\partial}{\partial\ov{z}_j}\right \rangle 
- \left \langle R^{TX} (\tfrac{\partial}{\partial z_j},
\tfrac{\partial}{\partial\ov{z}_i}) \tfrac{\partial}{\partial z_i}, 
\tfrac{\partial}{\partial\ov{z}_j}\right \rangle \right]\\
= \frac{4}{3}\, b_j \left[  \left \langle R^{TX} 
(\tfrac{\partial}{\partial z_i},\tfrac{\partial}{\partial\ov{z}_i})
\mR, \tfrac{\partial}{\partial\ov{z}_j}\right \rangle 
- \left \langle R^{TX} (\mR,\tfrac{\partial}{\partial z_i})
\tfrac{\partial}{\partial \ov{z}_i}, 
\tfrac{\partial}{\partial\ov{z}_j}\right \rangle \right] \\
 +\frac{1}{3}\,b_ib_j 
\left \langle R^{TX} (\mR,\tfrac{\partial}{\partial\ov{z}_i})
\mR, \tfrac{\partial}{\partial\ov{z}_j}\right \rangle
+ 4 \left \langle R^{TX} (\tfrac{\partial}{\partial z_i},
\tfrac{\partial}{\partial z_j}) \tfrac{\partial}{\partial \ov{z}_i}, 
\tfrac{\partial}{\partial \ov{z}_j}\right \rangle .
\end{multline}
Finally (\ref{g3}), (\ref{g14}), (\ref{g19}), (\ref{g26}) and (\ref{g27}), 
yield (\ref{g17}).
\end{proof}

\noindent
From (\ref{g10}), (\ref{g13}), (\ref{g15}), (\ref{1c52}) and (\ref{g17}), 
follows
\begin{multline}\label{0g27}
F_{1,2}= J_{1,2}P^N =  \Big(
2R^E (\tfrac{\partial}{\partial z_i},\tfrac{\partial}{\partial \ov{z}_i})
+4 \left \langle R^{TX} (\tfrac{\partial}{\partial z_i},
\tfrac{\partial}{\partial z_j}) \tfrac{\partial}{\partial \ov{z}_i}, 
\tfrac{\partial}{\partial \ov{z}_j}\right \rangle +  \Phi \Big)  P^N\\
+ P^N \Big(\left \langle (\nabla ^{X} \nabla ^{X}\mJ)_{(\mR,\mR)}
\tfrac{\partial}{\partial z_i},
 \tfrac{\partial}{\partial \ov{z}_i} \right \rangle 
+\frac{\sqrt{-1}}{4} 
\tr_{|TX} \Big(\nabla ^{X} \nabla ^{X}(J\mJ)\Big)_{(\mR,\mR)} \\
+ \frac{1}{9} |(\nabla_\mR^X \mJ) \mR |^2
+ \frac{4}{9} \left\langle(\nabla ^{X}_\mR\mJ)\mR,
\tfrac{\partial} {\partial z_i}\right\rangle b^+_i \cL_0^{-1}b_i 
\left\langle(\nabla ^{X}_\mR\mJ)\mR,\tfrac{\partial}
{\partial\ov{z}_i}\right\rangle \Big)  P^N.
\end{multline}


\subsection{The coefficients $b_{q,0}$\,}\label{s4.2}
In the rest of this Section we assume that  ${\bf J}=J$. 
A very useful observation is that 
\eqref{g9}, \eqref{g11} imply
\begin{equation}\label{g29}
\begin{split}
&\text{$\mJ=-2\pi \sqrt{-1} J$ and $a_i=2\pi$ in (\ref{0ue52}), $\tau=2\pi n$.
$\nabla ^{X}_U J$ is skew-adjoint} \\
&\text{and the tensor $\left\langle(\nabla ^{X}_\cdot J)\cdot,
\cdot \right\rangle$ 
is of the type $(T^{*(1,0)}X)^{\otimes 3}\oplus (T^{*(0,1)}X)^{\otimes 3}$.}
\end{split}
\end{equation}
Before computing $b_{q,0}$, we establish the relation between
the scalar curvature $r^X$ and $|\nabla ^{X}J| ^2$.
\begin{lemma}
\begin{align}\label{g28}
r^X= 8 \left \langle  R^{TX} (\tfrac{\partial}{\partial z_i},
\tfrac{\partial}{\partial \ov{z}_j})\tfrac{\partial}{\partial z_j},
\tfrac{\partial}{\partial \ov{z}_i} \right \rangle
-   \frac{1}{4}  |\nabla ^{X}J|^2.
\end{align}
\end{lemma}
\begin{proof}
By (\ref{g29}),
\begin{align}\label{g30}
|\nabla ^{X}J| ^2
=4 \Big \langle (\nabla ^{X}_{\tfrac{\partial}{\partial z_i}}J) e_j,
 (\nabla ^{X}_{\tfrac{\partial}{\partial \ov{z}_i}}J) e_j\Big \rangle
=8 \Big \langle (\nabla ^{X}_{\tfrac{\partial}{\partial z_i}} J)
\tfrac{\partial}{\partial z_j} ,
(\nabla ^{X}_{\tfrac{\partial}{\partial \ov{z}_i}}J)
\tfrac{\partial}{\partial \ov{z}_j}\Big \rangle . 
\end{align}
By (\ref{g0}),  (\ref{g9}) and (\ref{g29}),
\begin{multline}\label{g31}
 \Big \langle(\nabla ^{X}_{\tfrac{\partial}{\partial z_j}}J)
\tfrac{\partial}{\partial z_i},
(\nabla ^{X}_{\tfrac{\partial}{\partial \ov{z}_i}}J)
\tfrac{\partial}{\partial \ov{z}_j}\Big\rangle 
= 2\Big \langle(\nabla ^{X}_{\tfrac{\partial}{\partial z_j}}J)
\tfrac{\partial}{\partial z_i},\tfrac{\partial}{\partial z_k}\Big\rangle 
\Big \langle (\nabla ^{X}_{\tfrac{\partial}{\partial \ov{z}_i}}J)
\tfrac{\partial}{\partial \ov{z}_j}, \tfrac{\partial}{\partial \ov{z}_k}  
  \Big\rangle \\
= 2\Big \langle(\nabla ^{X}_{\tfrac{\partial}{\partial z_i}}J)
\tfrac{\partial}{\partial z_k}
- (\nabla ^{X}_{\tfrac{\partial}{\partial z_k}}J)
\tfrac{\partial}{\partial z_i},
 \tfrac{\partial}{\partial z_j}\Big\rangle
\Big \langle (\nabla ^{X}_{\tfrac{\partial}{\partial \ov{z}_i}}J)
\tfrac{\partial}{\partial \ov{z}_k}, \tfrac{\partial}{\partial \ov{z}_j}   
 \Big\rangle \\
=\Big\langle (\nabla ^{X}_{\tfrac{\partial}{\partial z_i}} J)
\tfrac{\partial}{\partial z_k} ,
(\nabla ^{X}_{\tfrac{\partial}{\partial \ov{z}_i}}J)
\tfrac{\partial}{\partial \ov{z}_k}\Big\rangle 
-\Big \langle(\nabla ^{X}_{\tfrac{\partial}{\partial z_k}}J)
\tfrac{\partial}{\partial z_i},
(\nabla ^{X}_{\tfrac{\partial}{\partial \ov{z}_i}}J)
\tfrac{\partial}{\partial \ov{z}_k}\Big\rangle .
\end{multline}
By \eqref{g30} and \eqref{g31},
\begin{align}\label{g32}
\Big\langle(\nabla ^{X}_{\tfrac{\partial}{\partial z_j}}J)
\tfrac{\partial}{\partial z_i},
(\nabla ^{X}_{\tfrac{\partial}{\partial \ov{z}_i}}J)
\tfrac{\partial}{\partial \ov{z}_j}\Big\rangle
=  \frac{1}{16}|\nabla ^{X}J| ^2.
\end{align}
Now, from (\ref{g11}), we get 
\begin{align}\label{g34} 
(\nabla ^{X}\nabla ^{X}J)_{(U,V)}J +(\nabla ^{X}_UJ)\circ (\nabla ^{X}_VJ)
+ (\nabla ^{X}_VJ)\circ (\nabla ^{X}_UJ)
+J (\nabla ^{X}\nabla ^{X}J)_{(U,V)}=0.
\end{align}
 thus from (\ref{g24}), (\ref{g29}) and (\ref{g34}),
for $u_1,u_2, u_{3}\in T^{(1,0)}X$, $\ov{v}_1,\ov{v}_2\in T^{(0,1)}X$,
\begin{equation}\label{g35} 
\begin{split}
&(\nabla ^{X}\nabla ^{X}J)_{(u_1,u_2)}u_{3},\,
(\nabla ^{X}\nabla ^{X}J)_{(\ov{v}_1,\ov{v}_2)}u_{3}\in T^{(0,1)}X,
\quad (\nabla ^{X}\nabla ^{X}J)_{(u_1,\ov{v}_2)}u_{3}\in T^{(1,0)}X, \\
&2\sqrt{-1} \left\langle (\nabla ^{X}\nabla ^{X}J)_{(u_1,\ov{v}_1)}
u_2 ,  \ov{v}_2\right\rangle
= \left\langle  (\nabla ^{X}_{u_1}J) u_2,  
(\nabla ^{X}_{\ov{v}_1}J)\ov{v}_2\right\rangle \, .  
\end{split}
\end{equation}
Formulas \eqref{g24} and \eqref{g35} yield
\begin{multline}\label{0g35}
\left\langle (\nabla ^{X}\nabla ^{X}J)_{(u_1,u_2)}\ov{v}_1, 
\ov{v}_2\right\rangle
= -\left\langle (\nabla ^{X}\nabla ^{X}J)_{(u_1,\ov{v}_1)}\ov{v}_2,
 u_2\right\rangle 
- \left\langle (\nabla ^{X}\nabla ^{X}J)_{(u_1,\ov{v}_2)}u_2,
 \ov{v}_1\right\rangle\\
= \frac{1}{2\sqrt{-1}} \left\langle (\nabla ^{X}_{u_1}J)u_2, 
(\nabla ^X_{\ov{v}_1} J)\ov{v}_2
- (\nabla ^X_{\ov{v}_2} J)\ov{v}_1\right\rangle.
\end{multline}
From (\ref{g24}), (\ref{g32}) and (\ref{0g35}), we deduce
\begin{multline}\label{g36}
\left \langle R^{TX} (\tfrac{\partial}{\partial z_i},
\tfrac{\partial}{\partial z_j}) \tfrac{\partial}{\partial \ov{z}_i}, 
\tfrac{\partial}{\partial \ov{z}_j}\right \rangle
=\frac{\sqrt{-1}}{2} \left \langle [R^{TX} (\tfrac{\partial}{\partial z_i},
\tfrac{\partial}{\partial z_j}), J] \tfrac{\partial}{\partial \ov{z}_i}, 
\tfrac{\partial}{\partial \ov{z}_j}\right \rangle\\
=\frac{\sqrt{-1}}{2} \left \langle 
\Big( (\nabla ^{X}\nabla ^{X}J)_{(\tfrac{\partial}{\partial z_i},
\tfrac{\partial}{\partial z_j})} 
- (\nabla ^{X}\nabla ^{X}J)_{(\tfrac{\partial}{\partial z_j},
\tfrac{\partial}{\partial z_i})} \Big)
\tfrac{\partial}{\partial \ov{z}_i}, 
\tfrac{\partial}{\partial \ov{z}_j}\right \rangle\\
=  \frac{1}{4} \Big\langle(\nabla ^{X}_{\tfrac{\partial}{\partial z_i}}J)
\tfrac{\partial}{\partial z_j} , 
(\nabla ^{X}_{\tfrac{\partial}{\partial \ov{z}_i}}J)
\tfrac{\partial}{\partial \ov{z}_j}\Big\rangle 
=   \frac{1}{32}|\nabla ^{X}J|^2.
\end{multline}
The scalar curvature $r^X$ of $(X,g^{TX})$ is
\begin{equation}\label{g38}
\begin{split}
r^X =& -\left \langle R^{TX} (e_i,e_j)e_i,e_j \right \rangle
= -4 \left \langle R^{TX} (\tfrac{\partial}{\partial z_i},e_j)
\tfrac{\partial}{\partial \ov{z}_i},e_j \right \rangle\\
=&-8 \left \langle  R^{TX} (\tfrac{\partial}{\partial z_i},
\tfrac{\partial}{\partial z_j})\tfrac{\partial}{\partial \ov{z}_i},
\tfrac{\partial}{\partial \ov{z}_j} \right \rangle 
- 8 \left \langle  R^{TX} (\tfrac{\partial}{\partial z_i},
\tfrac{\partial}{\partial \ov{z}_j})\tfrac{\partial}{\partial \ov{z}_i},
\tfrac{\partial}{\partial z_j} \right \rangle  .
\end{split}
\end{equation}
In conclusion, relations \eqref{g36} and \eqref{g38} imply \eqref{g28}.
\end{proof}

\noindent
From \eqref{g13} and \eqref{g29} we know
\begin{align}\label{g41}
\mO_1=\frac{2}{3}b_i \left\langle(\nabla ^{X}_{\overline z}\mJ)\overline z,
\tfrac{\partial}{\partial\ov{z}_i}\right\rangle
- \frac{2}{3} \left\langle(\nabla ^{X}_z\mJ)z,
\tfrac{\partial}{\partial z_i} \right\rangle b_i^+.
\end{align}
Hence by (\ref{g2}), (\ref{g14}), (\ref{g29}) and (\ref{g41}), 
\begin{equation}\label{g42}
\begin{split}
(\mO_1P^N )(Z,Z') =& \frac{2}{3} \Big( b_i  
\left\langle(\nabla ^{X}_{\overline z}\mJ)\overline z,\tfrac{\partial}
{\partial\ov{z}_i}\right\rangle P^N\Big) (Z,Z')\\
=& \frac{2}{3} \Big\{\Big( \frac{b_ib_j}{2\pi} 
\Big\langle (\nabla ^{X}_{\tfrac{\partial}
{\partial\ov{z}_j}}\mJ)\ov{z}',\tfrac{\partial}
{\partial\ov{z}_i} \Big\rangle 
+ b_i \left\langle(\nabla ^{X}_{\ov{z}'}\mJ)\ov{z}',
\tfrac{\partial}{\partial\ov{z}_i}\right\rangle  \Big) P^N\Big\} (Z,Z').
\end{split}
\end{equation}
By Theorem \ref{t3.4}, \eqref{g15} and \eqref{g42},  we have
\begin{align}\label{0g42}
&(\cL_0^{-1} P^{N^\bot} \mO_1 P^N)(Z,Z')\\
&\hspace*{10mm}= \frac{2}{3}\Big\{ \Big( \frac{b_ib_j}{16\pi ^2} \Big\langle 
(\nabla ^{X}_{\tfrac{\partial}{\partial\ov{z}_j}}\mJ)\ov{z}',
\tfrac{\partial}{\partial\ov{z}_i}\Big\rangle 
+ \frac{b_i}{4\pi} \left\langle(\nabla ^{X}_{\ov{z}'}\mJ)\ov{z}',
\tfrac{\partial}{\partial\ov{z}_i}\right\rangle \Big) P^N\Big\}(Z,Z'),
\nonumber\\
&P^N \mO_1 \cL_0^{-1} P^{N^\bot} \mO_1 P^N
=- \frac{2}{3} P^N \left\langle(\nabla ^{X}_z\mJ)z,
\tfrac{\partial}{\partial z_k} \right\rangle b_k^+ \cL_0^{-1} 
P^{N^\bot} \mO_1P^N .\nonumber
\end{align}
 (\ref{g0}), (\ref{g2}), (\ref{g14}), \eqref{g41} and (\ref{0g42}) imply
\begin{multline}\label{g43}
 \frac{2}{3} \Big(\left\langle(\nabla ^{X}_z\mJ)z,
\tfrac{\partial}{\partial z_k} \right\rangle b_k^+ \cL_0^{-1} 
P^{N^\bot} \mO_1P^N \Big )(Z,Z')
= \frac{4}{9} \left\{\left\langle(\nabla ^{X}_z\mJ)z,
\tfrac{\partial}{\partial z_k} \right\rangle \right. \\ 
 \left. \times \left(- \frac{b_i}{4\pi} \Big\langle 
(\nabla ^{X}_{\tfrac{\partial}{\partial\ov{z}_i}} \mJ)
\tfrac{\partial}{\partial\ov{z}_k} 
+(\nabla ^{X}_{\tfrac{\partial}{\partial\ov{z}_k}} \mJ)
\tfrac{\partial}{\partial\ov{z}_i} ,\ov{z}' \Big\rangle 
+  \left\langle(\nabla ^{X}_{\ov{z}'}\mJ)\ov{z}',
\tfrac{\partial}{\partial\ov{z}_k}\right\rangle \right) P^N\right\} (Z,Z')\\
= \left\{\left[ -\frac{b_i}{9\pi}  \left\langle(\nabla ^{X}_z\mJ)z,
\tfrac{\partial}{\partial z_k} \right\rangle 
\left\langle (\nabla ^{X}_{\tfrac{\partial}{\partial\ov{z}_i}} \mJ)
\tfrac{\partial}{\partial\ov{z}_k} 
+(\nabla ^{X}_{\tfrac{\partial}{\partial\ov{z}_k}} \mJ)
\tfrac{\partial}{\partial\ov{z}_i} ,\ov{z}' \right\rangle  \right. \right.\\
 -\frac{2}{9\pi}   \left\langle(\nabla ^{X}_{z}\mJ)
\tfrac{\partial}{\partial z_i} 
+(\nabla ^{X}_{\tfrac{\partial}{\partial z_i}}\mJ) z,
\tfrac{\partial}{\partial z_k} \right\rangle   
\left\langle (\nabla ^{X}_{\tfrac{\partial}{\partial\ov{z}_i}} \mJ)
\tfrac{\partial}{\partial\ov{z}_k} 
+(\nabla ^{X}_{\tfrac{\partial}{\partial\ov{z}_k}} \mJ)
\tfrac{\partial}{\partial\ov{z}_i} ,\ov{z}' \right\rangle \\
\left. \left.  + \frac{2}{9}\left\langle(\nabla ^{X}_{z}\mJ)z,
(\nabla ^{X}_{\ov{z}'}\mJ)\ov{z}'\right\rangle\right] P^N \right\} (Z,Z').
\end{multline}
Thanks to \eqref{g11}, (\ref{g14}), (\ref{g29}), (\ref{g30}) and (\ref{g31}) we obtain
\begin{multline}\label{g44} 
\frac{1}{9}  \Big |(\nabla_\mR^X \mJ) \mR \Big |^2P^N(Z,Z')= 
\frac{8\pi ^2}{9}   \left\langle(\nabla ^{X}_{z}J)z, 
(\nabla ^{X}_{\ov{z}}J)\ov{z} \right\rangle    P^N(Z,Z')  \\
=\frac{8\pi ^2}{9} \left\{\left\langle(\nabla ^{X}_{z}J)z,
\frac{b_ib_j}{4\pi ^2} (\nabla ^{X}_{\tfrac{\partial}{\partial \ov{z}_i}}J)
\tfrac{\partial}{\partial \ov{z}_j} 
+ \frac{b_i}{2 \pi} \Big[ 
(\nabla ^{X}_{\tfrac{\partial}{\partial \ov{z}_i}}J)
\ov{z}'+(\nabla ^{X}_{\ov{z}'}J)
\tfrac{\partial}{\partial \ov{z}_i} \Big] 
+ (\nabla ^{X}_{\ov{z}'}J)\ov{z}' \right\rangle P^N\right\}(Z,Z')\\
= \frac{8\pi ^2}{9} \left\{ \left[\left\langle \frac{b_ib_j}{4\pi ^2}
 (\nabla ^{X}_{\tfrac{\partial}{\partial \ov{z}_i}}J)
\tfrac{\partial}{\partial \ov{z}_j} 
+\frac{b_i}{2 \pi}\Big (
(\nabla ^{X}_{\tfrac{\partial}{\partial \ov{z}_i}}J)
\ov{z}'+(\nabla ^{X}_{\ov{z}'}J)
\tfrac{\partial}{\partial \ov{z}_i}\Big) ,
(\nabla ^{X}_{z}J)z \right\rangle\right.\right.\\
+\frac{b_i}{2\pi ^2}  \left\langle
(\nabla ^{X}_{\tfrac{\partial}{\partial z_j}}J)z
+ (\nabla ^{X}_{z}J)\tfrac{\partial}{\partial z_j}, 
(\nabla ^{X}_{\tfrac{\partial}{\partial \ov{z}_i}}J)
\tfrac{\partial}{\partial \ov{z}_j} 
+(\nabla ^{X}_{\tfrac{\partial}{\partial \ov{z}_j}}J)
\tfrac{\partial}{\partial \ov{z}_i}  \right\rangle  \\
+ \frac{1}{\pi} \left\langle 
(\nabla ^{X}_{\tfrac{\partial}{\partial z_i}}J)z
+(\nabla ^{X}_{z}J)\tfrac{\partial}{\partial z_i},
 (\nabla ^{X}_{\tfrac{\partial}{\partial \ov{z}_i}}J)\ov{z}'
+(\nabla ^{X}_{\ov{z}' }J)\tfrac{\partial}{\partial \ov{z}_i}   
\right\rangle\\
\left.\left.  + \left\langle (\nabla ^{X}_{z}J)z,
 (\nabla ^{X}_{\ov{z}'}J)\ov{z}'\right\rangle  
+ \frac{3}{16 \pi ^2}  |\nabla ^{X}J|^2 
 \right] P^N\right\}(Z,Z'). 
\end{multline}
Taking into account (\ref{g24}), (\ref{g35}) and $\left\langle [R^{TX}(\ov{z},z), J]
\tfrac{\partial}{\partial z_i} ,  
\tfrac{\partial}{\partial \ov{z}_i}\right\rangle=0$, we get
\begin{multline}\label{g45} 
\left\langle (\nabla ^{X}\nabla ^{X}J)_{(\mR,\mR)}
\tfrac{\partial}{\partial z_i}, 
 \tfrac{\partial}{\partial \ov{z}_i}\right\rangle
= \left\langle (\nabla ^{X}\nabla ^{X}J)_{(z,\ov{z})}
\tfrac{\partial}{\partial z_i} +(\nabla ^{X}\nabla ^{X}J)_{(\ov{z},z)}
\tfrac{\partial}{\partial z_i} ,  
\tfrac{\partial}{\partial \ov{z}_i}\right\rangle\\
=-\sqrt{-1} \left\langle (\nabla ^{X}_zJ)
\tfrac{\partial}{\partial z_i} ,  
(\nabla ^{X}_{\ov{z}}J)\tfrac{\partial}{\partial \ov{z}_i}\right\rangle.
\end{multline}
From (\ref{g14}),  (\ref{g30}) and (\ref{g45}),
\begin{multline}\label{g46} 
 \left\langle (\nabla ^{X}\nabla ^{X}\mJ)_{(\mR,\mR)}
\tfrac{\partial}{\partial z_i},  
\tfrac{\partial}{\partial \ov{z}_i}\right\rangle P^N (Z,Z')
= - 2 \pi
\left\langle  (\nabla ^{X}_zJ)
\tfrac{\partial}{\partial z_i} ,  
(\nabla ^{X}_{\ov{z}}J)\tfrac{\partial}{\partial \ov{z}_i}\right\rangle 
P^N(Z,Z')\\
=- 2 \pi \Big\{ \Big\langle  (\nabla ^{X}_zJ)
\tfrac{\partial}{\partial z_i} ,  \frac{b_j}{2 \pi} 
(\nabla ^{X}_{\tfrac{\partial}{\partial \ov{z}_j}}J)
\tfrac{\partial}{\partial \ov{z}_i}+ (\nabla ^{X}_{\ov{z}'}J)
\tfrac{\partial}{\partial \ov{z}_i} \Big\rangle P^N\Big\}(Z,Z')\\
=-  \Big\{\Big [\Big\langle b_j    
(\nabla ^{X}_{\tfrac{\partial}{\partial \ov{z}_j}}J)
\tfrac{\partial}{\partial \ov{z}_i} 
+ 2\pi (\nabla ^{X}_{\ov{z}'}J)\tfrac{\partial}{\partial \ov{z}_i}
, (\nabla ^{X}_zJ)\tfrac{\partial}{\partial z_i} \Big\rangle 
 +\frac{1}{4}  |\nabla ^{X}J|^2 \Big ] P^N\Big\}(Z,Z').
\end{multline}
Recall that the polynomial $J_{q,2q} (Z,Z')$ was defined in (\ref{c86}).
From $J\mJ=2\pi \sqrt{-1}$, (\ref{g15}), (\ref{0g27}), 
(\ref{g36}) and (\ref{g41})-(\ref{g46}), 
$J_{1,2}(Z,Z')$ is a polynomial on $z,\ov{z}'$, and each monomial of  $J_{1,2}$
has the same degree in $z$ and $\ov{z}'$; moreover 
\begin{align}\label{g48} 
J_{1,2}(0,0)
= \frac{1}{24}   |\nabla ^{X}J| ^2_{x_0} +2 R^E_{x_0} 
(\tfrac{\partial}{\partial z_i},\tfrac{\partial}{\partial \ov{z}_i})
+ \Phi_{x_0}.
\end{align}
Using (\ref{g7}), (\ref{1c52}) and the recurrence, 
we infer that each monomial of  $J_{q,2q}$
has the same degree in $z$ and $\ov{z}'$, and 
\begin{align}\label{g49} 
J_{q,2q}(0,0) = (J_{1,2}(0,0))^q.
\end{align}
In view of (\ref{g7}), (\ref{c86}), (\ref{1c54}), (\ref{g48}) and 
 (\ref{g49}) we obtain (\ref{0.7}).
\comment{\begin{rem}\label{t4.2} If ${\bf J}\neq J$, it seems that the 
monomial of  $J_{1,2}$ will involve the term $z^\alpha$ and ${\ov{z}'}^\beta$,
thus (\ref{g49}) will not hold in general.
\end{rem}}


\subsection{The coefficient $b_{0,1}$ }\label{s4.3}
By (\ref{1c54}), we need to compute $F_{0,2} (0,0)$.
By \eqref{g42} and \eqref{0g42}, we know that 
\begin{equation}\label{0g52}
(\mO_1 P^N)(Z,0)=0, \quad (\cL_0^{-1}P^{N^\bot} \mO_1 P^N)(0,Z^\prime)=0. 
\end{equation}
Thus the first and last two terms in \eqref{g51} are zero at $(0,0)$.
Thus we only need to compute $-(\cL_0^{-1}P^{N^\bot} \mO_2 P^N)(0,0)$,
 since the third and fourth terms  in \eqref{g51} 
are adjoint of the first two terms by Remark \ref{0t3.3}. 

Let $h_i(z)$ and $f_{ij}(z)$, $(i,j=1,\cdots, n)$ be arbitrary 
polynomials in $z$.
By Theorem \ref{t3.4}, (\ref{g2}), (\ref{g7}) and (\ref{g14}), we have
\begin{equation}\label{g52}
\begin{split}
&(b_ih_iP^N)(0,0) = -2 \frac{\partial h_i}{\partial z_i}(0),\, \quad
(b_ib_jf_{ij}P^N)(0,0) = 
4 \frac{\partial ^2 f_{ij}}{\partial z_i\partial z_j}(0),\\
&(\cL_0^{-1} b_if_{ij}b_jP^N)(0,0)= 
-\frac{1}{2\pi}\frac{\partial ^2 f_{ij}}{\partial z_i\partial z_j}(0).
\end{split}
\end{equation}
Owing to Theorem \ref{t3.4}, (\ref{g30}), (\ref{g32}), (\ref{g44})
and (\ref{g52}), 
\begin{multline}\label{g53}
-\frac{1}{9} \Big( \cL_0^{-1}  P^{N^\bot} 
\Big |(\nabla_\mR^X \mJ) \mR \Big |^2P^N\Big)(0,0)=- \frac{8}{9} \left\{
\Big[
\frac{b_ib_j}{32\pi } \left\langle(\nabla ^{X}_{z}J)z, 
(\nabla ^{X}_{\tfrac{\partial}{\partial \ov{z}_i}}J)
\tfrac{\partial}{\partial \ov{z}_j}  \right\rangle \right.\\
\left. +\frac{b_i}{8\pi} \Big\langle
(\nabla ^{X}_{\tfrac{\partial}{\partial z_j}}J)z
+ (\nabla ^{X}_{z}J)\tfrac{\partial}{\partial z_j}, 
(\nabla ^{X}_{\tfrac{\partial}{\partial \ov{z}_i}}J)
\tfrac{\partial}{\partial \ov{z}_j} 
+(\nabla ^{X}_{\tfrac{\partial}{\partial \ov{z}_j}}J)
\tfrac{\partial}{\partial \ov{z}_i}  \Big\rangle \Big] P^N\right\}(0,0)\\
=  \frac{1}{9 \pi}
\Big \langle(\nabla ^{X}_{\tfrac{\partial}{\partial z_i}}J)
\tfrac{\partial}{\partial z_j} +
(\nabla ^{X}_{\tfrac{\partial}{\partial z_j}}J)
\tfrac{\partial}{\partial z_i},
(\nabla ^{X}_{\tfrac{\partial}{\partial \ov{z}_i}}J)
\tfrac{\partial}{\partial \ov{z}_j}
 +2 (\nabla ^{X}_{\tfrac{\partial}{\partial \ov{z}_i}}J)
\tfrac{\partial}{\partial \ov{z}_j}\Big \rangle 
 = \frac{1}{16\pi } |\nabla ^X J  |^2\,,
\end{multline}
and by Theorem \ref{t3.4}, (\ref{g46}) and  (\ref{g52}), 
\begin{multline}\label{g54} 
 -\Big( \cL_0^{-1} P^{N^\bot} \left\langle (\nabla ^{X}\nabla ^{X}\mJ)_{(\mR,\mR)}
\tfrac{\partial}{\partial z_i},  
\tfrac{\partial}{\partial \ov{z}_i}\right\rangle P^N\Big) (0,0)\\
= \Big(\frac{b_j}{4\pi} \Big\langle  (\nabla ^{X}_z J)
\tfrac{\partial}{\partial z_i} ,   
(\nabla ^{X}_{\tfrac{\partial}{\partial \ov{z}_j}}J)
\tfrac{\partial}{\partial \ov{z}_i} \Big\rangle P^N\Big) (0,0)
= -\frac{1}{16\pi}|\nabla ^X J |^2.
\end{multline}
Observe that by (\ref{g14}), for a polynomial $g(z)$ in $z$, 
the constant term of $\frac{1}{P^N}\frac{b^\alpha}{2^{|\alpha|}} g(z)P^N$ is 
the constant term of $(\tfrac{\partial}{\partial z})^\alpha g$.  
Thus in the term $-\cL_0^{-1}  P^{N^\bot}\mO_2 P^N$, 
by \eqref{g14}, the contribution of $\frac{1}{2}b_i 
\sum_{|\alpha|=2}(\partial ^{\alpha}R^L)_{x_0}
(\mR,\tfrac{\partial}{\partial \ov{z}_i})\frac{Z^\alpha}{\alpha !}$ 
in $\mO_2$ consists of the terms whose total degree of $b_i$ and $\ov{z}_j$ 
is same as the degree of $z$. Hence we only need to consider 
the contribution from the terms where the degree of $z$ is $2$.
By (\ref{g22}), (\ref{g29}), (\ref{g35}), (\ref{0g35}) 
and $\left\langle [ R^{TX}(\ov{z}, z), \mJ]z,
\tfrac{\partial}{\partial \ov{z}_i}\right\rangle =0$, this term is 
\begin{multline}\label{g55} 
I_3= \frac{1}{4} b_i \Big [
 \left\langle (\nabla ^{X}\nabla ^{X}\mJ)_{(z,z)}\ov{z}
+(\nabla ^{X}\nabla ^{X}\mJ)_{(z,\ov{z})}z
+(\nabla ^{X}\nabla ^{X}\mJ)_{(\ov{z},z)}z,     
\tfrac{\partial}{\partial \ov{z}_i}\right\rangle\\
+\frac{1}{3} \left\langle R^{TX}(\ov{z},\mJ z)z +  R^{TX}(z,\mJ \ov{z})z, 
\tfrac{\partial}{\partial \ov{z}_i} \right\rangle \Big ]\\
=-\frac{\pi}{4} b_i \Big [ \Big\langle (\nabla ^{X}_zJ)z, 
3 (\nabla ^{X}_{\ov{z}}J)\tfrac{\partial}{\partial \ov{z}_i} 
- (\nabla ^{X}_{\tfrac{\partial}{\partial \ov{z}_i}}J)\ov{z}\Big\rangle
+\frac{4}{3}\left\langle R^{TX}(z,\ov{z})z, 
\tfrac{\partial}{\partial \ov{z}_i} \right\rangle  \Big ].
\end{multline}
Therefore, from (\ref{g14}), (\ref{g24}), (\ref{g32}), (\ref{g36}), (\ref{g52})
and (\ref{g55}), we get
\begin{multline}\label{g58} 
- \Big(\cL_0^{-1}  P^{N^\bot}I_3P^N\Big) (0,0)
=  \frac{\pi}{4} \Big\{\cL_0^{-1} b_i \Big [\frac{4}{3}
\left\langle R^{TX}(z, \tfrac{\partial}{\partial \ov{z}_j})z, 
\tfrac{\partial}{\partial \ov{z}_i} \right\rangle\\
 + \Big\langle  (\nabla ^{X}_zJ)z, 
3 (\nabla ^{X}_{\tfrac{\partial}{\partial \ov{z}_j} }J)
\tfrac{\partial}{\partial \ov{z}_i} 
- (\nabla ^{X}_{\tfrac{\partial}{\partial \ov{z}_i}}J)
\tfrac{\partial}{\partial \ov{z}_j}\Big\rangle  \Big ]
\frac{b_j}{2\pi}P^N\Big\} (0,0) \\
= - \frac{1}{16\pi} \Big[ \frac{4}{3} \left\langle 
R^{TX} (\tfrac{\partial}{\partial z_j},\tfrac{\partial}{\partial \ov{z}_j})
\tfrac{\partial}{\partial z_i}
+ R^{TX}(\tfrac{\partial}{\partial z_i},
\tfrac{\partial}{\partial \ov{z}_j})\tfrac{\partial}{\partial z_j},
\tfrac{\partial}{\partial \ov{z}_i}\right\rangle\\
+\Big\langle
(\nabla ^{X}_{\tfrac{\partial}{\partial z_i}} J)\tfrac{\partial}{\partial z_j}
+  (\nabla ^{X}_{\tfrac{\partial}{\partial z_j}} J)
\tfrac{\partial}{\partial z_i},
3(\nabla ^{X}_{\tfrac{\partial}{\partial \ov{z}_i}}J)
\tfrac{\partial}{\partial \ov{z}_j}
-(\nabla ^{X}_{\tfrac{\partial}{\partial \ov{z}_j}}J)
\tfrac{\partial}{\partial \ov{z}_i}\Big\rangle  \Big ]  \\
= - \frac{5}{192\pi } |\nabla ^X J |^2 - \frac{1}{6\pi}\left\langle
R^{TX} (\tfrac{\partial}{\partial z_i},\tfrac{\partial}{\partial \ov{z}_j})
\tfrac{\partial}{\partial z_j}, 
\tfrac{\partial}{\partial \ov{z}_i}\right\rangle.
\end{multline} 
Thanks to (\ref{g14}), (\ref{g36}) and (\ref{g52}) we have   
\begin{multline}\label{g59}
\frac{1}{3} \Big(P^{N^\bot}  \left \langle R^{TX} 
(\mR,\tfrac{\partial}{\partial z_i})
\mR, \tfrac{\partial}{\partial\ov{z}_i} \right \rangle P^N\Big) (0,0)\\
= \frac{1}{3} \Big(P^{N^\bot} 
\left \langle R^{TX} (z,\tfrac{\partial}{\partial z_i})
\tfrac{\partial}{\partial \ov{z}_j}
+ R^{TX}(\tfrac{\partial}{\partial \ov{z}_j}, 
\tfrac{\partial}{\partial z_i})z,
\tfrac{\partial}{\partial\ov{z}_i}\right \rangle 
\frac{b_j}{2\pi} P^N\Big) (0,0)\\
= -\frac{1}{3\pi} \Big [  \left \langle R^{TX} 
(\tfrac{\partial}{\partial z_j},\tfrac{\partial}{\partial z_i})
\tfrac{\partial}{\partial \ov{z}_j}
+ R^{TX}(\tfrac{\partial}{\partial \ov{z}_j}, 
\tfrac{\partial}{\partial z_i})\tfrac{\partial}{\partial z_j},
\tfrac{\partial}{\partial\ov{z}_i}\right \rangle \Big ] \\
= - \frac{1}{96\pi } |\nabla ^X J |^2 + \frac{1}{3\pi}\left\langle
R^{TX} (\tfrac{\partial}{\partial z_i},\tfrac{\partial}{\partial \ov{z}_j})
\tfrac{\partial}{\partial z_j}, 
\tfrac{\partial}{\partial \ov{z}_i}\right\rangle.
\end{multline} 
By (\ref{g29}), (\ref{g52}),  (\ref{g53}), (\ref{g54}), 
(\ref{g58}), (\ref{g59}), and the discussion above \eqref{g55}, 
\begin{multline}\label{g60}
-\Big(\cL_0^{-1}  P^{N^\bot} \mO_2  P^N\Big) (0,0)=
- \Big\{ \Big[\frac{b_ib_j }{24\pi}
\left \langle R^{TX} (z,\tfrac{\partial}{\partial\ov{z}_i})
z, \tfrac{\partial}{\partial\ov{z}_j}\right \rangle
+\frac{b_i}{4\pi}R^E (z,\tfrac{\partial}{\partial \ov{z}_i})\\
+  \frac{b_j}{3\pi} \left(  \left \langle R^{TX} 
(\tfrac{\partial}{\partial z_i},\tfrac{\partial}{\partial\ov{z}_i})
z, \tfrac{\partial}{\partial\ov{z}_j}\right \rangle 
- \left \langle R^{TX} (z,\tfrac{\partial}{\partial z_i})
\tfrac{\partial}{\partial \ov{z}_i}, 
\tfrac{\partial}{\partial\ov{z}_j}\right \rangle \right) 
\Big]P^N \Big\} (0,0) \\
+\frac{1}{3} \Big( P^{N^\bot}  \left \langle R^{TX}
(\mR,\tfrac{\partial}{\partial z_i})
\mR, \tfrac{\partial}{\partial\ov{z}_i} \right \rangle P^N \Big) (0,0)
 - \Big(\cL_0^{-1}  P^{N^\bot} I_3  P^N \Big) (0,0) \\
= -\frac{1}{6 \pi} \left \langle R^{TX} 
(\tfrac{\partial}{\partial z_i},\tfrac{\partial}{\partial\ov{z}_i})
 \tfrac{\partial}{\partial z_j}
+ R^{TX} (\tfrac{\partial}{\partial z_j},\tfrac{\partial}{\partial\ov{z}_i})
 \tfrac{\partial}{\partial z_i}, 
\tfrac{\partial}{\partial\ov{z}_j}\right \rangle 
+ \frac{1}{2\pi}R^E (\tfrac{\partial}{\partial z_i},
\tfrac{\partial}{\partial \ov{z}_i}) \\
+\frac{2}{3\pi} \left[  \left \langle R^{TX} 
(\tfrac{\partial}{\partial z_i},\tfrac{\partial}{\partial\ov{z}_i})
\tfrac{\partial}{\partial z_j}, 
\tfrac{\partial}{\partial\ov{z}_j}\right \rangle 
- \left \langle R^{TX} (\tfrac{\partial}{\partial z_j},
\tfrac{\partial}{\partial z_i}) \tfrac{\partial}{\partial \ov{z}_i}, 
\tfrac{\partial}{\partial\ov{z}_j}\right \rangle \right]\\
- \frac{7}{192\pi } |\nabla ^X J |^2 + \frac{1}{6\pi}\left\langle
R^{TX} (\tfrac{\partial}{\partial z_i},\tfrac{\partial}{\partial \ov{z}_j})
\tfrac{\partial}{\partial z_j}, 
\tfrac{\partial}{\partial \ov{z}_i}\right\rangle\\
= \frac{1}{2\pi} \left\langle
R^{TX} (\tfrac{\partial}{\partial z_i},\tfrac{\partial}{\partial \ov{z}_j})
\tfrac{\partial}{\partial z_j}, 
\tfrac{\partial}{\partial \ov{z}_i}\right\rangle
+\frac{1}{2\pi}R^E (\tfrac{\partial}{\partial z_i},
\tfrac{\partial}{\partial \ov{z}_i}).
\end{multline}
Formulas (\ref{g28}), (\ref{g60}) and the discussion 
at the beginning of Section  \ref{s4.3} yield finally 
\begin{equation}\label{g61}
\begin{split}
b_{0,1}(x_0)=&F_{0,2}(0,0) 
=- \Big(\cL_0^{-1}  P^{N^\bot} \mO_2  P^N \Big) (0,0)
- \Big ( \Big(\cL_0^{-1}  P^{N^\bot} \mO_2  P^N \Big) (0,0)\Big)^*\\
=& \,\frac{1}{8\pi}\Big[r^X_{x_0}+ \frac{1}{4} |\nabla ^X J |^2_{x_0}
+ 2\sqrt{-1} R^E_{x_0} (e_j,Je_j)\Big].
\end{split}
\end{equation}
The proof of Theorem \ref{t0.2} is complete.

\begin{rem} In the K\"ahler case, i.e. $J$ is integrable and $L,E$
are holomorphic, then $\mO_1=0$, and the above computation 
simplifies a lot.
\end{rem}

\section{Applications} \label{s5}

In this Section, we discuss various applications of our results.
In Section \ref{s5.1}, 
we study the density of states
function of $\Delta_{p,\Phi}$. 
In Section \ref{s5.3},  we  explain how to handle
the first-order pseudo-differential operator $D_b$
of Boutet de Monvel and Guillemin \cite{BoG} which was studied 
extensively by Shiffman and  Zelditch \cite{SZ02}.
In Section \ref{s5.31}, we prove a symplectic version of 
 the convergence of the Fubini-Study metric of 
an ample line bundle \cite{Tian}.
In Section \ref{s5.4}, we show how to handle the operator 
$\ov{\partial}+\ov{\partial}^*$ when $X$ is K\" ahler but ${\bf J}\neq J$.
Finally, in Sections \ref{s5.5}, \ref{s5.6}, we establish some generalizations 
for non-compact or singular manifolds.


\subsection{Density of states function}\label{s5.1}

Let $(X,\om)$ be a compact symplectic manifold of real dimension $2n$ 
and $(L,\nabla ^L, h^L)$ 
is a pre-quantum line bundle as in Section \ref{s1}.
Assume that $E$ is the trivial bundle $\bC$, 
$\Phi=0$ and $\mathbf{J}=J$. The latter 
means, by \eqref{0.1}, that $g^{TX}$ is the Riemannian metric 
associated to $\omega$ and $J$. We denote by 
$\vol (X) =\int_X \frac{\omega ^n}{n!}$ the Riemannian volume of $(X,g^{TX})$.
Recall that $d_p$ is defined in (\ref{0.0}).

 Our aim is to describe the asymptotic distribution of the energies
of the bound states as $p$ tends to infinity.
We define the spectrum counting function of $\Delta_{p}:=\Delta_{p,0}$ by
$N_{p}(\lambda)=\#\left\{i\,:\,\lambda_{i,p}\leqslant\lambda\right\}$
and the spectral density measure on $[-C_L,C_L]$ by
\begin{gather}
\nu_p=\frac{1}{d_p}\,\frac{d}{d\lambda}N_{p}(\lambda)\,,\quad
\lambda\in[-C_L,C_L]\,.\label{f2}
\end{gather}
Clearly, $\nu_p$ is a sum of Dirac measures supported on
$\spec\Delta_{p}\cap[-C_L,C_L]$. Set
\begin{equation}\label{f3}
\varrho:X\longrightarrow\bR\,,\quad\varrho(x)
=\frac{1}{24}\abs{\nabla^{X}J}^2\,.
\end{equation}
\begin{thm}\label{spectral-density}
The weak limit of the sequence $\{\nu_p\}_{p\geqslant1}$ is the direct image measure
$\varrho_{\ast}\Big(\dfrac{1}{\vol (X)}\dfrac{\om^n}{n!}\Big)$, that is,
for any continuous function $f\in\cC([-C_L,C_L])$, we have
\begin{equation}\label{f4}
\lim_{p\to\infty}\int_{-C_L}^{C_L}f\,d\nu_p
=\frac{1}{\vol (X)}\int_{X}(f\circ\varrho)\,\frac{\om^n}{n!}\:\:.
\end{equation}
\end{thm}
\begin{proof}
By  \eqref{0.4}, we have for $q\geqslant1$ (now $E$ is trivial):
$B_{q,p}(x)=\sum_{i=1}^{d_p} \lambda_{i,p}^q\,\abs{S^{p}_{i}(x)}^2$,
which yields by integration over $X$,
\begin{gather}
\frac{1}{d_p}\int_{X} B_{q,\,p}\,dv_{X}=
\frac{1}{d_p}\sum_{i=1}^{d_p} \lambda_{i,p}^q=
\int_{-C_L}^{C_L}\lambda^q\,d\nu_p(\lambda)\,,\label{f5}
\end{gather}
since $S^{p}_{i}$ have unit $L^2$ norm. On the other hand, 
\eqref{0.0}, \eqref{0.6} entail for $p\to \infty$,
\begin{align}\label{f6}
\frac{1}{d_p}\int_{X} B_{q,\,p}\,dv_{X}&=
\frac{p^n}{d_p}\int_{X} b_{q,0}\,dv_{X}+\frac{\cO(p^{n-1})}{d_p}\\
&= \frac{1}{\vol (X)}\int_{X}\varrho^q\,dv_{X} + \cO(p^{-1}).\nonumber
\end{align}
We infer from \eqref{f5} and \eqref{f6} that \eqref{f4} holds for
$f(\lambda)=\lambda^q$, $q\geqslant1$. Since this is obviously true for
$f(\lambda)\equiv1$, too, we deduce it holds for all polynomials.
Upon invoking the Weierstrass approximation theorem, we get \eqref{f4}
for all continuous functions on $[-C_L,C_L]$. This achieves the proof.
\end{proof}
\begin{rem} \label{t5.2} A function $\varrho$ satisfying \eqref{f4} is called spectral
density function. Its existence and uniqueness were demonstrated
by Guillemin-Uribe \cite{GU}.
As for the explicit formula of $\varrho$, the paper \cite{BU2} 
is dedicated to its computation. Our formula (\ref{f3}) is different from
\cite[Theorem 1.2]{BU2}
\footnote{In \cite[(3.7)]{BU2}, the leading term of 
$G_{0j}$ should be $\kappa ^{-1/2} b^{(1)}_j$ which was missed therein, 
as the principal terms of $\frac{\partial}{\partial s}$, 
$\frac{\partial}{\partial y^j}$ are $\partial_0$, $T^l_j \partial_l$ by
\cite[equation after (3.11)]{BU2}. Now, from \cite[(3.5)]{BU2},
$ b^{(1)}_j$ is $\frac{1}{2}\langle Jz,T^l_j \partial_l\rangle$.
 Thus $\mL_0$ in \cite[(3.8)]{BU2} is incorrect.
\cite[Theorem 1.2]{BU2} is $\varrho(x)=
-\frac{5}{24}\abs{\nabla^{X}J}^2$. }.
\end{rem}

An interesting corollary of (\ref{f3}) and (\ref{f4}) 
is the following result which was first stated in \cite[Cor. 1.3]{BU2}.
\begin{cor}\label{t5.3}
The spectral density function is identically zero if and only if 
$(X,\omega,J)$ is K{\"a}hler.
\end{cor}

\begin{rem}\label{t5.4}
Theorem \ref{spectral-density} can be slightly generalized. Assume namely
that ${\bf J}=J$ and 
$E$ is a Hermitian vector bundle as in Section \ref{s1} such that
$R^E=\eta\otimes\Id_E$, $\Phi=\varphi\Id_E$,
where $\eta$ is a $2$-form and $\varphi$ a real function on $X$.
Then there exists a spectrum density function satisfying \eqref{f4}
given by 
\begin{equation}\label{f8}
\varrho:X\longrightarrow\bR\,,\quad\varrho(x)=
\frac{1}{24}\abs{\nabla^{X}J}^2+\frac{\sqrt{-1}}{2}\eta(e_j,Je_j) +
\varphi\,.
\end{equation}
The proof is similar to the previous one, as
 $\tr_{E_x}B_{q,p}(x)=\sum_{i=1}^{d_p} \lambda_{i,p}^q\,\abs{S^{p}_{i}(x)}^2$.
\end{rem}

\subsection{Almost-holomorphic Szeg{\" o} kernels}\label{s5.3}

We use the notations and assumptions from Section \ref{s5.1}, 
especially, ${\bf J}=J$. Then $\tau=2\pi n$.

Let $Y= \{u\in L^*, |u|_{h^{L^*}}=1\}$ be the unit circle bundle in $L^*$. 
Then the smooth sections of $L^p$ can be identified to the smooth functions
$\cC^\infty (Y)_p= \{ f\in \cC^\infty(Y,\bC); 
f(ye ^{i\theta})$ $= e ^{ip\theta} f(y)\, \, {\rm for}\,  \,
e ^{i\theta}\in S^1, y\in Y \}$,
here $ye ^{i\theta}$ is the $S^1$ action on $Y$.

The connection $\nabla ^L$ on  $L$ induces a connection 
on the $S^1$-principal bundle $\pi : Y\to X$, 
and let $T^H Y \subset TY$ be the corresponding horizontal bundle.
 Let $g^{TY}= \pi ^ * g^{TX} \oplus d\theta ^2$ be the metric on 
$TY= T^HY\oplus TS ^1$, with $d\theta ^2$ 
the standard metric on $S^1= \bR/2\pi \bZ$.
Let $\Delta_Y$ be the Bochner-Laplacian on $(Y,g^{TY})$, then by construction, 
it commutes with the generator $\partial_\theta$ of the circle action, 
and so it commutes with the horizontal Laplacian
\begin{align}\label{f20}
 \Delta_h= \Delta_Y+ \partial_\theta ^2,
\end{align}
 then $\Delta_h$ on $\cC^\infty (Y)_p$  is identical with $\Delta ^{L^p}$
on $\cC^\infty (X,  L^p)$ (cf. \cite[\S 2.1]{BU1}).

In \cite[Lemma 14.11, Theorem A 5.9]{BoG}, \cite{BoS},
 \cite[(3.13)]{GU}, they construct 
a self-adjoint second-order pseudodifferential operator $Q$ on $Y$
such that 
\begin{align}\label{f21}
V= \Delta_h +\sqrt{-1} \tau \partial_\theta - Q
\end{align}
is a self-adjoint  pseudodifferential operator of order zero on $Y$,
and $V,Q$ commute with the $S^1$-action.
The orthogonal projection $\Pi$ onto the kernel of $Q$ is called 
the {\em Szeg\"o  projector} associated with the almost CR manifold $Y$. 
In fact,  the Szeg\"o 
projector is not unique or canonically defined, but the above construction
defines a canonical choice of $\Pi$  modulo smoothing operators. 
In the complex case, the construction produces the usual Szeg\"o projector
$\Pi$. 

We denote the operators on $\cC^\infty (X, L^p)$ corresponding to $Q$, $V$,
 $\Pi$ by $Q_p,V_p$, $\Pi_p$, especially, 
$V_p(x,y)= \frac{1}{2\pi} \int_0^{2\pi} e ^{-ip\theta} 
V(x e ^{i\theta},y)d\theta$. Then by (\ref{f21}),
\begin{align}\label{f22}
Q_p= \Delta ^{L^p}- p \tau -V_p.
\end{align}
By \cite[\S 4]{GU}, there exists $\mu_1>0$ such that for $p$ large,
\begin{align}\label{f23}
\spec Q_p \subset \{0\}\cup [\mu_1 p,+\infty[.
\end{align}
Since the operator $V_p$ is uniformly bounded in $p$, 
naturally, from (\ref{0.3}), (\ref{0.0}), we get 
\begin{align}\label{f24}
\dim \Ker Q_p=d_p =\int_X \td(TX)\ch(L^p).
\end{align}

Now we explain how to study the  Szeg\"o 
projector $\Pi_p$ \footnote{As Professor Sj\"ostrand pointed out to us, 
in general, $\Pi_p-P_{0,p}$ is not $\cO(p^{-\infty})$
 as $p\to \infty$, where $P_{0,p}$ is the smooth kernel
of the operator $\Delta_{0,p}$ (Definition \ref{d3.0}). 
This can also be seen from the presence of a contribution coming from $\Phi$ in the 
expression \eqref{0.6} of the coefficient $b_{0,2}$.}. This can  
be done from our point of view.
Recall $\wi{F}$ is the function defined after (\ref{0c3}).
Let $\Pi_p(x,x')$, $\wi{F}(Q_p)(x,x')$ be the smooth kernels 
of $\Pi_p$, $\wi{F}(Q_p)$ with respect to the volume form $dv_X(x')$. 

Note that $V_p$ is a 0-order pseudodifferential operator on $X$ induced from
a 0-order pseudodifferential operator on $Y$. 
Thus from (\ref{f22}), (\ref{f23}), we have the analogue of 
\cite[Proposition 3.1]{DLM} (cf. Proposition \ref{0t3.0}):
for any $l,m\in \bN$, there exists $C_{l,m}>0$ such that for $p\geqslant1$,
\begin{align}\label{f25}
|\wi{F}(Q_p)(x,x') -\Pi_p (x,x') |_{\cC^m(X\times X)}\leqslant C_{l,m} p^{-l}.
\end{align}
By finite propagation speed  \cite[\S 4.4]{T1}, we know that 
$\wi{F}(Q_p)(x, x')$ only depends on the restriction of $Q_p$ to
$B^X(x,\var)$, and is zero if $d(x, x') \geqslant\var$. It transpires
that the asymptotic of $\Pi_p(x,x')$ as $p\to \infty$ is localized on
a neighborhood of $x$. 
Thus we can  translate our analysis from $X$ to the manifold $\bR^{2n} \simeq
T_{x_0}X=:X_{0}$ as in Section \ref{s3.2}, especially, we extend 
$\nabla ^L$ to a
Hermitian connection $\nabla ^{L_{0}}$  on
$(L_0,h^{L_0})=(X_0\times L_{x_{0}},h^{L_{x_{0}}} )$ 
on $T_{x_0}X$ in such a way so that
we still have positive curvature $R ^{L_{0}}$; in addition 
$R^{L_{0}}=R ^{L}_{x_{0}}$ outside a compact set. 

Now, by using a micro-local partition of unity, one can still construct
the operator $Q^{X_0}$ as in 
\cite[Lemma 14.11, Theorem A 5.9]{BoG}, \cite{BoS}, \cite[(3.13)]{GU},
 such that 
$V^{X_0}$ differs from $V$ by a smooth operator in a neighborhood of $0$.
On $X_0$, and $Q^{X_0}$ still verifies (\ref{f23}).
Thus we can work on $\cC^\infty (X_0,\bC)$ as in Section \ref{s3.3}.
We rescale then the coordinates as in (\ref{c27}) 
and use the norm (\ref{u0}). The $V^{X_0}_p$ is a  
0-order pseudodifferential operator on $X_0$ induced from
a 0-order pseudodifferential operator on $Y_0$. This guarantees that the operator 
rescaled from $V^{X_0}_p$ will have the similar expansion as \eqref{c30}
with leading term $t^2 R_2$ in the sense of pseudo-differential operators. 

From (\ref{f24}) and \cite[(3.89)]{DLM}, 
similar to the argument in \cite[Theorem 3.18]{DLM}, we can also get the full 
off diagonal expansion for $\Pi_p$, 
which is an extension of \cite[Theorem 1]{SZ02}, where the authors obtain
(\ref{f26}) for $|Z|, |Z'|\leqslant  C/\sqrt{p}$ with $C>0$ fixed. 
More precisely, recalling that $P^N$ is the Bergman kernel of $\cL_0$ as in \eqref{ue62}, \eqref{g7}
we have:
\begin{thm} \label{tue17} 
There exist polynomials ${\bf j}_r(Z,Z')$ $(r\geqslant0)$ of $Z,Z'$ 
with the same parity with $r$,
and  ${\bf j}_0=1$, $C''>0$ such that 
for any $k,m,m'\in \bN$, there exist $N\in \bN, C>0$ such that for
$\alpha,\alpha'\in \bZ^{2n}$, $|\alpha|+|\alpha'|\leqslant m$,
$Z,Z'\in T_{x_0}X$, $|Z|, |Z'|\leqslant  \var$, $x_0\in X$, $p>1$,
\begin{multline}\label{f26}
\left |\frac{\partial^{|\alpha|+|\alpha'|}}
{\partial Z^{\alpha} {\partial Z'}^{\alpha'}}
\left (\frac{1}{p^n}  \Pi_p(Z,Z')
-\sum_{r=0}^k ({\bf j}_r P^N)(\sqrt{p} Z,\sqrt{p} Z') \kappa ^{-\frac{1}{2}}(Z)
\kappa ^{-\frac{1}{2}}(Z')
p^{-r/2}\right )\right |_{\cC ^{m'}(X)}\\
\leqslant C  p^{-(k+1-m)/2}  (1+|\sqrt{p} Z|+|\sqrt{p} Z'|)^N
\exp (- \sqrt{C''\mu_1 } \sqrt{p} |Z-Z'|)+ \cO(p^{-\infty}).
\end{multline}
\end{thm}
The term $\kappa ^{-\frac{1}{2}}$ in (\ref{f26}) 
comes from the conjugation of the operators as in (\ref{1c53}), 
$\cC ^{m'}(X)$ is the $\cC ^{m'}$ -norm for the parameter $x_0\in X$,
and we use the trivializations from Section \ref{s3.2},
the term $\cO(p^{-\infty})$ means that for any $l,l_1\in \bN$, 
there exists $C_{l,l_1}>0$ such that its $\cC^{l_1}$-norm is dominated 
by $C_{l,l_1} p^{-l}$.
We leave the details to the interested reader.

\subsection{Symplectic version of Kodaira Embedding Theorem}\label{s5.31}

Let $(X,\omega)$ be a compact symplectic manifold of real dimension  
$2n$ and let $(L,\nabla^L,h^L)$ be a pre-quantum line bundle 
and let $g^{TX}$ be a Riemannian metric on $X$
as in Section \ref{s1}. 

Recall that $\mH_p\subset\cC^\infty(X,L^p)$ is the span of those  
eigensections of $\Delta_p=\Delta ^{L^p}-\tau p$ corresponding 
to eigenvalues from $[-C_L,C_L]$.  
We denote by $\proj\mH^*_p$ the projective space associated to the dual of  
$\mH_p$ and we identify $\proj\mH^*_p$ with the Grassmannian of hyperplanes in  
$\mH_p$. The {\em base locus} of $\mH_p$ is the 
set $\operatorname{Bl}(\mH_p)=\{x\in X:s(x)=0\,\, 
\text{for all}\,s\in\mH_p\,\}$.  
As in algebraic geometry, we define the Kodaira map  
\begin{equation}\label{sz0} 
\begin{split}
&\Phi_p:X\smallsetminus\operatorname{Bl}(\mH_p)\longrightarrow\proj\mH^*_p\\ 
&\Phi_p(x)=\{s\in\mH_p:s(x)=0\} 
\end{split}
\end{equation} 
which sends $x\in X\smallsetminus\operatorname{Bl}(\mH_p)$ 
to the hyperplane of sections  vanishing at $x$.  
 Note that $\mH_p$ is endowed with the induced $L^2$ product \eqref{0c2}  
so there is a well--defined Fubini--Study metric 
$g_{FS}$ on $\proj\mH^*_p$ with the associated form $\omega_{FS}$.  
\begin{thm} \label{sym-Kodaira}
Let $(L,\nabla^L)$ be a pre--quantum line bundle over a compact symplectic  
manifold $(X,\omega)$. The following assertions hold true{\rm:} 
 
{\rm (i)} For large $p$, the Kodaira maps $\Phi_p:X\longrightarrow\proj 
\mH^*_p$ are well defined.  
 
{\rm (ii)} The induced Fubini--Study metric $\frac{1}{p}\Phi^*_p 
(\omega_{FS})$ converges in the $\cC^\infty$ topology 
to $\omega$\,{\rm;} for any  
$l\geqslant 0$ there exists $C_l>0$ such that  
\begin{equation} \label{sz1} 
\Big|\frac{1}{p}\,\,\Phi^*_p(\omega_{FS})
-\omega\Big|_{\cC^l}\leqslant\frac{C_l}{p} .
\end{equation} 

{\rm (iii)} For large $p$ the Kodaira maps $\Phi_p$ are embeddings.  
\end{thm} 
 
\begin{rem} 
1) Assume that $X$ is K\" ahler and $L$ is a holomorphic bundle.  
Then $\Delta_p$ is the twice the Kodaira-Laplacian 
and $\mH_p$ coincides with the  
space $H^0(X,L^p)$ of holomorphic sections of $L^p$. Then (i) and (iii) are  
simply the Kodaira embedding theorem. Assertion (ii) is due to Tian  
\cite[Theorem A]{Tian} as an answer to a conjecture of Yau. 
In \cite{Tian} the case  
$l=2$ is considered and the left--hand side of \eqref{sz1} is estimated by  
$C_l/\sqrt p$. Ruan \cite{Ru} proved the $\cC^\infty$ convergence and  
improved the bound to $C_l/p$. Both papers use the peak section method,  
based on $L^2$--estimates for ${\overline\partial}$. 
A proof for $l=0$ using the heat kernel  
appeared in Bouche \cite{Bou}. Finally, Zelditch deduced (ii) from the  
asymptotic expansion of the Szeg\" o  kernel \cite{Zelditch}.

2) Borthwick and Uribe \cite[Theorem 1.1]{BU1},  
Shiffman and Zelditch \cite[Theorems\,2, \,3]{SZ02} prove a different  
symplectic version of \cite[Theorem A]{Tian} when ${\bf J}=J$. 
Instead of $\mH_p$, they use the space  
$H^0_J(X,L^p):={\rm Im}(\Pi_p)$ (cf. \cite[p.601]{BU1}, \cite[\S 2.3]{SZ02},
\S \ref{s5.3}) 
of `almost holomorphic sections' proposed by Boutet de Monvel  
and Guillemin \cite{BoG}, \cite{BoS}. 
\end{rem} 
 
\begin{proof} 
Let us first give an alternate description of the map $\Phi_p$ which relates  
it to the Bergman kernel. Let $\{S^p_i\}^{d_p}_{i=1}$ be any orthonormal  
basis of $\mH_p$ with respect to the inner product \eqref{0c2}. Once we  
have fixed a basis, we obtain an identification $\mH_p\cong\mH^*_p\cong 
\bC^{d_p}$ and $\proj\mH^*_p\cong\bC\proj^{d_p-1}$. Consider the commutative  
diagram.  
\begin{equation} \label{sz1.1} 
\begin{CD} 
X\smallsetminus\operatorname{Bl}(\mH_p)@>\Phi_p>>\proj\mH^*_p\\ 
@VV{\Id}V     @VV{\cong}V\\ 
X\smallsetminus\operatorname{Bl}(\mH_p)@>\widetilde\Phi_p>>\bC\proj^{d_p-1} 
\end{CD} 
\end{equation} 
Then  
\begin{equation} \label{sz2} 
\Phi^*_p(\omega_{FS})=\widetilde\Phi^*_p\Big(\frac{\sqrt{-1}}{2\pi} 
\partial\overline\partial\log\sum^{d_p}_{j=1}|w_j|^2\Big), 
\end{equation} 
where $[w_1,\ldots,w_{d_p}]$ are homogeneous coordinates 
in $\bC\proj^{d_p-1}$.  
To describe $\widetilde\Phi_p$ in a neighborhood of a point  
$x_0\in X\smallsetminus\operatorname{Bl}(\mH_p)$, we choose
 a local frame  
$e_L$ of $L$ and write $S^p_i=f^p_ie^{\otimes p}_L$ 
for some smooth functions $f^p_i$. Then  
\begin{equation} \label{sz3} 
\widetilde\Phi_p(x)=[f^p_1(x);\ldots;f^p_{d_p}(x)],
\end{equation} 
and this does not depend on the choice of the frame $e_L$.  
 
(i) Let us choose an unit frame $e_L$ of $L$. Then $|S^p_i|^2= 
|f^p_i|^2|e_L|^{2p}=|f^p_i|^2$, hence 
\begin{equation*}
B_{0,p}=\sum^{d_p}_{i=1}|S^p_i|^2=\sum^{d_p}_{i=1}|f^p_i|^2 .
\end{equation*} 
Since $b_{0,0}>0$, the asymptotic expansion \eqref{0.6} shows that $B_{0,p}$ 
does not vanish on $X$ for $p$ large enough, so the sections $\{S^p_i\}^ 
{d_p}_{i=1}$ have no common zeroes. Therefore $\Phi_p$ and $\widetilde 
\Phi_p$ are defined on all $X$.  
 
(ii) Let us fix $x_0\in X$. We identify a small geodesic ball 
$B^X(x_0,\varepsilon)$  
to $B^{T_{x_0}X}(0,\varepsilon)$ by means of the exponential map and consider  
a trivialization of $L$ as in Section \ref{s3.2}, i.e. we trivialize $L$ 
by using  an unit frame $e_L(Z)$ which is parallel with respect 
to $\nabla ^L$ along $[0,1]\ni u\to uZ$
for $Z\in B^{T_{x_0}X}(0,\varepsilon)$. 
We can express the Fubini--Study metric as   
\begin{equation*} 
\frac{\sqrt{-1}}{2\pi}\partial\overline\partial
\log\Big(\sum^{d_p}_{j=1}|w_j|^2\Big)=\frac{\sqrt{-1}} 
{2\pi}\left[\frac{1}{|w|^2}\sum^{d_p}_{j=1}dw_j\wedge d\overline w_j 
-\frac{1}{|w|^4}\sum^{d_p}_{j,k=1}\overline w_jw_k\,dw_j\wedge d\overline 
w_k\right] ,
\end{equation*} 
and therefore, from \eqref{sz3}, 
\begin{multline} \label{sz4}
\Phi^*_p(\omega_{FS})(x_0)=\frac{\sqrt{-1}}{2\pi}
\left[\frac{1}{|f^p|^2}\sum^{d_p}_ 
{j=1}df^p_j\wedge d\overline{f^p_j} 
-\frac{1}{|f^p|^4}\sum^{d_p}_{j,k=1}\overline{f^p_j}f^p_k\,df^p_j\wedge d 
\overline{f^p_k}\right](x_0)\\ 
=\frac{\sqrt{-1}}{2\pi}\big[f^p(x_0,x_0)^{-1}d_xd_yf^p(x,y) 
-f^p(x_0,x_0)^{-2}d_xf^p(x,y)\wedge d_yf^p(x,y)\big]\vert_{x=y=x_0},
\end{multline}
where $f^p(x,y)=\sum^{d_p}_{i=1}f^p_i(x)\overline{f^p_i}(y)$ and  
$|f^p(x)|^2=f^p(x,x)$.  
Since 
\begin{equation} \label{0sz4}
P_{0,p}(x,y)=f^p(x,y)e^p_L(x)\otimes e^{p}_L(y)^*,
 \end{equation} 
thus $P_{0,p}(x,y)$ is $f^p(x,y)$ under our trivialization of $L$. 
By \eqref{0c30}, Theorem \eqref{0t3.6}, and \eqref{1c53}, we  
obtain  
\begin{multline} \label{0sz5}
\frac{1}{p}\,\,\Phi^*_p(\omega_{FS})(x_0)
=\frac{\sqrt{-1}}{2\pi}\Big[\,\frac{1}{F_{0,0}}d_xd_yF_{0,0} 
- \frac{1}{F_{0,0}^2}d_xF_{0,0}\wedge d_yF_{0,0}\Big](0,0)\\
-\frac{\sqrt{-1}}{2\pi}\frac{1}{\sqrt{p}}\Big[\,\frac{1}{F_{0,0}^2}
(d_xF_{0,1}\wedge d_yF_{0,0}+ d_xF_{0,0}\wedge d_yF_{0,1})\Big](0,0)
+\cO\big(1/p\big). 
\end{multline} 
Using again \eqref{g7},  \eqref{1c52}, we obtain 
\begin{equation} \label{0sz6}
\frac{1}{p}\,\,\Phi^*_p(\omega_{FS})(x_0)
=\frac{\sqrt{-1}}{4\pi}\sum^n_{j=1} a_jdz_j\wedge 
d \overline z_j|_{x_0}+\cO\big(1/p\big)=\omega(x_0)+\cO\big(1/p\big),
\end{equation} 
and the convergence takes place in the $\cC^\infty$ topology with respect to  
$x_0\in X$.  
 
(iii) Since $X$ is compact, we have to prove two things for $p$ sufficiently  
large: (a)~$\Phi_p$ are immersions and (b) $\Phi_p$ are injective.  
We note that (a) follows immediately from \eqref{sz1}. 

To prove (b) let us assume the contrary, namely that there exists 
a sequence of distinct points $x_p\neq y_p$ such that 
$\Phi_p(x_p)=\Phi_p(y_p)$. Relation \eqref{sz1.1} implies that 
$\widetilde\Phi_p(x_p)=\widetilde\Phi_p(y_p)$,
where $\widetilde\Phi_p$ is defined by any particular choice of basis.

The key observation is that Theorem \ref{t3.8} ensures the existence of
 a sequence  of {\em peak sections\/} at each point of $X$. 
The construction goes like follows.  
Let $x_0\in X$ be fixed. Since $\Phi_p$ is point base free for large $p$,  
we can consider the hyperplane $\Phi_p(x_0)$ of all sections of $\mH_p$  
vanishing at $x_0$. We construct then an orthonormal basis $\{S^p_i\}^{d_p}_ 
{i=1}$ of $\mH_p$ such that the first $d_p-1$ elements belong to $\Phi_ 
p(x_0)$. Then $S^p_{d_p}$ is a unit norm generator of the orthogonal  
complement of $\Phi_p(x_0)$, and will be denoted by $S^p_{x_0}$.
 This is a peak section at $x_0$. We note first that 
$|S^p_{x_0}(x_0)|^2=B_{0,p}(x_0)$ and $P_{0,p}(x,x_0)= 
S^p_{x_0}(x)\otimes S^p_{x_0}(x_0)^*$ and therefore  
\begin{equation}\label{sz4.1} 
S^p_{x_0}(x)=\frac{1}{B_{0,p}(x_0)}P_{0,p}(x,x_0)\cdot S^p_{x_0}(x_0).
\end{equation}   
From \eqref{1c53} we deduce that for a sequence $\{r_p\}$ with $r_p\to0$ and 
$r_p\sqrt{p}\to\infty$,  
\begin{equation} \label{sz5} 
\int_{B(x_0,r_p)}|S^p_{x_0}(x)|^2\,dv_X(x)=1-\cO(1/p)\,,
\quad \text{for $p\to\infty$}. 
\end{equation} 
Relation \eqref{sz5} explains the term `peak section': when $p$ grows, the  
mass of $S^p_{x_0}$ concentrates near $x_0$.  
Since  $\Phi_p(x_p)=\Phi_p(y_p)$
we can construct as before the peak section $S^p_{x_p}=S^p_{y_p}$ as the  
unit norm generator of the orthogonal 
complement of $\Phi_p(x_p)=\Phi_p(y_p)$. We fix in the sequel such a section
which peaks at both $x_p$ and $y_p$.

We consider the distance $d(x_p,y_p)$ between the two points $x_p$ and $y_p$.
By passing to a subsequence we have two possibilities: either 
$\sqrt{p}d(x_p,y_p)\to\infty$ as $p\to\infty$
or there exists a constant $C>0$ such that $d(x_p,y_p)\leqslant C/\sqrt{p}$ for all $p$.

Assume that the first possibility is true. 
For large $p$, we learn from relation \eqref{sz5} 
that the mass of $S^p_{x_p}=S^p_{y_p}$ 
(which is $1$) concentrates both in neighborhoods $B(x_p,r_p)$ and $B(y_p,r_p)$ with
$r_p=d(x_p,y_p)/2$ and approaches  
therefore 2 if $p\to\infty$. This is a contradiction which  
rules out the first possibility.

To exclude the second possibility we follow \cite{SZ02}. We identify as usual
$B^X(x_p,\varepsilon)$ to $B^{T_{x_p}X}(0,\varepsilon)$ so the point
$y_p$ gets identified to $Z_p/\sqrt{p}$ where $Z_p\in B^{T_{x_p}X}(0,C)$. 
We define then
\begin{equation}\label{sz6}
f_p:[0,1]\longrightarrow\bR\;,\quad f_p(t)
=\frac{|S_{x_p}^p(tZ_p/\sqrt{p})|^2}{B_{0,p}(tZ_p/\sqrt{p})}\;.
\end{equation}
We have $f_p(0)=f_p(1)=1$ (again because  $S^p_{x_p}=S^p_{y_p}$) 
and $f_p(t)\leqslant1$ by the  definition of the generalized Bergman kernel.
We deduce the existence of a point $t_p\in(0,1)$ such that $f''_p(t_p)=0$. 
Equations \eqref{1c53}, \eqref{sz4.1}, \eqref{sz6} imply the estimate
\begin{equation}\label{sz7}
f_p(t)=e^{-\frac{t^2}{4}\sum_j a_j|z_{p,j}|^2}\big(1+g_p(tZ_p)/\sqrt{p}\big)
\end{equation}
and  the $\cC^2$ norm of $g_p$ over $B^{T_{x_p}X}(0,C)$ is 
uniformly bounded in $p$. 
From \eqref{sz7}, we infer that 
$|Z_p|^2_0:=\frac{1}{4}\sum_j a_j|z_{p,j}|^2=\cO(1/\sqrt{p})$. 
Using a limited expansion $e^x=1+x+x^2\varphi(x)$ for $x=t^2|Z_p|^2_0$ 
in \eqref{sz7} and taking derivatives, we obtain 
$f''_p(t)=-2|Z_p|^2_0+\cO(|Z_p|^4_0)+\cO(|Z_p|^2_0/\sqrt{p})
=(-2+\cO(1/\sqrt{p}))|Z_p|^2_0$. 
Evaluating at $t_p$ we get $0=f''_p(t_p)=(-2+\cO(1/\sqrt{p}))|Z_p|^2_0$, 
which is a contradiction since by assumption $Z_p\neq0$. 
This finishes the proof of (iii).
\end{proof} 

\begin{rem} 
Let us point out complementary results which are analogues of
\cite[(1.3)--(1.5)]{BU1} for the spaces $\mH_p$.
Computing as in \eqref{sz4}
the pull-back $\Phi_p^*h_{FS}$ of the Hermitian metric 
$h_{FS}=g_{FS}-\sqrt{-1}\,\omega_{FS}$   
on $\proj \mH^*_p$, we get the similar inequality  
to (\ref{sz1}) for $g_{FS}$ and $\om(\cdot,J\cdot)$. 
Thus, $\Phi_p$ are asymptotically
symplectic and isometric. Moreover, arguing as in \cite[Proposition\,4.4]{BU1}
we can show that $\Phi_p$ are `nearly holomorphic'\,:
\be\label{sz9}
\frac{1}{p}\,\lVert\partial\Phi_p\rVert\geqslant C\,,\quad
\frac{1}{p}\,\lVert\overline\partial\Phi_p\rVert=\cO(1/p)\,,\quad\text{for some $C>0$}\,,
\ee
uniformly on $X$, where $\Vert\,\cdot\,\Vert$ is the operator norm.
\end{rem}

\subsection{Holomorphic case revisited}\label{s5.4}

In this Section, we assume that $(X,\omega, J)$ is K\"ahler 
and the vector bundles $E,L$ are holomorphic on $X$, 
and $\nabla ^E,\nabla ^L$ are the holomorphic  Hermitian 
connections on $(E,h^E)$, $(L,h^L)$. As usual, 
$\frac{\sqrt{-1}}{2 \pi} R^L=\om$. 

But we will work with an arbitrary ({\em non-K\"ahler})
Riemannian metric $g^{TX}$ on $TX$ compatible  with $J$.
That is, in general ${\bf J}\neq J$ in \eqref{0.1}. 
The use of non-K\"ahler metrics is useful for example in Section \ref{s5.6}. 
Set 
 \begin{align} \label{f11}
\Theta(X,Y)=  g^{TX}(JX,Y).
\end{align}
Then the 2-form $\Theta$ need not to be closed 
(the convention here is different to 
\cite[(2.1)]{B} by a factor $-1$). We denote by $T^{(1,0)}X$, $T^{(0,1)}X$ 
the holomorphic and anti-holomorphic tangent bundles as in Section \ref{s3.4}.
Let $\{e_i\}$ be an orthonormal frame of $(TX, g^{TX})$.

Let $g^{TX}_{\om}(\cdot,\cdot):= \om(\cdot,J\cdot)$ be the metric on $TX$ 
induced by $\om, J$. 
We will use a subscript $\om$ to indicate the objects corresponding to 
 $g^{TX}_{\om}$, especially  $r^{X}_{\om}$ is the scalar curvature of 
$(TX,g^{TX}_{\om})$, and  $\Delta_\om$ is the Bochner-Laplace 
 operator as in \eqref{0c1} associated  to $g^{TX}_{\om}$.

Let $\overline{\partial} ^{L^p\otimes E,*}$  be the formal adjoint of
the Dolbeault operator $\overline{\partial} ^{L^p\otimes E}$ 
on the Dolbeault complex
 $\Omega ^{0,\bullet}(X, L^p\otimes E)$ with the scalar product 
induced by $g^{TX}$, $h^L$, $h^E$ as in (\ref{0c2}).
Set $D_p = \sqrt{2}( \overline{\partial} ^{L^p\otimes E}
+ \overline{\partial} ^{L^p\otimes E,*})$.
Then $D_p^2= 
2( \overline{\partial} ^{L^p\otimes E}\overline{\partial} ^{L^p\otimes E,*}
+\overline{\partial} ^{L^p\otimes E,*}\overline{\partial} ^{L^p\otimes E})$ 
preserves the $\bZ$-grading of $\Omega ^{0,\bullet}(X, L^p\otimes E)$.
Then for $p$ large enough,
\begin{equation} \label{f12}
\Ker D_p  =\Ker D_p^2  = H^0 (X,L^p\otimes E).
\end{equation}
Here $D_{p}$ is not a spin$^c$ Dirac operator on 
 $\Omega ^{0,\bullet}(X, L^p\otimes E)$, and $D^2_p$ is not 
a renormalized Bochner--Laplacian as in (\ref{laplace}).

Let $P_p(x,x')$ $(x,x'\in X)$ be the smooth kernel of the orthogonal 
projection from $\cC ^\infty(X, L^p\otimes E)$ on $\Ker D_p^2$
with respect to the Riemannian volume form $dv_X(x')$ for $p$ large enough.
Recall that we denote by $\det_{\bC} $ the determinant function on the complex
 bundle $T^{(1,0)}X$. We denote by $|{\bf J}|= (-{\bf J}^2)^{-1/2}$,
then $\det_{\bC}|{\bf J}|= (2\pi)^{-n}\Pi_i\, a_i$
 under the notation in \eqref{0ue52}.
 Now we explain how to put it in the frame of our work.

\begin{thm}\label{nonkahler}
The smooth kernel $P_p(x,x')$
has a full off--diagonal asymptotic expansion analogous to \eqref{f26}
with ${\bf j}_0=\det_\bC |{\bf J}|$ as $p\to\infty$\,.
The corresponding term $b_{0,1}$ in the expansion \eqref{0.6} of $B_{0,p}(x):=P_p(x,x)$ 
is given by
\begin{equation} \label{af12}
b_{0,1}= \frac{\det_\bC |{\bf J}|}{8\pi}\Big[r^X_\om
 -2 \Delta_\om \Big(\log({\det}_\bC |{\bf J}|)\Big)
+ 4 R^E (w_{\om,j},\ov{w}_{\om,j})\Big].
\end{equation}
here $\{w_{\om,j}\}$ is an orthonormal basis of $(T^{(1,0)}X, g^{TX}_\om)$.
\end{thm}
\begin{proof}
As pointed out in \cite[Remark 3.1]{MM}, by \cite[Theorem 1]{BiV},
 there exist $\mu_0, C_L>0$ such that for any $p\in \bN$ and
any $s\in\Omega^{>0}(X,L^p\otimes E):
=\bigoplus_{q\geqslant 1}\Omega^{0,q}(X,L^p\otimes E)$,
\begin{equation}\label{main1}
\norm{D_{p}s}^2_{L^ 2}\geqslant(2p\mu_0-C_L)\norm{s}^2_{L^ 2}.
\end{equation}
Moreover  $\spec D^2_p \subset \{0\}\cup [2p\mu_0 -C_L,+\infty[$.

Let $S^{-B}$ denote the 1-form with values in antisymmetric
 elements of $\End(TX)$ which is such that if $U,V,W \in TX$, 
\begin{align}\label{f13}
\langle S^{-B}(U)V,W  \rangle = - \frac{\sqrt{-1}}{2} 
\Big( (\partial- \overline{\partial} )\Theta\Big)(U,V,W).
\end{align}
The Bismut connection $\nabla ^{-B}$ on $TX$ is defined by 
\begin{align}\label{f14}
\nabla ^{-B} = \nabla ^{TX} +  S^{-B}.
\end{align}
Then by \cite[Prop. 2.5]{B89}, $\nabla ^{-B}$ preserves the metric $g^{TX}$
 and the complex structure of $TX$. 
Let $\nabla ^{\det}$ be the holomorphic Hermitian connection on 
$\det (T^{(1,0)}X)$ with its curvature $R^{\det}$.
Then these two connections induce naturally 
an unique connection on $\Lambda (T^{*(0,1)}X)$
which preserves its $\bZ$-grading, and with the connections 
$\nabla ^L, \nabla ^E$, we get a connection $\nabla ^{-B,E_p}$ on 
$\Lambda (T^{*(0,1)}X)\otimes L^p\otimes E$. Let $\Delta ^{-B,E_p}$ 
be the Laplacian on $\Lambda (T^{*(0,1)}X)\otimes L^p\otimes E$ 
induced by $\nabla ^{-B,E_p}$ as in (\ref{0c1}). 
For any $v\in TX$ with decomposition $v=v_{1,0}+v_{0,1}
\in T^{(1,0)}X\oplus T^{(0,1)}X$,  let ${\overline v^\ast_{1,0}}\in T^{*(0,1)}X$
 be the metric dual of $v_{1,0}$. 
 Then $c(v)=\sqrt{2}({\overline v^\ast_{1,0}}\wedge-i_{v_{\,0,1}})$ 
 defines the Clifford action of $v$ on $\Lambda (T^{*(0,1)}X)$, where $\wedge$ and $i$
denote the exterior and interior product respectively.
We define a map $^c: \Lambda (T^*X)  \to C(TX)$, the Clifford bundle
of $TX$, by sending 
$e ^{i_1}\wedge\cdots \wedge e ^{i_j}$ to $c(e_{i_1})\cdots c(e_{i_j})$
for $i_1< \cdots < i_j$. For $B\in \Lambda^3(T^*X)$, set
$|B|^2= \sum_{i<j<k}|B(e_i,e_j,e_k)|^2$.
Then we can formulate  \cite[Theorem 2.3]{B89} as following,
\begin{align}\label{f16}
D_p^2=  \Delta ^{-B,E_p} + \frac{r^X}{4} + {^c(R^E +pR^L 
+ \frac{1}{2} R^{\det})} 
+ \frac{\sqrt{-1}}{2} {^c(\overline{\partial}\partial \Theta) }
- \frac{1}{8} |( \partial- \overline{\partial} )\Theta|^2.
\end{align}

We use now the connection $\nabla ^{-B,E_p}$ instead of $\nabla ^{E_p}$ 
in \cite[\S 2]{DLM}. Then by \eqref{main1}, \eqref{f16}, everything
goes through perfectly well and as in \cite[Theorem 3.18]{DLM},
so we can directly apply the result in \cite{DLM} to get the 
{\em full off-diagonal} asymptotic expansion of the Bergman kernel. 
As the above construction preserves the $\bZ$-grading on 
$\Omega ^{0,\bullet}(X, L^p\otimes E)$, we can also directly work
on $\cC ^{\infty} (X, L^p\otimes E)$. 

Now, we need to compute the corresponding  $b_{0,1}$.
Now $h^E_\om:= ({\det}_{\bC}|{\bf J}|)^{-1} h^E$ defines a metric on $E$,
let $R^E_\om$ be the curvature associated to the holomorphic Hermitian connection of $(E, h^E_\om)$,
then 
\begin{equation}\label{af17}
R^E_\om =R^E -\ov{\partial}\partial \log ({\det}_{\bC}|{\bf J}|).
\end{equation}
Thus
\begin{equation}\label{af18}
\sqrt{-1}R^E_\om(e_{\om,j},Je_{\om,j})
= 2 R^E_\om(w_{\om,j},\ov{w}_{\om,j})
=\sqrt{-1}R^E( e_{\om,j},Je_{\om,j}) -\Delta_\om  \log ({\det}_{\bC}|{\bf J})\,.
\end{equation}
Let $\left\langle\,\cdot,\cdot \right \rangle_\om$ be the Hermitian product on 
$\cC^\infty (X, L^p\otimes E)$
induced by $g^{TX}_\om, h^L, h^E_\om$. Then 
\begin{equation}\label{af19}
(\cC^\infty (X, L^p\otimes E), \left\langle \quad \right \rangle_\om)
= (\cC^\infty (X, L^p\otimes E), \left\langle \quad \right \rangle)\,,
 \quad dv_{X,\om} =({\det}_{\bC}|{\bf J}|) dv_{X}.
\end{equation}
Observe that $H^0(X, L^p\otimes E)$ does not depend on $g^{TX}$,  $h^L$ or $h^E$.
If $P_{\om,p}(x,x')$, ($x,x'\in X$) denotes the smooth kernel of the orthogonal projection from
$(\cC^\infty (X, L^p\otimes E), \left\langle\,\cdot,\cdot \right \rangle_\om)$
onto $H^0(X, L^p\otimes E)$ with respect to $dv_{X,\om}(x)$,
we have
\begin{equation}\label{af20}
P_{p}(x,x')= ({\det}_{\bC}|{\bf J}|(x')) P_{\om,p}(x,x').
\end{equation}
Now for the kernel $ P_{\om,p}(x,x')$, we can apply Theorem \ref{t0.1} (or \cite[Theorem 1.3]{DLM})
since $g^{TX}_\om(\cdot, \cdot)= \om(\cdot, J\cdot)$ is a K\"ahler metric on $TX$,
and \eqref{af12} follows from \eqref{0.5} and \eqref{af18}.
\end{proof}

\begin{rem} Certainly, the argument in this Subsection goes through 
the orbifold case as in \cite[\S 4.2]{DLM}.
\end{rem}

\subsection{Generalizations to non-compact manifolds} \label{s5.5}
 
Let $(X,\Theta)$ be a K\"ahler manifold and $(L,h^L)$ 
be a holomorphic  Hermitian line bundle over $X$.
As in Section \ref{s5.4}, let $R^L, R^{\det}$ be the curvatures 
of the holomorphic Hermitian connections $\nabla ^L,\nabla^{\det}$ 
on $L$, $\det (T^{(1,0)}X)$, and let $J^L\in \End(TX)$ such that
$\frac{\sqrt{-1}}{2\pi}R^L(\cdot,\cdot)=\Theta (J^L\cdot,\cdot)$. The space 
of holomorphic sections of $L^p$ which are $L^2$ with respect to the norm 
given by \eqref{0c2} is denoted by $H^0_{(2)}(X,L^p)$. 
Let $P_p(x,x')$, $(x,x'\in X)$ be the Schwartz kernel of 
the orthogonal projection $P_p$ from the $L^2$ section of $L^p$
onto $H^0_{(2)}(X,L^p)$ with respect to the Riemannian volume form $dv_X(x')$
associated to $(X,\Theta)$. 
Then by the ellipticity of the Kodaira-Laplacian and Schwartz kernel theorem, 
we know $P_p(x,x')$ is $\cC^\infty$.
Choose an orthonormal basis $(S^p_i)_{i\geqslant 1}$ 
of $H^0_{(2)}(X,L^p)$. For each local holomorphic frame $e_L$ we
have $S^p_i=f^p_i e^{\otimes p}_L$ for some local holomorphic functions 
$f^p_i$. 
Then $B_p(x):=P_p(x,x)=\sum_{i\geqslant 1}|S^p_i(x)|^2
=\sum_{i\geqslant 1}|f^p_i(x)|^2|e^{\otimes p}_L|^2$ is a smooth function.
We have the following generalization of Theorem \ref{t0.1}. 
\begin{thm} \label{noncompact}
Assume that $(X,\Theta)$ is a complete K\"ahler manifold. 
Suppose that there exist $\varepsilon>0\,,\,C>0$ such that
one of the following assumptions 
holds true\,{\rm:} 
\begin{gather}
\text{ $\sqrt{-1}R^L>\varepsilon\Theta\,,\,\sqrt{-1}R^{\det}>-C\Theta$. }\label{i}\\
\text{$L=\det (T^{*(1,0)}X)$, $h^L$ is induced by $\Theta$ and 
$\sqrt{-1}R^{\det}<-\varepsilon\Theta$. }\label{ii}
\end{gather} 
 The kernel $P_p(x,x')$
has a full off--diagonal asymptotic expansion analogous to \eqref{f26}
with ${\bf j}_0=\det_\bC |J^L|$ as $p\to\infty$,
uniformly for any $x,x'\in K$, a compact set of $X$.
Especially there exist coefficients $b_r\in\cC^\infty(X)\,,\,r\in\bN$, 
such that for any compact set $K\subset X$, any $k,l\in\bN$, there exists
 $C_ {k,l,K}>0$ such that for $p\in \bN$,
\begin{equation} \label{ell1} 
\Big|\frac{1}{p^n}B_p(x)-\sum^k_{r=0}b_r(x)p^{-r}\Big|_{\cC^l(K)}\leqslant 
C_{k,l,K}\,p^{-k-1}. 
\end{equation} 
Moreover, $b_0= (\det J^L)^{1/2}$ and $b_1$ equals $b_{0,1}$ from \eqref{af12}.
\end{thm} 
\begin{proof} By the argument in Section \ref{s3.1}, if
the Kodaira--Laplacian $\Box^{L^p}=\frac{1}{2}\Delta_p:=\frac{1}{2}\Delta_{p,0}$
acting on sections of $L^p$ has a spectral gap as in \eqref{0.3},
then we can localize the problem, and we get directly (\ref{ell1}) 
from Section  \ref{s3.3}. 
Observe that $D^2_p|_{\Omega^{0,\bullet}}= \Delta_p$.
In general, on a non-compact manifold, we define a
self-adjoint extension of $D^2_p$ by 
\begin{equation*}
\begin{split}
\Dom\,D^2_p=\big\lbrace u\in\Dom\,&\db^{L^p}\cap\,\Dom\,\db^{L^p*} \,:\,
\db^{L^p}u\in\Dom\,\db^{L^p*}\,,\:\db^{L^p*} u\in \Dom\db^{L^p} \big\rbrace\,,\\
D^2_p\,u&=2(\db^{L^p}\db^{L^p*} +\db^{L^p*}\db^{L^p})\,u\,,\quad \text{for $u\in\Dom D^2_p$}\,.
\end{split}
\end{equation*}

The quadratic form associated to $D^2_p$ is the form $Q_p$ given by 
\begin{equation}\label{ell2} 
\begin{split} 
&\Dom\,Q_p:=\Dom\, \db^{L^p}\cap\;\Dom\, \db^{L^p*}\\ 
Q_p(u,v)=&2 \big\langle \db^{L^p}u\,,\db^{L^p}v\big\rangle
+ 2 \big\langle \db^{L^p*}u\,,\db^{L^p*}v
\big\rangle\,,\quad u,v\in \Dom\,Q_p\,. 
\end{split} 
\end{equation} 
In the previous formulas $\db^{L^p}$ is the weak maximal extension 
of $\db^{L^p}$ to $L^2$ forms and $\db^{L^p*}$ is its Hilbert space adjoint. 
We denote by $\Omega^{0,\bullet}_0(X,L^p)$ the space of smooth compactly 
supported forms and by $L^{0,\bullet}_2(X,L^p)$ 
the corresponding $L^2$-completion.

Under one of the hypotheses \eqref{i} or \eqref{ii} there 
exists  $\mu>0$ such that for $p$ large enough
\begin{equation}\label{ell4}
Q_p(u)\geqslant \mu p\norm{u}^2\,,\quad u\in\Dom Q_p\cap L^{0,q}_2(X,L^p)\,
 \,{\rm for} \, \, q>0.
\end{equation}
Indeed, the estimate holds for $u\in\Omega^{0,q}_0(X,L^p)$ since 
the Bochner-Kodaira formula \cite[Prop. 3.71]{BeGeV} reduces to
$Q_p(u)\geqslant 2 \big\langle\,(p R^L+ R^{\det})(w_i,\ov{w}_j)
 \ov{w}^j\wedge i_{\ov{w}_i} u\,,u\big\rangle\, $,
for $u\in\Omega^{0,q}_0(X,L^p)$, where $\{w_i\}$ is an orthonormal frame of 
$T^{(1,0)}X$.
But this implies \eqref{ell4} in general, since $\Omega^{0,\bullet}_0(X,L^p)$
is dense in $\Dom Q_p$ with respect to the graph norm, as the metric is
complete.

Next, consider $f\in\Dom\,\Delta_p\cap L^{0,0}_2(X,L^p)$ and set $u=\db^{L^p}f$.
It follows from the definition of the Laplacian and \eqref{ell4} that
\begin{equation}\label{ell5}
\lVert\Delta_p f\rVert^2= 2\big\langle\, \db^{L^p*}u\,,\db^{L^p*}u\big\rangle
=Q_p(u)\geqslant \mu p\norm{u}^2=\mu p\big\langle\,\Delta_p f\,,f\big\rangle\,.
\end{equation}
This clearly implies $\spec(\Delta_p)\subset\{0\}\cup[p\mu,\infty[$ for large $p$.
\end{proof} 

Theorem \ref{noncompact} permits an immediate generalization of Tian's
convergence theorem. Tian \cite[Theorem 4.1]{Tian} already proved 
a non--compact version for convergence in the $\cC^2$ topology 
and convergence rate $1/\sqrt{p}$\,. 
Another easy consequence are  
holomorphic Morse inequalities for the space $H^0_{(2)}(X,L^p)$.
  
Observe that the quantity $\sum_{i\geqslant 1}|f^p_i(x)|^2$ 
is not globally defined, but the current
\begin{equation}\label{ell0}
\omega_p=\frac{\sqrt{-1}}{2\pi}\,\partial\overline\partial\log
\Big(\sum_{i\geqslant 1}|f^p_i(x)|^2\Big)  
\end{equation}
is well defined globally on $X$. Indeed, since 
$R^L=-\partial\overline\partial\log|e_L|_{h^L}^2$
we have 
\begin{equation}\label{ell0.1}
\frac{1}{p}\omega_p-\frac{\sqrt{-1}}{2\pi}R^L=\frac{\sqrt{-1}}{2\pi p}\,
\partial\overline\partial\log B_p \,  .
\end{equation}
If $\dim H^0_{(2)}(X,L^p)<\infty$ we have by \eqref{sz0} that 
$\omega_p=\Phi^*_p(\omega_{FS})$ where $\Phi_p$ is defined as in 
\eqref{sz0} with $\cH_p$ replaced by $H^0_{(2)}(X,L^p)$.

\begin{cor}  Assume one of the hypotheses \eqref{i} or \eqref{ii} holds. 
Then\,{\rm:}

{\rm(a)} 
for any compact set $K\subset X$ the restriction $\omega_p|_{K}$ is a smooth
$(1,1)$-form for sufficiently large $p${\rm;} 
moreover, for any $l\in\bN$ there exists a  
constant $C_{l,K}$ such that  
\begin{equation*} 
\Big|\frac{1}{p}\omega_p-\frac{\sqrt{-1}}{2\pi} R^L\Big|_{\cC^l(K)}
\leqslant\frac{C_{l,K}}{p}\;; 
\end{equation*} 

{\rm(b)} the Morse inequalities hold in bidegree $(0,0)$\,{\rm:} 
\begin{equation} \label{ell6} 
\liminf_{p\longrightarrow\infty}p^{-n}\dim H^0_{(2)}(X,L^p)\geqslant 
\frac{1}{n!}\int_X\Big(\frac{\sqrt{-1}}{2\pi}R^L\Big)^n .
\end{equation} 
In particular, if $\dim H^0_{(2)}(X,L^p)<\infty$, the manifold $(X,\Theta)$ has  
finite volume.  
\end{cor}  
 
\begin{proof}
Due to \eqref{ell1}, $B_p$ doesn't vanish on any given compact set $K$  
for $p$ sufficiently large. Thus, (a) is a consequence of \eqref{ell1} 
and \eqref{ell0.1}.

Part (b) follows from Fatou's lemma, applied on $X$ with the measure $\Theta^n/n!$  
to the sequence $p^{-n}B_p$ which converges pointwise to 
$(\det J^L)^{1/2}=
\big(\frac{\sqrt{-1}}{2\pi}R^L\big)^n/\Theta ^n$ on $X$.   
\end{proof} 
\begin{rem} 

Under the hypothesis \eqref{ii}, the inequality \eqref{ell6} is  
\cite[Theorem 1.1]{NT} of Nadel--Tsuji, where Demailly's  
holomorphic inequalities on compact sets $K\subset X$ were used.  
The volume estimate is essential in their compactification theorem of complete K\"ahler 
manifolds with negative Ricci curvature (a generalization of the fact that arithmetic  
varieties can be complex--analytically compactified). The Morse inequalities 
\eqref{ell6} were also used by Napier--Ramachandran \cite{Na-R} to show that
some quotients of the unit ball in $\bC^n$ ($n>2$) having a strongly pseudoconvex end
have finite topological type (for the compactification of such quotients see also
\cite{MY}).    
\end{rem}

\begin{rem}\label{non-comp and non-kahler}
The statement of Theorem \ref{noncompact} still holds true for a {\em non-K\"ahler complete} metric 
$\Theta$ satisfying \eqref{i} and having bounded torsion 
$T=[i(\Theta),\partial\Theta]$, i.e. 
$\abs{T}\leqslant C$, where $\abs{T}$ is the norm with respect to
$\Theta$ (that is, Theorem \ref{nonkahler} has a non-compact version analogous
to Theorem \ref{noncompact}).
As in the proof of Theorem \ref{noncompact}, the localization argument in Section \ref{s3.1}
goes through provided we can prove the existence of the 
spectral gap of the Kodaira-Laplacian $\Box^{L^p}$\,.
This follows by applying  the generalized  Bochner-Kodaira-Nakano formula of
Demailly \cite[Theorem 0.3]{dem2} with torsion term as in \cite[Theorem 1]{BiV}.
\end{rem} 

Another generalization is a version of Theorem \ref{t0.1} 
for covering manifolds.
Let $\tx$ be a paracompact smooth manifold, such that there is a discrete group
$\g$ acting freely on $\tx$ with a compact quotient $X=\tx/\g$.
Let $\pi_{\Gamma}:\tx\longrightarrow X$ be the projection.  
Assume that there exists a $\g$--invariant pre--quantum line bundle 
$\tl$ on $\tx$  and a $\g$--invariant connection $\nabla^{\tl}$ such that 
$\tom=\frac{\sqrt{-1}}{2\pi}(\nabla^{\tl})^2$ is non--degenerate.
We endow $\tx$ with a $\g$--invariant Riemannian metric $g^{T\tx}$. 
Let $\tj$ be an $\g$-invariant almost complex structure on $T\tx$ which is separately 
compatible with $\tom$ and $g^{T \tx}$. 
Then $\tbj$, $g^{T \tx}$, $\tom$, $\tj$, $\tl$, $\te$ are the pull-back of the corresponding 
objects in Section \ref{s1} by the projection $\pi_{\Gamma}:\tx \to X$. 
Let $\Phi$ be a smooth Hermitian section of $\End(E)$,
and $\tphi=\Phi\circ \pi_{\Gamma}$.
Then  the renormalized Bochner-Laplacian $\tdel$ is
 \begin{equation*}
\tdel=\Delta^{\tl^p\otimes\te}-p\,(\tau\circ \pi_{\Gamma})+\tphi
\end{equation*}
which is an essentially self--adjoint operator. It is shown in \cite[Corollary 4.7]{MM} that
\begin{equation}\label{0ell6}
\spec\tdel\subset [-C_L,C_L]\cup[2p\mu_0-C_L,+\infty[\,,
\end{equation}
where $C_L$ is the same constant as in Section \ref{s1} and $\mu_0$ is introduced in 
\eqref{0.21}.
Let $\widetilde{\mH}_p$ be the eigenspace of $\tdel$ with the eigenvalues in
$[-C_L, C_L]$:
\be\label{0.31}
\widetilde{\mH}_p=\operatorname{Range}E\big([-C_L, C_L],\tdel\big)\,,
\ee
where $E(\,\cdot\,,\tdel)$ is the spectral measure of $\tdel$. 
From \cite[Corollary 4.7]{MM},
the von Neumann dimension of $\widetilde{\mH}_p$ equals $d_p=\dim \mH_p$\/.
Finally, we define the generalized Bergman kernel $\widetilde{P}_{q,p}$ 
of $\tdel$ as in Definition \ref{d3.0}. 
Unlike most of the objects on $\tx$,
$\widetilde{P}_{q,p}$ is not $\g$--invariant.

By \eqref{0ell6} and the proof of Proposition \ref{0t3.0}, 
the analogue of \eqref{0c7} still holds 
on any compact set $K\subset \tx$. By the finite propagation speed 
as the end of Section \ref{s3.1}, we have:
\begin{thm}\label{t0.11} We fix 
$0<\var_0 < \inf_{x\in X}\{\text{injectivity radius of $x$}\}$.  
For  any compact set $K\subset\tx$ and $k,l\in \bN$, there exists
$C_{k,\,l,\,K}>0$ such that for $x,x'\in K$, $p\in \bN$,
\begin{equation}\label{0.61}
\begin{split}
&\Big |\widetilde{P}_{q,p}(x,x') 
- P_{q,p}(\pi_{\Gamma}(x),\pi_{\Gamma}(x'))\Big |_{\cC ^l(K\times K)}
\leqslant C_{k,\,l,\,K}\: p^{-k-1}\, ,\quad  {\rm if}\,\, d(x,x')< \var_0,\\
&\Big |\widetilde{P}_{q,p}(x,x')\Big |_{\cC ^l(K\times K)}
\leqslant C_{k,\,l,\,K}\: p^{-k-1}\, ,\quad 
{\rm if}\,\,  d(x,x')\geqslant\var_0.
\end{split}
\end{equation}
Especially, $\widetilde{P}_{q,p}(x,x)$ has the same asymptotic expansion as 
$B_{q,p}(\pi_{\Gamma} (x))$ in Theorem \ref{t0.1} on any compact set $K\subset\tx$.
\end{thm}
\begin{rem}
Theorem \ref{t0.11} works well for coverings of non-compact manifolds.
Let $(X,\Theta)$ be a complete K\"ahler manifold, $(L,h^L)$ be a 
holomorphic line bundle on $X$ and let $\pi_{\Gamma}:\tx\to X$ be a Galois 
covering of $X=\tx/\g$. Let $\wi{\Theta}$ and $(\wi{L},h^{\wi{L}})$ be the 
inverse images of $\Theta$ and $(L,h^L)$ through $\pi_{\Gamma}$. 
If $(X,\Theta)$ and $(L,h^L)$ satisfy one of the 
conditions \eqref{i} or \eqref{ii}, $(\tx,\wi{\Theta})$ and 
$(\wi{L},h^{\wi{L}})$ have the same properties. We obtain therefore as 
in \eqref{ell6} (by integrating over a fundamental domain):
\be
\liminf_{p\longrightarrow\infty}p^{-n}\dim_{\g} H^0_{(2)}(\tx,\tl^p)\geqslant 
\frac{1}{n!}\int_X\Big(\frac{\sqrt{-1}}{2\pi}R^L\Big)^n .
\ee
where $\dim_{\g}$ is the von Neumann dimension of the $\g$--module 
$H^0_{(2)}(X,L^p)$. Such type of inequalities imply as in \cite{TCM:01} 
weak Lefschetz theorems \`a la Nori. 
\end{rem}

\subsection{Singular polarizations} \label{s5.6}

Let $X$ be a compact complex manifold. A {\em singular K\"ahler metric} on 
$X$ is a closed, strictly positive $(1,1)$-current $\omega$. 
This means there exist locally strictly plurisubharmonic functions 
$\varphi\in L^1_{loc}$ such that $\sqrt{-1}\partial\db\varphi=\omega$. 
 
If the cohomology class of $\omega$ in $H^2(X,\bR)$ is integral, 
there exists a holomorphic line bundle $(L,h^L)$, endowed with a  
singular Hermitian metric, such that $\frac{\sqrt{-1}}{2\pi}R^L=\omega$  
in the sense of currents. We call $(L,h^L)$ a {\em singular polarization}
of $\omega$. If we change the metric $h^L$, the curvature of the new metric 
will be in the same cohomology class as $\omega$. In this case we speak of a  
polarization of $[\omega]\in H^2(X,\bR)$.
Our purpose is to define an appropriate notion of polarized 
section of $L^p$, possibly by changing the metric of $L$, 
and study the associated Bergman kernel. 
 
First recall that a Hermitian metric $h^L$ is called {\em singular }
if it is given in local  
trivialization by functions $e^{-\varphi}$ with  
$\varphi\in L^1_{\mathrm{loc}}$.  
The curvature current $R^L$ of $h^L$ is well defined and  
given locally by the currents $\partial\db\varphi$.  
 
By the approximation theorem of Demailly \cite[Theorem 1.1]{dem},
we can assume that $h^L$ is smooth outside a proper analytic set 
$\Sigma\subset X$. Using this fundamental fact, we will introduce 
in the sequel the {\em generalized Poincar\'e metric} on $X\smallsetminus\Sigma$. 
Let $\pi:\widetilde{X}\longrightarrow X$ be a resolution of singularities such that $\pi: 
\widetilde{X}\smallsetminus\pi^{-1}(\Sigma)\longrightarrow 
X\smallsetminus \Sigma$ is biholomorphic and $\pi^{-1}(\Sigma)$ 
is a divisor with only simple normal crossings. 
Let $g^{T\wi{X}}_0$ be an arbitrary smooth $J$-invariant  metric  
on $\widetilde X$ and $\Theta '(\cdot,\cdot)=g^{T\wi{X}}_0(J\cdot,\cdot)$ 
the corresponding $(1,1)$-from. The generalized Poincar\'e metric on 
$X\smallsetminus \Sigma=\widetilde X \smallsetminus\pi^{-1}(\Sigma)$ 
is defined in \cite[\S 3]{Zu:79} by
\begin{equation}\label{poin} 
\Theta_{\var_0}=\Theta '
-\varepsilon_0\sqrt{-1}\sum_i\partial\db\log(-\log\|\sigma_i\|^2_i)^2\,, 
\quad \text{$0<\varepsilon_0\ll 1$ fixed}, 
\end{equation}   
where $\pi^{-1}(\Sigma)= \cup_i \Sigma_i$ is the decomposition into 
irreducible components $\Sigma_i$ of $\pi^{-1}(\Sigma)$
and each $\Sigma_i$ is non-singular; $\sigma_i$ are sections 
of the associated  holomorphic line bundle $[\Sigma_i]$ which vanish to
 first order on $\Sigma_i$, and $\|\sigma_i\|_i$ is the norm for a smooth 
Hermitian metric on $[\Sigma_i]$ such that $\|\sigma_i\|_i<1$. 
\begin{lemma}\label{lem-poin}
{\rm(i)} The generalized Poincar\'e metric \eqref{poin} is a complete Hermitian metric
of finite volume. It has bounded torsion and Ricci curvature.

{\rm(ii)} If $(E,h^E)$ is a holomorphic bundle over $X$ with smooth Hermitian metric
and $H_{(2)}^0(X\smallsetminus\Sigma,E)
=\big\lbrace u\in L^{0,0}_2(X\smallsetminus \Sigma, E\,,\,\Theta_{\var_0}\,,h^E):
\db^{E}u=0\big\rbrace$ then $H_{(2)}^0(X\smallsetminus\Sigma,E)=H^0(X,E)$.
\end{lemma}
\begin{proof}
To describe the metric more precisely we denote by  
$\bD$ the unit disc in $\bC$ and by $\bD^*=\bD\smallsetminus\{0\}$.  
On the product $(\bD^*)^l\times\bD^{n-l}$ we  
introduce the metric 
\begin{equation} \label{compl12,16} 
\omega_P=\frac{\sqrt{-1}}{2}\sum^l_{k=1}\frac{dz_k\wedge d\overline z_k} 
{|z_k|^2(\log|z_k|^2)^2}+\frac{\sqrt{-1}}{2}\sum^n_{k=l+1}dz_k\wedge d 
\overline z_k. 
\end{equation} 
For any point $p\in\pi^{-1}(\Sigma)$ there exists  
a coordinate neighbourhood $U$ of $p$ isomorphic to $\bD^n$ in which  
$(X\smallsetminus\pi^{-1}(\Sigma))\cap U=\{z=(z_1,\ldots,z_n):z_1\neq 0,\ldots,z_l\neq 0 
\}$. Such coordinates are called special. We endow $(X\smallsetminus\pi^{-1}(\Sigma)) 
\cap U\cong(\bD^*)^l\times\bD^{n-l}$ with the metric \eqref{compl12,16}.  
Now, a calculation in special coordinates as in \cite[Prop.\,3.4]{Zu:79}
show that the metrics \eqref{poin} and \eqref{compl12,16} are equivalent.
From this the first assertion of (i) follows.

We wish to show that there 
exist a constant $C>0$ such that 
\begin{equation}\label{zar-ell9}  
\sqrt{-1}R^{\det}>-C\Theta_{\varepsilon_0}\,,\; |T_{\varepsilon_0}|<C\,. 
\end{equation} 
where $T_{\varepsilon_0}=[\Theta_{\varepsilon_0},\partial\Theta_{\varepsilon_0}]$ 
is the torsion operator of $\Theta_{\varepsilon_0}$ and $|T_{\varepsilon_0}|$ is its norm with 
respect to $\Theta_{\varepsilon_0}$. 
Now $\partial\Theta_{\varepsilon_0}=\partial\Theta'$ by \eqref{poin}, so
it extends smoothly over $\tx$, and thus we get the second relation of  
\eqref{zar-ell9}.

We turn now to the first condition of \eqref{zar-ell9}. We have
\begin{equation}\label{zar-ell10}
\Theta_{\varepsilon_0}=\Theta '+ 2\sqrt{-1} \varepsilon_0
\sum_i \Big(\frac{R^{[\Sigma_i]}}{\log\|\sigma_i\|^2_i}
+\frac{\partial\log\|\sigma_i\|^2_i\wedge
\overline\partial\log\|\sigma_i\|^2_i}{(\log\|\sigma_i\|^2_i)^2}\Big )\, .
\end{equation}
The terms $R^{[\Sigma_i]}/\log\|\sigma_i\|^2_i$ tend to zero as we approach $\Sigma$
so they can be absorbed in $\Theta'$ and do not contribute to the singularity of 
$\Theta_{\varepsilon_0}$ near $\Sigma$\,.
To examine the last term let us localize to a point $x_0\in\Sigma$\,.
We choose special coordinates in a neighborhood $U$ of $x_0$ in which
$\Sigma_j$ has the equation $z_j=0$ for $j=1,\dots,k$ and $\Sigma_j$, $j>k$,
do not meet $U$.
Then for $1\leq i\leq k$, $\|\sigma_i\|^2_i=u_i|z_i|^2$ for some positive smooth function $u_i$
 on $U$ and
\begin{equation}\label{zar-ell11}
\frac{\partial\log\|\sigma_i\|^2_i\wedge\overline\partial\log\|\sigma_i\|^2_i}{(\log\|\sigma_i\|^2_i)^2}=
\frac{dz_i\wedge d\overline{z}_i+v_i}{|z_i|^2(\log\|\sigma_i\|^2_i)^2}
\end{equation}
where $v_i$ is a smooth $(1,1)$--\,form on $U$. Without loss of generality we may assume that
$\Theta'$ is the Euclidean metric on $U$ so that $\Theta'^n$ is the Euclidean 
volume element. Then there exists a smooth function $\beta$ such that 
\begin{equation}\label{zar-ell12}
\Theta_{\varepsilon_0}^n=
\left(1+\frac{1+\beta(z)}{\prod_i|z_i|^2(\log\|\sigma_i\|^2_i)^2}\right)\,
\Theta'^n\,=:\gamma(z)\Theta'^n\,.
\end{equation}
and consequently
\begin{multline}\label{ell13}
\sqrt{-1}\,R^{\det}=-\sqrt{-1}\,\partial\overline\partial\log\gamma(z)=
-\sqrt{-1}\,\Big(\frac{\partial\overline\partial\gamma(z)}{\gamma(z)}-
\frac{\partial\gamma(z)\wedge\overline\partial\gamma(z)}{\gamma(z)^2}\Big)\\
\geqslant-\sqrt{-1}\,\frac{\partial\overline\partial\gamma(z)}{\gamma(z)}\;.
\end{multline}
A brute force calculation of $-\sqrt{-1}\,\partial\overline\partial\gamma(z)/\gamma(z)$
and comparison to the singularities of $\Theta_{\varepsilon_0}$ given by \eqref{zar-ell11} 
show that
$\sqrt{-1}\,R^{\det}>-C\Theta_{\varepsilon_0}$ for some positive constant $C$\,.
This achieves the proof of \eqref{zar-ell9}.

Let us prove (ii). First observe that $\Theta_{\varepsilon_0}$ dominates 
the euclidian metric in special coordinates near $\pi^{-1}(\Sigma)$,
being equivalent with \eqref{compl12,16}. Therefore it dominates some positive multiple
of any smooth Hermitian metric on $\tx$. We deduce that, given a smooth Hermitian
metric $\Theta''$ on $X$, there exists a constant $c>0$ such that 
$\Theta_{\varepsilon_0}\geqslant c\Theta''$ on $X\smallsetminus\Sigma$.
It follows that elements of $H_{(2)}^0(X\smallsetminus\Sigma,E)$ are $L^2$
integrable with respect to the smooth metrics $\Theta''$ and $h^E$ over $X$, which entails
they extend holomorphically to sections of $H^0(X,E)$ by \cite[Lemme\,6.9]{De:82}.
We have therefore $H_{(2)}^0(X\smallsetminus\Sigma,E)\subset H^0(X,E)$. The reverse inclusion
follows from the finiteness of the volume of $X\smallsetminus\Sigma$ in the Poincare metric.  
\end{proof}

We can construct as in \cite[\S 4]{Ta:94} 
a singular Hermitian line bundle  
$(\widetilde{L},h^{\widetilde L})$ on $\widetilde{X}$ which is strictly  
positive and $\widetilde{L}|_{\widetilde{X}\smallsetminus\pi^{-1}(\Sigma)}\cong 
\pi^*(L^{p_0})$, for some $p_0\in\bN$.
We introduce on $L|_{X\smallsetminus\Sigma}$ the metric 
$(h^{\widetilde L})^{1/p_0}$
whose curvature extends to a strictly positive $(1,1)$--current on $\widetilde X$. Set
\begin{subequations}
\begin{align}\label{ell7}  
&h^L_{\varepsilon}=(h^{\widetilde L})^{1/p_0}\,
\prod_i(-\log\|\sigma_i\|^2_i))^\var\,, \quad 0<\varepsilon\ll 1\,,\\
&\label{ell8} 
H^0_{(2)}(X\smallsetminus \Sigma,L^p)
=\big\lbrace u\in L^{0,0}_2(X\smallsetminus \Sigma, 
L^p\,,\,\Theta_{\var_0}\,,h^L_{\varepsilon}):\db^{L^p}u=0\big\rbrace .
\end{align} 
\end{subequations}
The space $H^0_{(2)}(X\smallsetminus \Sigma,L^p)$ is the space of
 $L^2$-holomorphic sections relative to the metrics $\Theta_{\var_0}$ on 
$X\smallsetminus \Sigma$ and 
$h^L_\varepsilon$ on $L|_{X\smallsetminus \Sigma}$. 
Since $(h^{\widetilde L})^{1/p_0}$ is bounded away from zero (having
plurisubharmonic weights),
the elements of this space are $L^2$ integrable with respect to the Poincar\'e metric 
and a smooth metric $h^L_{*}$ of $L$ over whole $X$. By Lemma \ref{lem-poin} (ii)
we have $H^0_{(2)}(X\smallsetminus \Sigma,L^p)\subset H^0(X,L^p)$.
(Here we cannot infer the other inclusion since $h^L$ might be infinity on $\Sigma$.)
The space $H^0_{(2)}(X\smallsetminus \Sigma,L^p)$ is our space of polarized sections of $L^p$.  
\begin{cor}\label{moi} 
Let $(X,\omega)$ be a compact complex manifold with a singular K\"ahler 
metric with integral cohomology class. Let $(L,h^L)$ be a singular 
polarization of $[\omega]$ with strictly positive curvature current having 
singular support along a proper analytic set $\Sigma$\,.  
Then 
the Bergman kernel of the space of polarized sections \eqref{ell8} 
has the asymptotic expansion as in Theorem \ref{noncompact} 
for $X\smallsetminus\Sigma$. 
\end{cor}  
\begin{proof}
We will apply Remark \ref{non-comp and non-kahler} to the 
non--K\"ahler Hermitian manifold
$(X\smallsetminus \Sigma,\Theta_{\var_0})$ equipped with the Hermitian bundle
$(L|_{X\smallsetminus \Sigma}, h^L_\varepsilon)$. 
Thus we have to show that there
exist constants $\eta>0$, $C>0$ such that 
\begin{align}\label{ell9}  
&\sqrt{-1}R^{(L|_{X\smallsetminus \Sigma},\, 
h^L_\varepsilon)}>\eta\Theta_{\varepsilon_0}\,,
\;\sqrt{-1}R^{\det}>-C\Theta_{\varepsilon_0}\,,\; |T_{\varepsilon_0}|<C\,. 
\end{align} 
where $T_{\varepsilon_0}=[i(\Theta_{\varepsilon_0}),\partial\Theta_{\varepsilon_0}]$ 
is the torsion operator of $\Theta_{\varepsilon_0}$ and $|T_{\varepsilon_0}|$ is its norm with 
respect to $\Theta_{\varepsilon_0}$. 
The first one results for all $\varepsilon$ small enough from \eqref{poin}, \eqref{ell7} 
and the fact that the curvature of $(h^{\widetilde L})^{1/p_0}$
extends to a strictly positive $(1,1)$--\,current on $\widetilde X$
(dominating a small positive multiple of $\Theta'$ on $\widetilde X$).
The second and third relations were proved in \eqref{zar-ell9}. 
This achieves the proof of Corollary \ref{moi}.
\end{proof}

\begin{rem} 
(a) Corollary \ref{moi} gives an alternative proof of the characterization  
of Moishezon manifolds given by Ji--Shiffman \cite{JS:93}, Bonavero
\cite{Bo:93} and Takayama \cite{Ta:94}. Indeed,  
any Moishezon manifold possesses a strictly positive singular polarization  
$(L,h^L)$.    
Conversely, suppose $X$ has such a polarization. Then as in \eqref{ell6}, we  
have $\dim H^0_{(2)}(X\smallsetminus \Sigma,L^p)\geqslant C p^n$ for some $C>0$
and $p$ large enough. 
Since $H^0_{(2)} 
(X\smallsetminus  \Sigma,L^p)\subset H^0(X,L^p)$, it follows that $L$ is big and  
$X$ is Moishezon.  

(b) By \cite[Proposition 6.6. (f)]{De:96}, or \cite{Ts:92}, any big line bundle $L$ on a 
projective manifold carries a singular Hermitian metric
having  strictly positive curvature current with singularities along a proper analytic set.

(c) The results of this section hold also for reduced compact complex spaces $X$
possessing a holomorphic line bundle $L$ with singular Hermitian metric $h^L$
having positive curvature current (see \cite{Ta:94} for definitions). This is just a matter of 
desingularizing $X$. As space of polarized sections we obtain  $H^0_{(2)}(X\smallsetminus \Sigma,L^p)$
where $\Sigma$ is an analytic set containing the singular set of $X$.
\end{rem} 

\subsection*{Acknowledgments}
We express our hearty thanks to Professors Jean-Michel Bismut, Jean-Michel Bony
and Johannes Sj\"ostrand for useful conversations. It's a pleasure to 
acknowledge our intellectual debt to Xianzhe Dai and Kefeng Liu.
We are grateful Professor Paul Gauduchon for the reference \cite{Ga04}
on the Hermitian scalar curvature.
Our collaboration was partially supported by the European Commission 
through the Research Training Network ``Geometric Analysis''.


\providecommand{\href}[2]{#2}

\end{document}